\documentclass[a4paper,10pt]{elsarticle}
\usepackage{jcomp}
\usepackage{framed,multirow}
\usepackage[utf8]{inputenc}
\usepackage{bm}
\usepackage{amsmath}
\usepackage{mathrsfs}
\usepackage{graphicx}
\usepackage{epsfig, setspace}
\usepackage{url}
\usepackage{graphicx, transparent, color}
\usepackage{algorithm} 
\usepackage{algorithmicx}
\usepackage{etex}
\usepackage[bookmarks=false]{hyperref}
\usepackage{amsmath}
\usepackage{xcolor}
\usepackage{rotating} 
\usepackage{pdflscape} 
\usepackage{silence}
\usepackage{ulem,cancel}
\newtheorem{remark}{\bf Remark}[section]
\newtheorem{example}{Example}[section]
\usepackage{graphicx}
\usepackage{float}
\usepackage{subfigure}
\usepackage{subfigmat}

\usepackage{tikz}
\usepackage{capt-of}
\usepackage{setspace,lipsum}
\newcommand{\half}{\frac{1}{2}}
\usepackage{setspace,lipsum}
\begin{document}

\begin{frontmatter}

\title{Implicit gradients based novel finite volume scheme for compressible single and multi-component flows}

\author[AA_address]{Amareshwara Sainadh Chamarthi \cortext[cor1]{Corresponding author. \\ 
E-mail address: skywayman9@gmail.com (Amareshwara Sainadh  Ch.).}}
\author[AA_address]{Steven H.\ Frankel}
\author[BB_address]{Abhishek Chintagunta}
\address[AA_address]{Faculty of Mechanical Engineering, Technion - Israel Institute of Technology, Haifa, Israel}
\address[BB_address]{School of Mechanical Engineering, Vellore Institute of Technology, Vellore, India}

\begin{abstract}
This paper introduces a novel approach to compute the numerical fluxes at the cell boundaries in the finite volume approach. Explicit gradients used in deriving the reconstruction polynomials are replaced by high-order gradients computed by compact finite differences, referred to as implicit gradients in this paper. The new finite volume scheme has superior dispersion and dissipation properties in comparison to the compact reconstruction approach. These implicit gradients are re-used in viscous flux computation and post-processing, which further improves efficiency. A problem-independent shock capturing approach via Boundary Variation Diminishing (BVD) algorithm is used to suppress oscillations for the simulation of flows with shocks and material interfaces. Several numerical test cases are carried out to verify the proposed finite volume method's capability using the implicit gradient method for single and multicomponent flows. Significant improvements are observed by computing the gradients implicitly for the viscous flows.
\end{abstract}


\end{frontmatter}
\section{Introduction}\label{sec-1}

Compressible flow simulations involving discontinuities, such as shock waves and material interfaces, and turbulence is a challenging task. The choice of numerical approaches significantly affects the simulation quality, and different target applications place different demands on the numerical scheme. For compressible flows, even for smooth initial data, the flow may develop sudden changes in pressure and density mathematically treated as discontinuities. These discontinuities will result in spurious oscillations (Gibbs phenomenon), which can be mitigated by introducing additional dissipation. Even for incompressible flows, numerical schemes can fail if there is discontinuous initial data \cite{Zhang2009}.
On the other hand, the numerical scheme should have low dispersion and low dissipation to capture the small scales in amplitude and phase accurately. This numerical dissipation needed to stabilize the discontinuities will significantly impact the fine scales of turbulence. These conflicting requirements pose a significant challenge in developing numerical schemes for simulating such flows \cite{Pirozzoli2011}. These requirements are similar to the Kobayashi-Maru test, a fictional no-win scenario \cite{enwiki:1016021340}. Unlike in the fictional world where Kirk cheats, researchers in the real world try to develop meaningful solutions as nature cannot be fooled.

Of the several numerical methods proposed to solve the compressible flows, the finite volume method \cite{leveque2002finite} is more widely used due to its ability to handle arbitrary meshes and maintains local conservation. The numerical schemes proposed in the paper are based on the finite volume approach. The development of high-resolution shock-capturing schemes to solve compressible flows dates back to the pioneering work of the Godunov scheme\cite{godunov1959}, which is a piece-wise constant reconstruction. A brief review of existing high-order reconstruction schemes is given below.

The high-order accurate weighted essentially non-oscillatory (WENO) scheme introduced by Liu et al. \cite{Liu1994} and improved by Jiang and Shu \cite{jiang1995} is the most widely used approach for capturing discontinuities as well as complicated smooth features. The underlying principle of the WENO schemes is to combine lower-order interpolation schemes to get a higher-order scheme. The adaptive stencil procedure results in a non-oscillatory interpolation across discontinuities, and in the smooth regions of the solution, it will revert to its higher-order approximation. However, it was demonstrated by Henrick et al. \cite{Henrick2005} that the scheme may lose accuracy near critical points and show sub-optimal convergence in smooth regions and proposed an improvement. Borges et al. \cite{Borges2008} suggested an alternative approach with the same benefits without additional expense for the fifth-order scheme and later extended the approach to higher-order schemes \cite{castro2011high}. Johnsen et al. \cite{johnsen2010assessment} compared several high-order shock-capturing schemes. It was shown that these schemes could significantly dissipate fine scales of turbulence due to the inherent upwind-biased interpolation. Martin et al. \cite{martin2006bandwidth} optimized the WENO scheme's stencils and added an additional candidate stencil to make the background scheme symmetric rather than upwind-biased. Hu et al. \cite{Hu2010} proposed an adaptive WENO scheme that adapts between central and upwind schemes by smoothly blending the smoothness indicators of the optimal higher-order stencil and the lower order upwind stencils.

Contrary to the WENO-like combination of stencils, Fu et al. \cite{Fu2016,fu2018new} has proposed a family of Targeted ENO (TENO) schemes. TENO schemes recognized that the nonlinear adaptation is unnecessary in the smooth regions, and the reconstruction is reduced to that of linear background scheme by using a cut-off function. The cut-off function ensures that the non-smooth candidate stencils are completely discarded from the flux reconstruction. They have also optimized the linear background scheme through an optimization procedure of Hu et al. \cite{hu2012dispersion}. Such optimization results in the reduction of the formal order of accuracy but provides improved dispersion and dissipation properties. 

Compared to the computation of the single-species flows, numerical methods for multi-component flows require additional effort to prevent the appearance of numerical instabilities and oscillations in pressure at material interfaces. Johnsen and Colonius \cite{johnsen2006implementation} showed that WENO reconstruction should be carried for primitive variables to prevent oscillations and maintain pressure equilibrium. Coralic and Colonius \cite{coralic2014finite} further extended the approach to genuinely high-order and included viscous effects. He et al. \cite{he2017characteristic} used the monotonicity-preserving scheme to suppress numerical oscillations, stabilize computation, and ensure volume fraction's boundedness. Ory and Frankel \cite{haimovich2017numerical} applied the TENO reconstruction procedure, which showed improved results over the WENO.

A class of finite difference schemes presented by Lele \cite{lele1992compact} have significantly higher spectral resolution than the schemes mentioned above. It has been shown through Fourier analysis that these compact (or spatially implicit) schemes have a spectral-like resolution for short waves. Unlike the approaches mentioned earlier, compact schemes involve the computation of unknowns spanning the entire grid. Several attempts have been proposed in the literature to combine the compact schemes' advantages with those of the WENO schemes. One such earliest approach is that of Adams and Shariff \cite{adams1996high} who introduced a hybrid compact-ENO scheme where asymmetric coefficients were used for the compact scheme unlike that of \cite{lele1992compact}. This approach is further improved by Pirozzoli \cite{Pirozzoli2002} where a hybrid compact-WENO scheme is developed by replacing the non-conservative compact scheme with the conservative compact reconstruction scheme. Ren et al.\cite{Ren2003} developed a characteristic-based compact-WENO hybrid scheme where the interpolation is carried out on the characteristic variables. A characteristic-based reconstruction reduces the oscillations but incurs higher cost due to the block-tridiagonal matrix's inversion. Ghosh and Baeder \cite{ghosh2012compact} introduced a new class of upwind biased compact-reconstruction WENO schemes called CRWENO, which combined the compact upwind schemes and the WENO schemes. Despite the advantages of the WENO schemes, other shock-capturing approaches have been extensively studied in the literature. The monotonicity Preserving (MP) scheme proposed by Suresh and Hyunh \cite{suresh1997accurate} can effectively suppress numerical oscillations across the discontinuity and also preserve the accuracy. MP scheme has been optimized further and has been used for Large Eddy Simulations (LES) \cite{fang2013optimized, fang2017investigation}. Recently, Li et al.\cite{li2021low} improved the TENO schemes by filtering the non-smooth stencils by an extra nonlinear limiter such as MP limiter \cite{suresh1997accurate} or a total variation diminishing (TVD) limiters\cite{vanleer1979,arora1997well} instead of being completely discarded.

Recently a novel algorithm called Boundary Variation Diminishing (BVD) was proposed by Sun et al. \cite{sun2016boundary} which combines a non-polynomial reconstruction scheme, THINC (Tangent of Hyperbola for INterface Capturing), for discontinuous regions and an unlimited polynomial based reconstruction for the smooth regions of flows. The proposed methodology adaptively chooses the scheme with minimum Total Boundary Variation (TBV), reducing the numerical dissipation. Following \cite{sun2016boundary}, the BVD approach has been significantly improved and extended to more complex and challenging problems \cite{deng2018high,deng2020implicit,deng2019fifth,deng2020constructing}. Following their idea, Chamarthi and Frankel \cite{chamarthi2021high} have presented a new algorithm named HOCUS (High-Order Central Upwind Scheme), which combined the MP scheme and a linear-compact scheme using the BVD principle. In their approach, unlike the previous BVD implementations, both the schemes are polynomial based.

In the above mentioned finite volume methods, the reconstructed cell interface values used for the convective part of the governing equations are readily available by explicit formulas, whether they are computed using explicit or compact reconstruction approaches. Due to this, the gradients of the variables (e.g., the velocity components and temperature) required for the viscous part are re-computed. Shen et al. \cite{shen2009high,shen2010large} developed fourth- and sixth-order accurate conservative finite difference schemes for the computation of viscous terms to be used in conjunction with WENO schemes. Recently, Wang et al. \cite{wang2017compact} have proposed a variational reconstruction (VR) method for the computation of gradients for the unstructured finite volume method. Inspired by their work, Nishikawa \cite{nishikawa2018green} has proposed the Implicit Green–Gauss (IGG) method for the computation of the gradients. Both VR and IGG can be considered an extension of a compact finite difference scheme to unstructured grids. Inspired by their work, this paper proposes a new, finite volume approach for reconstructing the cell interface values using an implicit gradient method. In this paper, compact finite difference schemes of Lele \cite{lele1992compact} are used in place of the second-order gradients used in the Legendre polynomials for the derivation of the unlimited third-order reconstruction by van Leer \cite{van1977towards} (readers may also refer to Equations 3.1 and 4.1 in Ref. \cite{van1985upwind}). The proposed method has the following advantages:

\begin{description}
	\item[(a)] It has superior dissipation and dispersion properties as compared to the compact reconstruction approach of Ref. \cite{Pirozzoli2002} (please see Remark \ref{remark:3.1} and Equations \ref{eqn:upwind-compact} in the main text),
	\item[(b)] The gradients of the variables (density, pressure and velocity components) necessary for viscous fluxes are readily available with very high-order accuracy, which are reused, thereby improving the computational efficiency. These gradients are also reused for post-processing statistics like enstrophy, and
	\item[(c)] Finally, the implicit gradient approach combined with the shock-capturing approach via the BVD algorithm gave superior results of flows with shocks, material interfaces and small scale features than the approach presented in \cite{chamarthi2021high}.
\end{description}

The rest of the paper is organized as follows. Section \ref{sec-2} introduces the governing equations of single and multi-component flows. In Section \ref{sec-3}, a brief description of the finite volume method is presented and introduces the novel reconstruction schemes along with the implementation details. Numerical results and discussion are presented in Section \ref{sec-4}, and finally, in Section \ref{sec-5}, we provide concluding remarks and suggestions for future work.

\section{Governing equations}\label{sec-2}
\subsection{Single component flows}
The compressible Navier–Stokes (NS) equations for simulating single-component flows in two-dimensional (2-D) Cartesian coordinates can be expressed as:

\begin{equation}\label{CNS-base}
\frac{\partial \mathbf{Q}}{\partial t}+\frac{\partial \mathbf{F^c}}{\partial x}+\frac{\partial \mathbf{G^c}}{\partial y}+\frac{\partial \mathbf{F^v}}{\partial x}+\frac{\partial \mathbf{G^v}}{\partial y}=\mathbf{S},
\end{equation}
where  $\textbf{Q}$ are the conservative variable vector, $\mathbf{F^c}$, $\mathbf{G^c}$, and $\mathbf{F^v}$, $\mathbf{G^v}$, are convective and viscous flux vectors in each coordinate direction, respectively. In Equation (\ref{CNS-base}), the superscripts $c$ and $v$ refer to the convective and viscous terms, respectively and the compressible Euler equations are obtained by excluding the viscous fluxes. The conservative state vector and convective and viscous flux vectors vectors are expressed as:
\begin{equation}
\mathbf{Q}=\left[\begin{array}{c}
\rho \\
p u \\
\rho v \\
E
\end{array}\right],
\mathbf{F^c}=\left[\begin{array}{c}
\rho u \\
\rho u^{2}+p \\
\rho u v \\
(E+p) u
\end{array}\right],
\mathbf{G^c}=\left[\begin{array}{c}
\rho v \\
\rho u v \\
\rho v^{2}+p \\
(E+p) v
\end{array}\right],
\end{equation}

\begin{equation}\label{eqn-visc}
\mathbf{F^v}=\left[\begin{array}{c}
0 \\
-\tau_{x x} \\
-\tau_{y x} \\
-\tau_{x x} u-\tau_{x y} v+q_{x}
\end{array}\right],
\mathbf{G^v}=\left[\begin{array}{c}
0 \\
-\tau_{x y} \\
-\tau_{y y} \\
-\tau_{y x} u-\tau_{y y} v+q_{y}
\end{array}\right],
\end{equation}
here, $\rho$ is the density, whereas $u$ and $v$ are the velocity components in the $x-$ and $y-$ directions, respectively. The total energy per unit volume of the fluid is expressed as $E = \rho (e + (u^2+v^2)/2)$, where $e$ is the specific internal energy. The system of equations is closed with the ideal gas equation of state which relates the thermodynamic pressure $p$ and the total energy:

\begin{equation}\label{eqn:pressure}
p = (\gamma -1) (E - \rho \frac{(u^2+v^2)}{2}),
\end{equation}
where $\gamma$ is the ratio of specific heats of the fluid. The components of the viscous stress tensor $\tau$ and the heat flux $q$ are given by

\begin{equation}\label{eqn:5-stress}
\tau_{x x}=\frac{2}{3} \mu\left(2 \frac{\partial u}{\partial x}-\frac{\partial v}{\partial y}\right), \quad \tau_{x y}=\tau_{y x}=\mu\left(\frac{\partial u}{\partial y}+\frac{\partial v}{\partial x}\right), \quad \tau_{y y}=\frac{2}{3} \mu\left(2 \frac{\partial v}{\partial y}-\frac{\partial u}{\partial x}\right)
\end{equation}

\begin{equation}\label{eqn:6-heat}
\begin{aligned}
q_{x} &=-k \frac{\partial T}{\partial x}=-\frac{\gamma \mu}{\operatorname{Pr}(\gamma-1)} \frac{\partial}{\partial x}\left(\frac{p}{\rho}\right), \\
q_{y} &=-k \frac{\partial T}{\partial y}=-\frac{\gamma \mu}{\operatorname{Pr}(\gamma-1)} \frac{\partial}{\partial y}\left(\frac{p}{\rho}\right)
\end{aligned}
\end{equation}
in which $T$  is the temperature, $Pr$ is the Prandtl number, and $\mu$ is the dynamic viscosity.

\subsection{Multicomponent flows}

In this study, the compressible multi-component flows described by the quasi-conservative five equation model of Allaire et al. \cite{allaire2002five} including viscous effects \cite{coralic2014finite} are considered. For a system consisting of two fluids, it consists of two continuity equations, momentum and an energy equation, respectively, and an advection equation for the volume fraction of one of the two fluids and are expressed as:

\begin{equation}\label{5eqn-base}
\frac{\partial \mathbf{Q}}{\partial t}+\frac{\partial \mathbf{F^c}}{\partial x}+\frac{\partial \mathbf{G^c}}{\partial y}+\frac{\partial \mathbf{F^v}}{\partial x}+\frac{\partial \mathbf{G^v}}{\partial y}=\mathbf{S},
\end{equation}
where the state vector, flux vectors and the source term, $\mathbf{S}$, are

\begin{equation}
\mathbf{Q}=\left[\begin{array}{c}
\alpha_{1} \rho_{1} \\
\alpha_{2} \rho_{2} \\
\rho u \\
\rho v \\
E \\
\alpha_{1}
\end{array}\right], \quad \mathbf{F^c}=\left[\begin{array}{c}
\alpha_{1} \rho_{1} u \\
\alpha_{2} \rho_{2} u \\
\rho u^{2}+p \\
\rho v u \\
(E+p) u \\
\alpha_{1} u
\end{array}\right], \quad \mathbf{G^c}=\left[\begin{array}{c}
\alpha_{1} \rho_{1} v \\
\alpha_{2} \rho_{2} v \\
\rho u v \\
\rho v^{2}+p \\
(E+p) v \\
\alpha_{1} v
\end{array}\right],  \quad \mathbf{S}=\left[\begin{array}{c}
0 \\
0 \\
0 \\
0 \\
0 \\
\alpha_{1} \nabla \cdot \mathbf{u}
\end{array}\right],
\end{equation}

\begin{equation}
\mathbf{F^v}=\left[\begin{array}{c}
0 \\
-\tau_{x x} \\
-\tau_{y x} \\
-\tau_{x x} u-\tau_{x y} v
\end{array}\right],
\mathbf{G^v}=\left[\begin{array}{c}
0 \\
-\tau_{x y} \\
-\tau_{y y} \\
-\tau_{y x} u-\tau_{y y} v
\end{array}\right],
\end{equation}
where $\rho_1$ and $\rho_2$ corresponds to the densities of fluids 1 and 2, $\alpha_{1}$ and $\alpha_{2}$ are the volume fractions of the fluids 1 and 2, $\rho$, $u$,$v$, $p$ and $E$ are the density, x- and y- velocity components, pressure, total energy per unit volume of the mixture respectively. The five-equation model is incomplete in the vicinity of the material interface where the fluids are in a mixed state while using a diffuse interface approach for mathematical modelling. In order to complete the closure a set of mixture rules for various properties of the fluids are to be defined. The mixture rules for the volume fractions of the two fluids $\alpha_1$ and $\alpha_2$, density and mixture rule for the ratio of specific heats $\gamma$ of the mixture are given by:

\begin{equation}
\alpha_{2}=1-\alpha_{1}
\end{equation}

\begin{equation}
\rho=\rho_{1} \alpha_{1}+\rho_{2} \alpha_{2}
\end{equation}

\begin{equation}
\frac{1}{\gamma-1}=\frac{\alpha_{1}}{\gamma_{1}-1}+\frac{\alpha_{2}}{{\gamma}_{2}-1}
\end{equation}
where $\gamma_1$ and $\gamma_2$ are the specific heat ratios of the fluids 1 and 2.  Finally, Under the isobaric assumption, the equation of state given by Equation (\ref{eqn:pressure}) is used to close the system. The viscous terms are defined as in the single-component flows. The advection equation for the volume fraction is written in non-conservative form, given by Equation (\ref{eqn-source}), following Johnsen and Colonius \cite{johnsen2006implementation}. This choice ensures the advection equation is consistently coupled with the governing equations, leading to an oscillation-free behaviour at the material interface.

\begin{equation}\label{eqn-source}
\frac{\partial \alpha_{1}}{\partial t}+\nabla \cdot\left(\alpha_{1} \mathbf{u}\right)=\alpha_{1} \nabla \cdot \mathbf{u}
\end{equation}

\section{Numerical methods}\label{sec-3}

As explained in the introduction, the study of high-resolution numerical schemes for compressible flow problems will be proposed in the context of the finite volume approach in this section. The equations described in Section \ref{sec-2} are therefore discretized using the finite volume method on a uniform ($\Delta x= \Delta y= h$) Cartesian grid. The resulting semi-discrete form of the equations is given for cell $(j, i)$ by:

\begin{equation}\label{eqn-differencing}
\begin{aligned}
\frac{d \mathbf{\hat Q}_{j, i}}{d t}=&-\frac{1}{\Delta x}\left[\left(\mathbf {\hat{F}^c}_{j+ 1 / 2}-\mathbf {\hat{F}^c}_{j- 1 / 2}\right)-\left(\mathbf {\hat{F}^v}_{j+ 1 / 2}-\mathbf {\hat{F}^v}_{j- 1 / 2}\right)\right] \\
&-\frac{1}{\Delta x}\left[\left(\mathbf {\hat{G}^c}_{i+ 1 / 2}-\mathbf {\hat{G}^c}_{i- 1 / 2}\right)-\left(\mathbf {\hat{G}^v}_{i+ 1 / 2}-\mathbf {\hat{G}^v}_{i- 1 / 2}\right)\right] +\mathbf{S}_{j, i} \\
=& \mathbf{Res}\left(\mathbf{\hat Q}_{j, i}\right)
\end{aligned}
\end{equation}
where $\mathbf{\hat Q}$ is the vector of state variables, $\mathbf {\hat{F}^c}$, $\mathbf {\hat{G}^c}$ and $\mathbf {\hat{F}^v}$, $\mathbf {\hat{G}^v}$ are regarded as numerical approximation of the convective and viscous fluxes in the x, and y directions respectively, $\mathbf{S}_{j, i}$ is the source term, and $\mathbf{Res}(\mathbf{\hat Q}_{j, i})$ is the residual function. In the following subsections, we provide the details of these computations. In subsection \ref{sec-3.1} the computation of convective fluxes is explained which includes the novel implicit gradient method (\ref{sec-3.1}), shock-capturing algorithm (\ref{sec-3.1.1}), details of Riemann solver (\ref{sec-3.1.2}), source term discretization for multicomponent flows are described in \ref{sec-3.1.3}. In subsection \ref{sec-3.2} the computation of viscous fluxes is explained and finally temporal integration is presented in subsection \ref{sec-3.3}.


%
%
%

\subsection{Spatial discretization of convective fluxes}\label{sec-3.1}
\begin{figure}
\centering
 \includegraphics[width=0.9\textwidth]{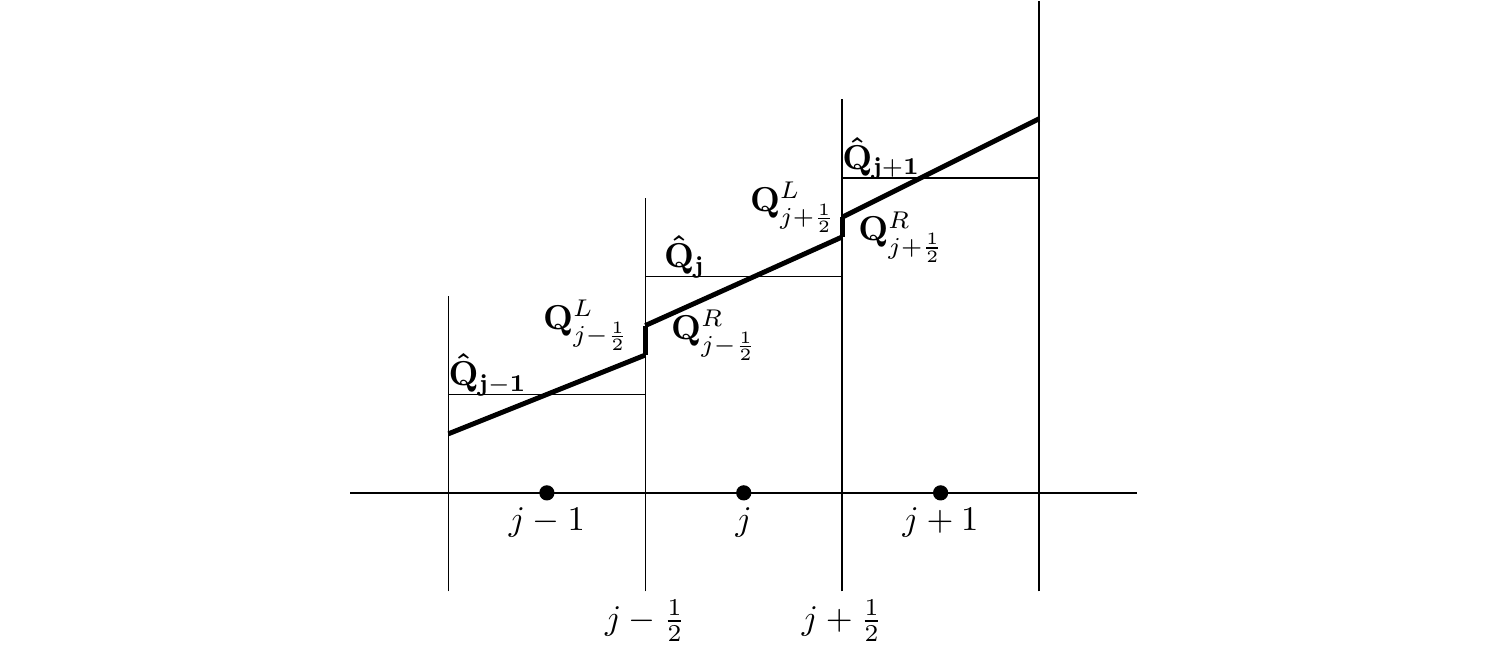}
\caption{Piecewise linear representation within cells}
\label{recon}
\end{figure}

In this section, we present the computation of convective fluxes. As it can be seen from Equation (\ref{eqn-differencing}), if the numerical approximation of flux at the cell interface has high order approximation, then one can achieve high order approximations to the corresponding to $\mathbf{\hat Q}$. These numerical fluxes, $\mathbf {\hat{F^c}}_{j+ 1 / 2}$, are computed by a Riemann solver. There are several types of Riemann solvers in literature  \cite{deledicque2007exact, ivings1998riemann,roe1981approximate, batten1997choice, einfeldt1988godunov, toro1994restoration, osher1982upwind}, and the canonical form of Riemann flux can be written as:

\begin{equation}\label{Numerical Flux}
\mathbf {\hat{F}^c}_{j+ 1 / 2}={F}^{\rm Riemann}_{j+\frac{1}{2}}(\mathbf{Q}_{j+\frac{1}{2}}^{L},\mathbf{Q}_{j+\frac{1}{2}}^{R}), 
\end{equation}
\begin{equation}
\mathbf{F}^{\rm Riemann}_{j+\frac{1}{2}} = \frac{1}{2}({\mathbf{F}^L_{j+\frac{1}{2}}} + {\mathbf{F}^R}_{j+\frac{1}{2}}) -
 \frac{1}{2} | {\mathbf{A}_{j+\frac{1}{2}}}|({\mathbf{Q}^R_{j+\frac{1}{2}}}-{\mathbf{Q}^L_{j+\frac{1}{2}}}),
\label{eqn:Riemann}
\end{equation}
where $L$ and $R$ are adjacent values of a cell interface, as show in Fig. \ref{recon}, and $|{\mathbf{A}_{j+\frac{1}{2}}}|$ denotes the characteristic signal velocity evaluated at the cell interface in a hyperbolic equation or in the case of Euler equations, the inviscid Jacobian. The procedure of obtaining the values at the interface from cell-averaged cell center variables is called reconstruction procedure. It is obvious from the Equation (\ref{eqn:Riemann}) a core problem is how to reconstruct the left- and right-side values, $\mathbf{Q}_{j+\frac{1}{2}}^{L}$ and $\mathbf{Q}_{j+\frac{1}{2}}^{R}$, for cell boundaries, which can fundamentally influence the numerical solution. The values ${\mathbf{Q}^L}_{j+\frac{1}{2}}$ should be biased to left and similarly ${\mathbf{Q}^R}_{j+\frac{1}{2}}$ should be biased right for upwinding. Representing these numerical approximations of $L$ and $R$ at cell interface as a piecewise constant is equivalent to first-order approximation, i.e. 
\begin{equation}\label{eqn:first-order}
\begin{aligned}
{\mathbf{Q}^L}_{j+\frac{1}{2}} &= \mathbf{\hat{Q}_j}\\
{\mathbf{Q}^R}_{j+\frac{1}{2}} &= \mathbf{\hat{Q}_{j+1}}. 
\end{aligned}
\end{equation}

Therefore, a linear approximation of the solution, shown in Fig. \ref{recon}, is a second-order spatial approximation, while a quadratic representation on each cell leads to a third-order spatial approximation. By considering a general local representation, as explained in \cite{van1977towards}, these approximation can be expressed in terms of Legendre polynomials, valid for $x_{j-1 / 2} \leq x \leq x_{j+1 / 2}$:

\begin{equation}\label{eqn:legendre}
\mathbf{Q}(x)=\mathbf{\hat{Q}}_{j}+\frac{\mathbf{Q}'_j}{\Delta x_j}\left(x-x_{j}\right)+\frac{3  \mathbf{Q}''_j}{2 \Delta x_j^{2}} \kappa\left[\left(x-x_{j}\right)^{2}-\frac{\Delta x_j^{2}}{12}\right] 
\end{equation}
where $\mathbf{\hat{Q}}_{j}$ is the cell-averaged value and ${ \mathbf{Q}'_j}$, ${\mathbf{Q}''_j}$ are the estimations of the first and second derivatives within the cell $j$. The Equation (\ref{eqn:legendre}) is the basis for the Monotonic Upstream-centered Scheme for Conservation Laws (MUSCL) scheme or popularly known as kappa scheme of van Leer \cite{van1977towards, van1985upwind}. The WENO schemes are an extension of the MUSCL scheme to an arbitrary order of accuracy (see derivations in Ref. \cite{balsara2016efficient}). For the numerical approximations of the Riemann problem we need the values at the cell interfaces only. By setting $x = x_j \pm \Delta x/2$ within a cell $j$ gives us the interface values:

\begin{equation}\label{eqn:left_right}
\begin{aligned}
\mathbf{Q}_{j+1 / 2}^{L} &=\mathbf{\hat {Q}}_{j}+\frac{1}{2} \mathbf{Q}'_j+\frac{\kappa}{4} \mathbf{Q}''_j \\
\\
\mathbf{Q}_{j-1 / 2}^{R} &=\mathbf{\hat {Q}}_{j}-\frac{1}{2}  \mathbf{Q}'_j+\frac{\kappa}{4}  \mathbf{Q}''_j \quad \text{or} \ \mathbf{Q}_{j+1 / 2}^{R} &=\mathbf{\hat {Q}}_{j+1}-\frac{1}{2} \mathbf{Q}'_{j+1}+\frac{\kappa}{4} \mathbf{Q}''_{j+1}
\end{aligned}
\end{equation}

In order to define these approximation completely, the derivatives ${ \mathbf{Q}'_j}$ and ${\mathbf{Q}''_j}$ have to be estimated. By using $\kappa=\frac{1}{3}$, which gives third order accuracy for cell-averaged solutions, and substituting the following explicit central differences for the derivatives in Equation (\ref{eqn:left_right}),

\begin{equation}\label{eqn:expligrad}
\begin{aligned}
 \mathbf{Q}'_j &= \frac{\mathbf{\hat {Q}}_{j+1} - \mathbf{\hat {Q}}_{j-1}}{2}\\
 \mathbf{Q}''_j &= \mathbf{\hat {Q}}_{j+1} - 2\mathbf{\hat {Q}}_{j} + \mathbf{\hat {Q}}_{j-1},
\end{aligned}
\end{equation}
we obtain the following third order reconstruction polynomials,

\begin{equation}
\begin{aligned} 
{\mathbf{Q}^L}_{j+\frac{1}{2}}&= \frac{1}{6}\left(-\mathbf{\hat {Q}}_{j-1} + 5\mathbf{\hat {Q}}_{j} + 2\mathbf{\hat {Q}}_{j+1} \right) \\
{\mathbf{Q}^R}_{j-\frac{1}{2}} &= \frac{1}{6}\left(2\mathbf{\hat {Q}}_{j-1} + 5\mathbf{\hat {Q}}_{j} - \mathbf{\hat {Q}}_{j+1} \right). \
\end{aligned}
\label{eqn:3linear}
\end{equation}

In this paper the derivatives, ${ \mathbf{Q}'}$ and ${ \mathbf{Q}''}$, in Equation (\ref{eqn:left_right}) are computed by using high-order compact finite difference schemes proposed in \cite{lele1992compact} and therefore the reconstruction procedure is named as \textit{implicit gradient} (IG) finite volume approach.  Lele \cite{lele1992compact} proposed two kinds of central compact schemes. One is a linear cell-centered compact scheme, staggered grid, and the other is cell-node, collocated grid, (see Ref. \cite{lele1992compact,nagarajan2003robust,boersma2005staggered} for discussion on these two variants) compact scheme. For the current approach cell-node compact scheme is appropriate for the evaluation of the derivatives, which can be written in a general form for the first derivatives as:

\begin{equation} \label{eqn:pade}\small
        \beta \mathbf{ {Q}}'_{j-2} + \alpha \mathbf{ {Q}}'_{j-1} + \mathbf{ {Q}}'_{j} + \alpha \mathbf{ {Q}}'_{j+1} + \beta \mathbf{ {Q}}'_{j+2} =
        c \frac{\mathbf{\hat {Q}}_{j+3} - \mathbf{\hat {Q}}_{j-3}}{6 \Delta x} + b \frac{\mathbf{\hat {Q}}_{j+2} - \mathbf{\hat {Q}}_{j-2}}{4 \Delta x} + a \frac{\mathbf{ \hat{Q}}_{j+1} - \mathbf{\hat {Q}}_{j-1}}{2 \Delta x}
\end{equation}

The left hand side of the Equation (\ref{eqn:pade}) contain the spatial derivatives $\mathbf{{Q}}'_j$ while the right hand side contains the function values $\mathbf{\hat {Q}}$ at the cell center $x_j$, respectively. Compact finite difference schemes of different orders of accuracy are derived by matching the Taylor series coefficients with different constraints on the parameters $\alpha$, $\beta$, $a$, $b$ and $c$ and are listed in \cite{lele1992compact}. In this work, we considered difference schemes for the first derivatives with the following parameters,

\begin{subequations}
     \begin{alignat}{4}
&\beta = 0, \quad &a_1=\frac{2}{3}(\alpha+2), \quad &b_1=\frac{1}{3}(4 \alpha-1), \quad &c=0 \label{eq:param4} \\
&\beta = \frac{-1+3 \alpha}{12}, \quad &a=\frac{2}{9}(8 - 3 \alpha), \quad &b=\frac{1}{18}(-17 +57 \alpha), \quad &c=0 \label{eq:param6}.
     \end{alignat}
   \end{subequations}
By substituting $\alpha$ = $\frac{5}{14}$ in Equation (\ref{eq:param4}) we obtain the optimised fourth-order compact derivative (see Fig. 2  in Ref. \cite{lele1992compact}) and is denoted by \textbf{CD4}. For $\alpha$ = $\frac{1}{3}$ in Equation (\ref{eq:param6})we obtain the sixth-order compact derivative and is denoted by  \textbf{CD6}, and can be written as:

\begin{subequations}\label{eqn:ddx}
     \begin{alignat}{1}
     \frac{5}{14} \mathbf{ {Q}}_{j-1}^{\prime}+\mathbf{ {Q}}_{j}^{\prime}+\frac{5}{14} \mathbf{ {Q}}_{j+1}^{\prime}&=\frac{b_1}{4 \Delta x}\left(\mathbf{\hat {Q}}_{j+2}-\mathbf{\hat {Q}}_{j-2}\right)+\frac{a_1}{2 \Delta x}\left(\mathbf{\hat {Q}}_{j+1}-\mathbf{\hat {Q}}_{j-1}\right) \label{eq:cd4} \\
\frac{1}{3} \mathbf{ {Q}}_{j-1}^{\prime}+\mathbf{ {Q}}_{j}^{\prime}+\frac{1}{3} \mathbf{ {Q}}_{j+1}^{\prime}&=\frac{1}{36 \Delta x}\left(\mathbf{\hat {Q}}_{j+2}-\mathbf{\hat {Q}}_{j-2}\right)+\frac{7}{9 \Delta x}\left(\mathbf{\hat {Q}}_{j+1}-\mathbf{\hat {Q}}_{j-1}\right) \label{eq:cd6} ,
     \end{alignat}
   \end{subequations}
where $j= 1, 2 , 3, .....,N-1$. Unlike the second-order central differences given by the Equations (\ref{eqn:expligrad}), which depend only on values at $j-1$, $j$ and $j+1$, compact finite differences depend on all the nodal values of the domain and therefore mimic the global dependence of the spectral methods. This global dependence results in a tridiagonal system of equations that the Thomas algorithm can easily invert. Near the boundary cells, lower-order one-sided difference formulas are used to approximate derivatives $\mathbf{{Q}_0}'$ and $\mathbf{{Q}_N}'$. The following third-order formulas are considered for both the CD4 and CD6 schemes in the present work

\begin{equation}
\mathbf{{Q}'_0}+2\mathbf{{Q}'_1} = \frac{1}{\Delta x}(\frac{-5}{2}\mathbf{\hat {Q}}_{0}+2\mathbf{\hat {Q}}_1+\frac{1}{2}\mathbf{\hat {Q}}_{2})
\end{equation}  

\begin{equation}
\mathbf{{Q}'_N}+2\mathbf{{Q}'_{N-1}} = \frac{1}{\Delta x}(\frac{5}{2}\mathbf{\hat {Q}}_{N}+2\mathbf{\hat {Q}}_{N-1}+\frac{1}{2}\mathbf{\hat {Q}}_{N-2})
\end{equation}  


For the computation of the second derivatives (Hessians), ${\mathbf{Q}''_j}$, we compute the \textit{derivative of the first derivatives} obtained from the Equations (\ref{eq:cd4}) and (\ref{eq:cd6}) as shown below:
\begin{subequations}\label{eqn:ddx2}
     \begin{alignat}{1}
     \frac{5}{14} \mathbf{ {Q}''}_{j-1}+\mathbf{ {Q}''}_{j}+\frac{5}{14} \mathbf{ {Q}''}_{j+1}&=\frac{b_1}{4 \Delta x}\left(\mathbf{ {Q}'}_{j+2}-\mathbf{ {Q}'}_{j-2}\right)+\frac{a_1}{2 \Delta x}\left(\mathbf{ {Q}'}_{j+1}-\mathbf{ {Q}'}_{j-1}\right) \label{eq:cd42} \\
\frac{1}{3} \mathbf{ {Q}''}_{j-1}+\mathbf{ {Q}''}_{j}+\frac{1}{3} \mathbf{ {Q}''}_{j+1}&=\frac{1}{36 \Delta x}\left(\mathbf{ {Q}'}_{j+2}-\mathbf{{Q}'}_{j-2}\right)+\frac{7}{9 \Delta x}\left(\mathbf{ {Q}'}_{j+1}-\mathbf{ {Q}'}_{j-1}\right) \label{eq:cd62} ,
     \end{alignat}
   \end{subequations}

After obtaining the derivatives ${ \mathbf{Q}'}$ and ${ \mathbf{Q}''}$ from solving Equations (\ref{eqn:ddx}) and (\ref{eqn:ddx2})  they are substituted into the Equations (\ref{eqn:left_right}) to obtained the left and right reconstructed values necessary for the Riemann problem, completing the implicit gradient approach.  The left and right interface values obtained from the implicit gradient approach are still \textit{upwind} in nature and $\kappa$ =1/3 for both the schemes. The schemes given by Equations (\ref{eqn:IG}) are denoted by \textbf{IG4}, and \textbf{IG6} and are obtained by evaluating the gradients and Hessians according to CD4, and CD6, respectively.

\begin{equation}\label{eqn:IG}
\small
\begin{cases} 
 \mathbf{Q}_{j+1 / 2}^{L,IG} &=\mathbf{\hat {Q}}_{j}+\frac{1}{2} \mathbf{Q}'_j+\frac{1}{12} \mathbf{Q}''_j \\
\\											
\mathbf{Q}_{j-1 / 2}^{R,IG} &=\mathbf{\hat {Q}}_{j}-\frac{1}{2}  \mathbf{Q}'_j+\frac{1}{12}  \mathbf{Q}''_j 
\end{cases}
\begin{aligned}
\rightarrow &\text{$\textbf{Q}'_j$ and $\textbf{Q}''_j$ computed}\\& \text{by Eqns (\ref{eq:cd4}) and (\ref{eq:cd42}) is \textbf{IG4} and Eqns (\ref{eq:cd6}) and (\ref{eq:cd62}) is \textbf{IG6}}.
\end{aligned}
\end{equation}

\begin{remark}\label{remark:3.1}
\normalfont It is important to note that the left and right interface values obtained from the \textit{implicit gradient} approach explained above are different from the compact reconstruction polynomials (denoted as C5 in this paper and given by Equations (\ref{eqn:upwind-compact}) ) discussed by Pirozzoli  \cite{Pirozzoli2002} despite their \textit{implicit in space} approach required to compute them. The reconstruction approach by C5 will only give us the $L$ and $R$ interface values. However, the gradients of velocity components and temperature necessary for the viscous fluxes (see Equations (\ref{eqn:6-heat}) and (\ref{eqn:5-stress})) are to be computed separately. The significant advantage of IG schemes over the C5 approach is that these gradients are readily available to be used in viscous flux computation (explained in Section \ref{sec-3.2}), which makes them efficient.
\end{remark}

\begin{subequations}
     \begin{alignat}{1}
\frac{1}{2}  \mathbf {{Q}}^{L, C5}_{j-\frac{1}{2}}+  \mathbf {{Q}}^{L, C5}_{j+\frac{1}{2}} + \frac{1}{6}  \mathbf {{Q}}^{L, C5}_{j+\frac{3}{2}}&= \frac{1}{18}\mathbf{\hat {Q}}_{j-1}+  \frac{19}{18}\mathbf{\hat {Q}}_{j} + \frac{5}{9}\mathbf{\hat {Q}}_{j+1} \label{eq:left}\\
\frac{1}{6}  \mathbf {{Q}}^{R, C5}_{j-\frac{1}{2}}+ \mathbf {{Q}}^{R, C5}_{j+\frac{1}{2}} + \frac{1}{2}  \mathbf {{Q}}^{R, C5}_{j+\frac{3}{2}}&= \frac{5}{9}\mathbf{\hat {Q}}_{j}+ \frac{19}{18}\mathbf{\hat {Q}}_{j+1} + \frac{1}{18}\mathbf{\hat {Q}}_{j+2} \label{eq:CR}
     \end{alignat}
     \label{eqn:upwind-compact}
   \end{subequations}
   
\begin{figure}[H]
\centering
\subfigure[Dispersion]{\includegraphics[width=0.45\textwidth]{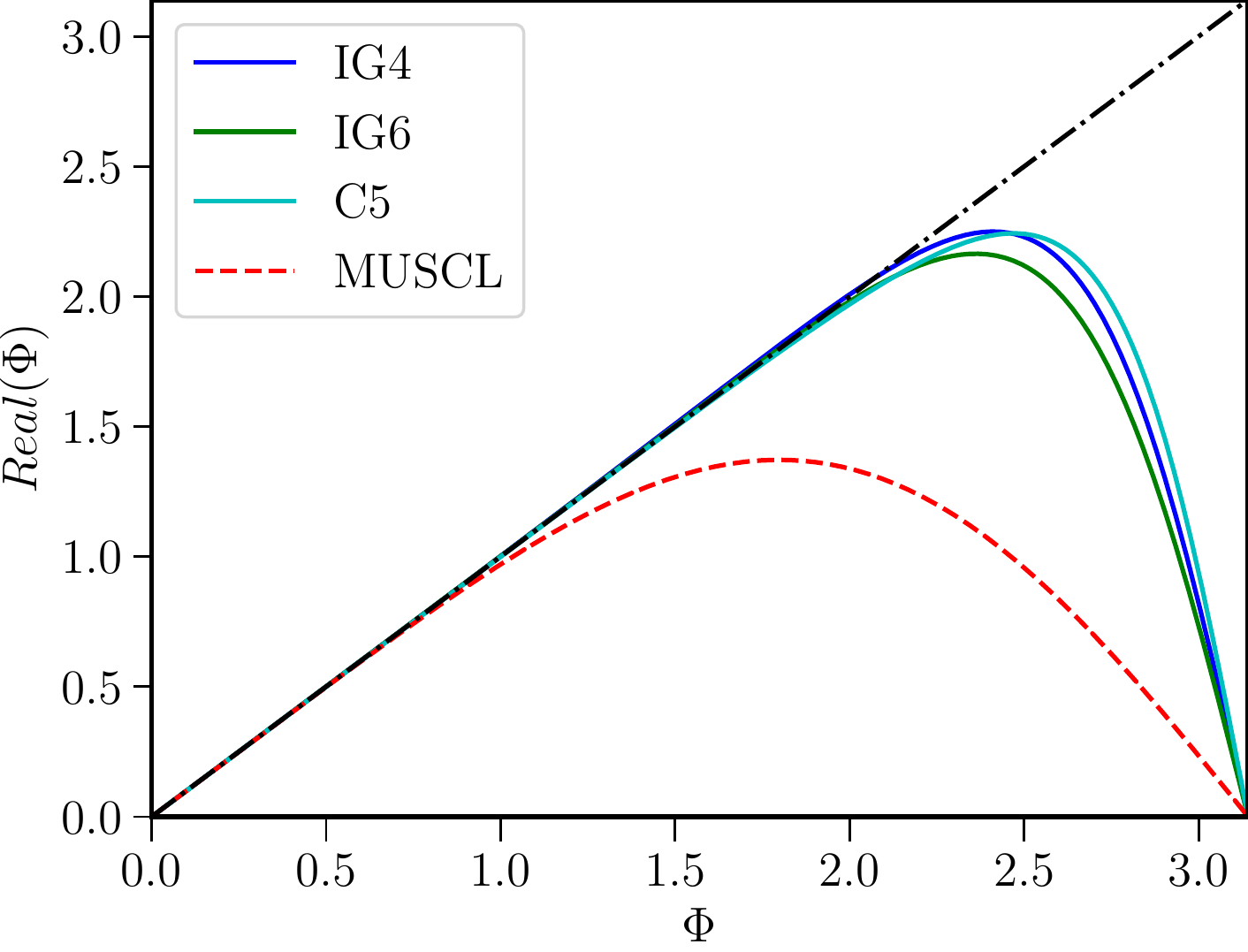}
\label{fig:dispersion}}
\subfigure[Dissipation]{\includegraphics[width=0.45\textwidth]{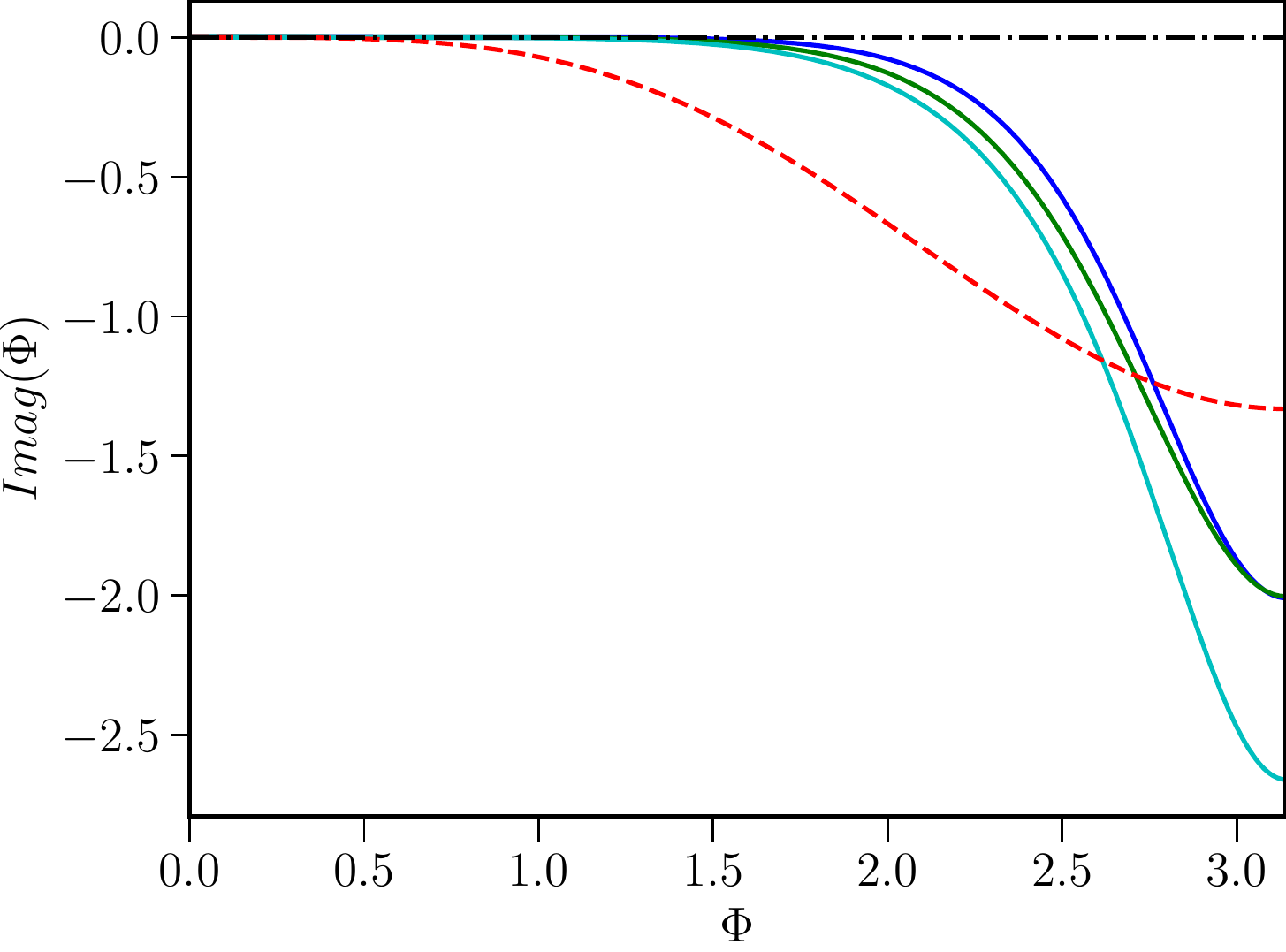}
\label{fig:dissipation}}
\caption{Dispersion and Dissipation properties of the linear upwind schemes.}
\label{fig_disp}
\end{figure}

\begin{remark} \label{remark-3.2}
\normalfont We analysed the dispersion and dissipation properties of the proposed schemes in this paper using Fourier analysis. Fig. \ref{fig_disp} shows the dispersion and dissipation properties of the IG4, IG6, MUSCL and C5 schemes. Dispersion properties of the CD4 are better than CD6 and are reflected in IG4 and IG6 schemes (see Fig. 2  in Ref. \cite{lele1992compact} regarding CD4 and CD6). It can be seen that the proposed schemes have superior dissipation in comparison with the C5 scheme. The dispersion property of the IG4 scheme is very close to that of the C5 scheme. These properties are also confirmed later by the inviscid Taylor-Green simulations (Example \ref{ex:TGV} in numerical results). Comparison is also made with the MUSCL scheme, Equation (\ref{eqn:3linear}) which uses explicit gradients. It can be seen that explicit gradients are much more dispersive and dissipative than the implicit gradients considered here.
\end{remark}

\begin{remark}\label{eqn:whichgradeints}
\normalfont  Both conservative, $\mathbf{Q}$, and primitive variables, $\mathbf{U}$, can be used for the evaluation of the gradients in the implicit gradient approach for the single component flows. For shock-capturing purposes (will be explained in Section \ref{sec-3.1.1}) and ease of reusing the gradients in viscous fluxes (explained momentarily in Section \ref{sec-3.2}), \textit{primitive variables} are used for reconstruction in this paper for both single and multi-component flows.  For the single component flows primitive variable vector is, $\mathbf{U}$ = $(\rho, u, v, p)^T$ and for multi-component flows it is, $\mathbf{U}$ = $(\rho_{1} \alpha_{1},\rho_{2} \alpha_{2}, u, v, p,\alpha_1)^T$, respectively. For brevity, only the gradients in the x direction are shown below,

\begin{gather}
 \textbf{for single component flows}  \rightarrow \frac{\partial \rho}{\partial x}, \frac{\partial u}{\partial x}, \frac{\partial v}{\partial x}, \frac{\partial p}{\partial x}  \nonumber\\
\nonumber \\
 \textbf{for multi-component flows}  \rightarrow \frac{\partial \rho_{1} \alpha_{1}}{\partial x}, \frac{\partial \rho_{2} \alpha_{2}}{\partial x}, \frac{\partial u}{\partial x}, \frac{\partial v}{\partial x}, \frac{\partial p}{\partial x}, \frac{\partial  \alpha_{1}}{\partial x}  \nonumber.
\end{gather}
\end{remark}


\subsubsection{Shock-capturing via BVD algorithm}\label{sec-3.1.1}
The novel schemes, IG4 and IG6, derived in the earlier section, are linear in nature and therefore lead to oscillations for flows involving material interfaces and shocks. In this section, we describe the shock-capturing scheme using the BVD algorithm, previously presented in \cite{chamarthi2021high}, extended to the IG schemes for both single and multi-component flows. By comparing two different polynomials, the BVD algorithm selects the reconstruction polynomial with minimum numerical dissipation by evaluating the Total Boundary Variation (TBV) given by the Equation (\ref{Eq:TBV}) for each cell for each \textit{primitive variable}:

\begin{equation}\label{Eq:TBV}
\mathbf{TBV}_{j}=\big|\mathbf{U}_{j-\frac{1}{2}}^{L}-\mathbf{U}_{j-\frac{1}{2}}^{R}\big|+\big|\mathbf{U}_{j+\frac{1}{2}}^{L}-\mathbf{U}_{j+\frac{1}{2}}^{R} \big|
\end{equation} 

For a given cell, the terms on the right-hand side of Equation (\ref{Eq:TBV}) represent the amount of numerical dissipation introduced in the numerical flux in Equation (\ref{eqn:Riemann}). The BVD algorithm compares the TBVs of the concerned polynomials and selects the one that is least dissipative at an interface. In the smooth regions of the flow, the IG6 or IG4 linear schemes will be used, and in the presence of discontinuities, the BVD algorithm will turn to the MP5 scheme. The combination of IG4 and MP5 is denoted as \textbf{IG4MP} and, the combination of IG6 and MP5 is denoted as \textbf{IG6MP} in this paper, respectively. Both schemes are denoted together as \textbf{IGMP} schemes. 

\begin{remark}
\onehalfspacing
\normalfont The choice of MP5 scheme as a candidate polynomial for shock-capturing has been arrived at by testing several different combinations of linear and nonlinear schemes. It is also possible to use a third-order reconstruction with the minmod limiter and THINC, and readers are referred to Appendix A and B of Ref.\cite{chamarthi2021high}.
\end{remark}

Physically, the density and pressure should be positive, and failure to preserve these variables' positivity may cause a blow-up of the numerical solutions. The reconstructed states, $L$ and $R$, obtained by IG4MP, IG6MP and even MP5 methods, may sometimes lead to unphysical values for density and pressure and also lead to unboundedness of volume fractions in multi-component flow simulations. In the present paper, we use the order reduction strategy of Titarev and Toro \cite{titarev2004finite}. We note that in the test cases we have addressed, the IG4MP and IG6MP schemes failed in very few cells, even after using the BVD algorithm, and in relatively few time steps, this choice is justified.  If the MP5 reconstruction also produced unphysical values in the worst-case scenario, the reconstructed states are reduced to first-order (FO) accurate reconstructions given by Equations (\ref{eqn:first-order}). In the following, reconstruction procedure only in the $x$-direction is discussed. Due to the dimension-by-dimension approach, the other directions are handled the same way and the procedure is summarized below:

\begin{description}
\item[Step i.] Evaluate the interface values by using two different reconstruction procedures:
\begin{enumerate}[(a)]
\item Linear implicit gradient reconstruction given by Equation (\ref{eqn:IG}), either IG4 or IG6 and
\item MP5 scheme \cite{suresh1997accurate} given by the following equations Equations (\ref{eqn:mp5}). For brevity, we only explain the procedure for the left interface values, $\mathbf{U}^{L,MP5}_{j+\half}$, since the right interface values, $\mathbf{U}^{R,MP5}_{j+\half}$ can be obtained via symmetry. The steps involved are as follows.

\begin{equation}\label{eqn:mp5}
\mathbf{U}^{L,MP5}_{j+1/2}=\left\{\begin{array}{ll}
\mathbf{U}_{j+1 / 2}^{\text {Linear }} & \text { if }\left(\mathbf{U}_{j+1 / 2}^{\text {Linear }}-\mathbf{\hat{U}}_{j}\right)\left(\mathbf{U}_{j+1 / 2}^{\text {Linear }}-\mathbf{U}_{j+1 / 2}^{M P}\right) \leq 10^{-20} \\
\mathbf{U}_{j+1 / 2}^{\text {Nonlinear }} & \text { otherwise }
\end{array}\right.
\end{equation}
where
\begin{equation}\label{eqn:alpha}
\begin{aligned}
\mathbf{U}_{j+1 / 2}^{\text {Linear }} &=\frac{1}{60}(2 \mathbf{\hat{U}}_{j-2} - 13 \mathbf{\hat{U}}_{j-1} + 47 \mathbf{\hat{U}}_{j} + 27 \mathbf{\hat{U}}_{j+1} - 3 \mathbf{\hat{U}}_{j+2}),  \\
\mathbf{U}_{j+1 / 2}^{\text {Nonlinear }} &=\mathbf{U}_{j+1 / 2}^{\text {Linear }}+\operatorname{minmod}\left(\mathbf{U}_{j+1 / 2}^{\min }-\mathbf{U}_{j+1 / 2}^{\text {Linear }}, \mathbf{U}_{j+1 / 2}^{\max }-\mathbf{U}_{j+1 / 2}^{\text {Linear }}\right), \\
\mathbf{U}_{j+1 / 2}^{M P} &=\mathbf{U}_{j+1 / 2}^{\text {Linear }}+\operatorname{minmod}\left[\mathbf{\hat{U}}_{j+1}-\mathbf{\hat{U}}_{j}, \mathbf{\alpha}\left(\mathbf{\hat{U}}_{j}-\mathbf{\hat{U}}_{j-1}\right)\right] \\
\mathbf{U}_{j+1 / 2}^{\min } &=\max \left[\min \left(\mathbf{\hat{U}}_{j}, \mathbf{\hat{U}}_{j+1}, \mathbf{U}_{j+1 / 2}^{M D}\right), \min \left(\mathbf{\hat{U}}_{j}, \mathbf{U}_{j+1 / 2}^{U L}, \mathbf{U}_{j+1 / 2}^{L C}\right)\right] \\
\mathbf{U}_{j+1 / 2}^{\max } &=\min \left[\max \left(\mathbf{\hat{U}}_{j}, \mathbf{\hat{U}}_{j+1}, \mathbf{U}_{j+1 / 2}^{M D}\right), \max \left(\mathbf{\hat{U}}_{j}, \mathbf{U}_{j+1 / 2}^{U L}, \mathbf{U}_{j+1 / 2}^{L C}\right)\right] \\
\mathbf{U}_{j+1 / 2}^{M D} &=\frac{1}{2}\left(\mathbf{\hat{U}}_{j}+\mathbf{\hat{U}}_{j+1}\right)-\frac{1}{2} d_{j+1 / 2}^{M} \\
\mathbf{U}_{j+1 / 2}^{U L} &=\mathbf{\hat{U}}_{j}+4\left(\mathbf{\hat{U}}_{j}-\mathbf{\hat{U}}_{j-1}\right) \\
\mathbf{U}_{j+1 / 2}^{L C} &=\frac{1}{2}\left(3 \mathbf{\hat{U}}_{j}-\mathbf{\hat{U}}_{j-1}\right)+\frac{4}{3} d_{j-1 / 2}^{M} \\
d_{j+1 / 2}^{M} &=\operatorname{minmod}\left(4 d_{j}-d_{j+1}, 4 d_{j+1}-d, d_{j},d_{j+1}\right), \\
d_{j} &=\mathbf{\hat{U}}_{j-1}-2 \mathbf{\hat{U}}_{j}+\mathbf{\hat{U}}_{j+1}
\end{aligned}
\end{equation}

where,
\begin{equation}
minmod(a,b) = \half \left ( sign(a)+sign(b) \right ) min(|a|,|b|), 
\end{equation}

\noindent Similar to the Ref. \cite{chamarthi2021high}, for both single- and multi-component flows, the above mentioned MP5 reconstruction involves the transforming of primitive variables into characteristic variables, $\mathbf{W}$, as explained below:

\begin{enumerate}
\item Compute the Roe averages at the interface $(x_{j+\frac{1}{2}})$ by using neighbouring cells, $(x_j)$ and $(x_{j+1})$. Compute the left and right eigenvectors $\bm{L_{n}}$ and $\bm{R_{n}}$, ( \textcolor{black}see Appendix \ref{sec-appa}). The primitive variables are then transformed to characteristic variables by the following equation: \\
\begin{equation}
	\bm{W}_{j} = \bm{L}_{\bm{n}_{j+\frac{1}{2},i}} \bm{U}_{j}.
\end{equation}
\item Carry out the reconstruction of characteristic variables by MP5 scheme, through Equations ( \textcolor{black}{\ref{eqn:mp5}}), and obtain left- and right- interface values denoted by $\tilde{\bm{W}}_{j+\frac{1}{2},i}^L$ and $\tilde{\bm{W}}_{j+\frac{1}{2},i}^R$. \\
\item After obtaining $\tilde{\bm{W}}_{j+\frac{1}{2},i}^L$ and $\tilde{\bm{W}}_{j+\frac{1}{2},i}^R$ from the MP5 reconstruction the primitive variables are then recovered by projecting the characteristic variables back to physical fields: \\
\begin{equation}
\begin{aligned}
	{\bm{U}}_{j+\frac{1}{2},i}^{L,MP5} &= \bm{R}_{\bm{n}_{j+\frac{1}{2},i}} \tilde{\bm{W}}_{j+\frac{1}{2},i}^L, \\
  {\bm{U}}_{j+\frac{1}{2},i}^{R,MP5} &= \bm{R}_{\bm{n}_{j+\frac{1}{2},i}} \tilde{\bm{W}}_{j+\frac{1}{2},i}^R.
\end{aligned}
\end{equation}

\end{enumerate}

\end{enumerate}
\item[Step ii.] Calculate the TBV values for each cell $I_{j}$ by using the implicit gradient reconstruction for IG4,
\begin{equation}\label{Eq:TBVIG4}
\mathbf{TBV}_{j}^{IG4}=\big|\mathbf{U}_{j-\frac{1}{2}}^{L,IG4}-\mathbf{U}_{j-\frac{1}{2}}^{R,IG4}\big|+\big|\mathbf{U}_{j+\frac{1}{2}}^{L,IG4}-\mathbf{U}_{j+\frac{1}{2}}^{R,IG4} \big|
\end{equation} 
for IG6,
\begin{equation}\label{Eq:TBVIG6}
\mathbf{TBV}_{j}^{IG6}=\big|\mathbf{U}_{j-\frac{1}{2}}^{L,IG6}-\mathbf{U}_{j-\frac{1}{2}}^{R,IG6}\big|+\big|\mathbf{U}_{j+\frac{1}{2}}^{L,IG6}-\mathbf{U}_{j+\frac{1}{2}}^{R,IG6} \big|
\end{equation} 

and  MP5 scheme:
\begin{equation}\label{Eq:TBVSC}
\mathbf{TBV}_{j}^{MP5}=\big|\mathbf{U}_{j-\frac{1}{2}}^{L,MP5}-\mathbf{U}_{j-\frac{1}{2}}^{R,MP5}\big|+\big|\mathbf{U}_{j+\frac{1}{2}}^{L,MP5}-\mathbf{U}_{j+\frac{1}{2}}^{R,MP5} \big|.
\end{equation}

\item[Step iii.] Now, \textcolor{black}{all the interface values at $j-\frac{3}{2}$, $j-\frac{1}{2}$, $j+\frac{1}{2}$ and $j+\frac{3}{2}$, both $L$ and $R$} are modified according to the following algorithm, for each primitive variable, to obtain non-oscillatory results. For the IG4MP scheme we have

\begin{equation}
\text{if} \; \mathbf{TBV}^{MP5} < \mathbf{TBV}^{IG4} \;  \left\{\begin{matrix}
\mathbf{U}^{K,IG4}_{j-\frac{3}{2}} = \mathbf{U}^{{K},MP5}_{j-\frac{3}{2}} \\ 
\mathbf{U}^{K,IG4}_{j-\frac{1}{2}} = \mathbf{U}^{{K},MP5}_{j-\frac{1}{2}} \\ 
\mathbf{U}^{K,IG4}_{j+\frac{1}{2}} = \mathbf{U}^{{K},MP5}_{j-\frac{1}{2}} \\ 
\mathbf{U}^{K,IG4}_{j+\frac{3}{2}} = \mathbf{U}^{{K},MP5}_{j+\frac{3}{2}}
\end{matrix}\right.
\end{equation}

and for the IG6MP scheme we obtain,

\begin{equation}
\text{if} \; \mathbf{TBV}^{MP5} < \mathbf{TBV}^{IG6} \;  \left\{\begin{matrix}
\mathbf{U}^{K,IG6}_{j-\frac{3}{2}} = \mathbf{U}^{{K},MP5}_{j-\frac{3}{2}} \\ 
\mathbf{U}^{K,IG6}_{j-\frac{1}{2}} = \mathbf{U}^{{K},MP5}_{j-\frac{1}{2}} \\ 
\mathbf{U}^{K,IG6}_{j+\frac{1}{2}} = \mathbf{U}^{{K},MP5}_{j-\frac{1}{2}} \\ 
\mathbf{U}^{K,IG6}_{j+\frac{3}{2}} = \mathbf{U}^{{K},MP5}_{j+\frac{3}{2}}
\end{matrix}\right.
\end{equation}

where $K$ = $L$ or $R$.

\item[Step iv.] For each reconstructed state, $K$, we check the following conditions for positivity of pressure and density for single-component flows

\begin{equation}\label{eqn:single-positivity}
\begin{aligned}
\rho^{K} & \leq 0 \quad \text{\textbf{or}} \quad p^{K} \leq 0
\end{aligned}
\end{equation}

and for multi-component flows following conditions, which includes additional conditions for phasic densities and boundedness of volume fractions, are checked:
\begin{equation}\label{eqn:multi-positivity}
\begin{aligned}
\left(\rho_{1} \alpha_{1}\right)^{K}  \leq 0 \ \text{\textbf{or}} \quad \left(\rho_{2} \alpha_{2}\right)^{K}  \leq 0 \ \text{\textbf{or}} \quad \left(\rho_{1} \alpha_{1}\right)^{K}+\left(\rho_{2} \alpha_{2}\right)^{K} \leq 0 \ \text{\textbf{or}} \\ \quad \left(\alpha_{1}\right)^{K}  <0 \ \text{\textbf{or}} \quad \left(\alpha_{1}\right)^{K} >1  \ \text{\textbf{or}} \quad p^{K} \leq 0
\end{aligned}
\end{equation}
 If conditions given by Equations (\ref{eqn:single-positivity}) and (\ref{eqn:multi-positivity}) are not satisfied we decrease the order of reconstruction and the procedure is as follows,
\begin{equation}
\mathbf{U}^K=\left\{\begin{array}{ll}
\mathbf{U}^{K, \text { IG4MP}} \text { or } \mathbf{U}^{K, \text { IG6MP}}& \text { default } \\
\mathbf{U}^{K, \text { MP5}}& \text { if } \mathbf{U}^{K, \text { IG4MP}} \text { or } \mathbf{U}^{K, \text { IG6MP}} \textbf { fails } \\
\mathbf{U}^{K, \text { FO}} & \text { if } \mathbf{U}^{K, \text { MP5}} \textbf { fails } 
\end{array}\right.
\end{equation}
where, $\mathbf{U}^{K, \text { FO}}$ is the first order approximation computed using Equation. (\ref{eqn:first-order}).

\item[Step v.] Evaluate the conservative variables, $(\mathbf{Q}_{j+\frac{1}{2}}^{L},\mathbf{Q}_{j+\frac{1}{2}}^{R})$, from the primitive variables, $(\mathbf{U}_{j+\frac{1}{2}}^{L},\mathbf{U}_{j+\frac{1}{2}}^{R})$ obtained from the above procedure, and compute the interface flux $\mathbf{F}^{\rm Riemann}_{j+\frac{1}{2}}$
\end{description}

\subsubsection{Riemann solver}\label{sec-3.1.2}

Riemann solvers approximate the convective flux after the obtaining reconstructed states at the interface as explained in the earlier section. This section illustrates how the HLLC Riemann solver \cite{batten1997choice,toro1994restoration} approximates convective fluxes. For simplicity, only the HLLC approximations for both single and multicomponent flows in the x-direction are illustrated in this section. The HLLC flux in $x$-direction is given by:

\begin{equation}\label{eq:HLLC}
\mathbf{F}^{\rm Riemann}= \mathbf{F}^{HLLC}=\left\{\begin{array}{ll}
\mathbf{F}_{L} & , \text { if } \quad 0 \leq S_{L} \\
\mathbf{F}_{* L} & , \text { if } \quad S_{L} \leq 0 \leq S_{*} \\
\mathbf{F}_{* R} & , \text { if } S_{*} \leq 0 \leq S_{R} \\
\mathbf{F}_{R} & , \text { if } \quad 0 \geq S_{R}
\end{array}\right.
\end{equation}
\begin{equation}
\mathbf{F}_{* K}=\mathbf{F}_{K}+S_{K}\left(\mathbf{Q}_{* K}-\mathbf{Q}_{K}\right)
\end{equation}
where $L$ and $R$ are the left and right states respectively. With $K$ = $L$ or $R$, the star state for single-species flow is defined as:

\begin{equation}\label{eqn-hllc-single}
\mathbf{Q}_{* K}=\left(\frac{S_{\mathrm{K}}-u_{\mathrm{K}}}{S_{\mathrm{K}}-S_{*}}\right)\left[\begin{array}{c}
\rho_{K} \\
\rho_{\mathrm{K}} S_{*} \\
\rho_{K} v_{K} \\
E_{k}+\left(S_{*}-u_{K}\right)\left(\rho_{K} S_{*}+\frac{p_{K}}{S_{K}-u_{K}}\right)
\end{array}\right]
\end{equation}

For two-component flow with five-equation model is defined as:

\begin{equation}\label{eqn-hllc-multi}
\mathbf{Q}_{* \mathrm{K}}=\left(\frac{S_{\mathrm{K}}-u_{\mathrm{K}}}{S_{\mathrm{K}}-S_{*}}\right)\left(\begin{array}{c}
\left(\alpha_{1} \rho_{1}\right)_{\mathrm{K}} \\
\left(\alpha_{2} \rho_{2}\right)_{\mathrm{K}} \\
\rho_{\mathrm{K}} S_{*} \\
\rho_{\mathrm{K}} v_{\mathrm{K}} \\
E_{\mathrm{K}}+\left(S_{*}-u_{\mathrm{K}}\right)\left(\rho_{\mathrm{K}} S_{*}+\frac{p_{\mathrm{K}}}{S_{\mathrm{K}}-u_{\mathrm{K}}}\right) \\
\alpha_{1_{\mathrm{K}}}
\end{array}\right)
\end{equation}
 In the above expressions, the waves speeds  $S_L$ and $S_R$ can be obtained as suggested by Einfeldt \cite{einfeldt1988godunov}, $ S_L = min(u_{L}-c_L, \tilde{u}-\tilde{c}) \ \text{and} \ S_R = max(u_{R}+c_R, \tilde{u} +\tilde{c})$, where $\tilde {u}$ and $\tilde {c}$ are the Roe averages from the left and right states. Batten et al. \cite{batten1997choice} provided a closed form expression for $S_*$ which is as follows:

\begin{align}
 \label{eqn:HLLCmiddlewaveestimate}
 S_* = \frac{p_R - p_L + \rho_Lu_{nL}(S_L - u_{nL}) - \rho_Ru_{nR}(S_R - u_{nR})}{\rho_L(S_L - u_{nL}) -\rho_R(S_R - u_{nR})}.
\end{align}

\subsubsection{Source term for multi-component flows}\label{sec-3.1.3}
The HLLC solver is also used for the computation of the non-conservative equation in the multi-component flows to ensure consistency with the other conservative equations. The approximated velocity components at any cell interface for the approximation of source term in \textcolor{black}{Equation (\ref{eqn-source})} is as follows, for instance, the x component of the velocity is computed as:
\begin{equation}
\begin{aligned}
{{u}_{j+\half}} &=\frac{1+\operatorname{sgn}\left(S_{*}\right)}{2}\left[u_{L}+s_{-}\left(\frac{S_{L}-u_{L}}{S_{L}-S_{*}}-1\right)\right] +\frac{1-\operatorname{sgn}\left(S_{*}\right)}{2}\left[u_{R}+s_{+}\left(\frac{S_{R}-u_{R}}{S_{R}-S_{*}}-1\right)\right]
\end{aligned}
\end{equation}
where
\begin{equation}
s_{-}=\min \left(0, s_{L}\right), \quad s_{+}=\max \left(0, s_{R}\right)
\end{equation}
Similarly, the $y-$ component of the velocity can also be computed. Finally, the source term can be evaluated within a computational cell $j,i$ using,
\begin{equation}\label{eqn-sourceterm}
\begin{aligned}
\left(\alpha_{1} \nabla \cdot \mathbf{u}\right)_{j, i} &=\left(\alpha_{1}\right)_{j, i}\left[ \frac{1}{\Delta x}\left(u_{j+\frac{1}{2}, i}-u_{j-\frac{1}{2}, i}\right)+\frac{1}{\Delta y}\left(v_{j, i+\frac{1}{2}}-v_{j, i-\frac{1}{2}}\right)\right]
\end{aligned}
\end{equation}

\textcolor{black}{\subsection{Spatial discretization of viscous fluxes}\label{sec-3.2}}

In this section, we describe the discretization of the viscous fluxes. The majority of the research focus of compressible flows is on the convective fluxes to resolve discontinuities, but it is essential to evaluate viscous fluxes with high order accuracy. Shen et al. \cite{shen2009high, shen2010large} developed a set of fourth and sixth order conservative viscous flux formulations for compressible flows. Shen et al. \cite{shen2009high} (see their Table 3) computed the derivatives of the primitive variables at each cell center by the following fourth order finite differences:
\begin{equation}\label{eqn:4e}
\left(\frac{\partial \mathbf{U}}{\partial x}\right)_{i}=\frac{8}{12 \Delta x}\left(\mathbf{U}_{i+1}-\mathbf{U}_{i-1}\right)-\frac{1}{12 \Delta x}\left(\mathbf{U}_{i+2}-\mathbf{U}_{i-2}\right)
\end{equation}
The discretization of viscous fluxes in this paper reuses the gradients that are already computed for convective fluxes, only x-direction fluxes are presented, and is as follows:

\begin{description}
\item[Step i.] In this paper, the gradients of primitive variables are \textit{reused as they are already computed} by the Equations (\ref{eqn:ddx}) as explained in Remark \ref{eqn:whichgradeints}. Moreover, the temperature gradients are also computed from the already existing pressure and density gradients, as shown below

\begin{equation}
\frac{\partial T}{\partial x} = \frac{\partial}{\partial x}\left(\frac{p}{\rho}\right)=\frac{1}{\rho} \frac{\partial p}{\partial x}-\frac{p}{\rho^{2}} \frac{\partial \rho}{\partial x}
\end{equation}

Using the same gradients for both convective and viscous fluxes is the \textit{main advantage} of proposed numerical schemes that improve the accuracy of the compressible flow simulations' and overall efficiency.  

\item[Step ii.] Once the gradients of the necessary variables are obtained it is straightforward to compute the fluxes, $\mathbf{\tilde{F}}_v$, at each cell center according to the Equations (\ref{eqn-visc}) using the relations (\ref{eqn:5-stress}) and (\ref{eqn:6-heat}).
\item[Step iii.] The fluxes, $\mathbf{\tilde{F}}_v$, are then interpolated to the cell interfaces by using the sixth-order interpolation at the internal points, and a fourth-order interpolation at the boundary points given by Equations \eqref{sixth} and \eqref{fourth}, respectively, to obtain the numerical viscous fluxes, $\mathbf {\hat{F}^v}_{j+ 1 / 2}$: 

\begin{equation} \label{sixth}
\begin{split}
\mathbf {\hat{F}^v}_{j+ 1 / 2} = \frac{1}{60} & \left(\frac{\partial \mathbf{\tilde{F}}_v}{\partial x}\right )_{j-2} - \frac{8}{60} \left(\frac{\partial \mathbf{\tilde{F}}_v}{\partial x}\right )_{j-1} +  \frac{37}{60} \left(\frac{\partial \mathbf{\tilde{F}}_v}{\partial x}\right )_{j} + \\
&\frac{37}{60} \left(\frac{\partial \mathbf{\tilde F}_v}{\partial x}\right )_{j+1} - \frac{8}{60} \left(\frac{\partial \mathbf{\tilde{F}}_v}{\partial x}\right )_{j+2} + \frac{1}{60} \left(\frac{\partial \mathbf{\tilde F}_v}{\partial x}\right )_{j+3} 
\end{split}
\end{equation}

\begin{equation} \label{fourth}
\small
\mathbf {\hat{F}^v}_{j+ 1 / 2} = \frac{1}{12} \left[ -\left(\frac{\partial \mathbf{ \tilde F}_v}{\partial x}\right )_{j+2} + 7 \left(\frac{\partial \mathbf{\tilde F}}{\partial x}\right )_{j+1} + 7\left(\frac{\partial \mathbf{\tilde F}_v}{\partial x}\right )_{j} -\left(\frac{\partial \mathbf{\tilde F}_v}{\partial x}\right )_{j-1}   \right ]
\end{equation}

\end{description}

The numerical method described above can be easily extended to multi-dimensional (2D and 3D) problems via dimension by dimension approach. The complete numerical algorithm for the compressible viscous flows is summarized below:
\begin{algorithm}\label{algo:IGMP} 
\onehalfspacing
\caption{Complete algorithm for the compressible viscous flows.}
	\begin{algorithmic}[1]
        \State{Compute the first and second derivatives of all the primitive variables in both $x-$ and $y-$ directions as outlined in Section \ref{sec-3.1}. Equations (\ref{eqn:ddx}) and (\ref{eqn:ddx2}).}
        \State{To compute the contribution of convective fluxes in the Equation (\ref{eqn-differencing}) in the $x-$ direction:}
        \begin{itemize} 
        \item Evaluate the reconstructed states outlined in Section \ref{sec-3.1.1} through \textbf{Steps i - v}.
        \item Using the reconstructed states the numerical fluxes are obtained by Riemann solver outlined in \ref{sec-3.1.2}. Equations (\ref{eqn-hllc-single}) for single-, and Equations (\ref{eqn-hllc-multi}) for multi-component flows, respectively.
        \item For multi-component flow evaluate the contribution of the source term by Equation (\ref{eqn-sourceterm}).
        \end{itemize}
        \State{To compute the contribution of  viscous fluxes in the Equation (\ref{eqn-differencing}) in the $x-$ direction, \textbf{Steps i - iii} in Section \ref{sec-3.1.2} are followed.}
        \State{To evaluate the contribution of the fluxes in $y-$ direction repeat \textbf{steps (2)-(3)} in this algorithm.}
        \State{Evaluate the residual in Equation (\ref{eqn-differencing}) and perform time integration as described in the next section.}
	\end{algorithmic} 
\end{algorithm}

{\subsection{Temporal integration}\label{sec-3.3}
Finally, the semi-discrete approximation of the equations of motion is temporally integrated. The conserved variables are  integrated in time using the following third-order TVD Runge–Kutta scheme \cite{jiang1995}:

\begin{eqnarray} \label{rk}
\mathbf{\hat Q^{(1)}} & = & \mathbf{\hat Q^n} + \Delta t \mathbf{Res}(\mathbf{\hat Q^n}), \nonumber \\
\mathbf{\hat Q^{(2)}} & = & \frac{3}{4} \mathbf{\hat Q^n} + \frac{1}{4} \mathbf{\hat Q^{(1)}} + \frac{1}{4}\Delta t \mathbf{Res}(\mathbf{\hat Q^{(1)}}) ,\\
\mathbf{\hat Q^{n+1}} & = & \frac{1}{3} \mathbf{\hat Q^n} + \frac{2}{3} \mathbf{\hat Q^{(2)}} + \frac{2}{3}\Delta t \mathbf{Res}(\mathbf{\hat Q^{(2)}}), \nonumber
\end{eqnarray}
where $\mathbf{Res}$ is the residual that is evaluated on the right-hand side of Equation(\ref{eqn-differencing}). The superscripts $\mathbf{n}$ and $\mathbf{n+1}$ denote the current and the subsequent time-steps, and superscripts $\mathbf{(1)-(2)}$ corresponds to intermediate steps, respectively. The time step $\Delta t$ is taken as suggested by Titarev and Toro \cite{titarev2004finite}:

\begin{equation}
\Delta{t}=\text{CFL} \times min \left( \text{$min_{cells}$} \left (\frac{\Delta{x}}{|u+c|}, \frac{\Delta{y}}{|v+c|} \right ), \text{$min_{cells}$} \left (\frac{\Delta{x}^2}{\mu}, \frac{\Delta{y}^2}{\mu} \right )  \right ),
\end{equation}
where $c$ is the speed of sound and given by $c=\sqrt{\gamma{p}/\rho}$. The value of ${\alpha}$ in Equation (\ref{eqn:alpha}) puts a restriction on the CFL (Courant-Friedrichs-Lewy) number such that CFL $\leq$ 1/(1+${\alpha}$) (see Suresh and Hyunh \cite{suresh1997accurate} for further details). The values considered various schemes for parameter ${\alpha}$ for both single and multi-component flows are shown in Table \ref{tab:alpha}. Time integration is performed with a CFL = 0.2 for all single component flows and 0.1 for multi-component flows.

\begin{table}[H]
 \centering
 \footnotesize
 \caption{Parameter ${\alpha}$ for MP5, IG6MP and IG4MP schemes.}
  \begin{tabular}{ccccccc}
 \hline
     & MP5         & IG6MP &  IG4MP & \\
 \hline
   Single        & 4 &  7 &   7 \\
   Multi          & 4 &  4 &   4 \\
   \hline
  \end{tabular}%
 \label{tab:alpha}%
\end{table}%

\section{Results}\label{sec-4}
									
In this section, the proposed numerical scheme is tested for both single and multi-component flows.

\subsection{One-dimensional Euler equations}

In this subsection, we consider the test cases for the one-dimensional Euler equations.

\begin{example}\label{sod}{Shock tube problems}
\end{example}

The two shock-tube problems proposed by Sod \cite{sod1978survey} and Lax \cite{lax1954weak}, are solved by the proposed schemes, IG6MP and IG4MP. The solutions are obtained by setting the specific heat ratio to be $\gamma=1.4$, and are compared with that of results from exact Riemann solver \cite{toro2009riemann}. The initial conditions for the Sod test case and the Lax problem are given by the following initial conditions \eqref{sod_prob}, and \eqref{lax_prob}, respectively.

\begin{align}\label{sod_prob}
(\rho,u,p)=
\begin{cases}
(0.125,\ \ 0,\ \ 0.1),&\quad 0<x<0.5,\\
(1,\ \ 0,\ \ 1),&\quad 0.5 \leq x<1.
\end{cases}
\end{align}

\begin{align}\label{lax_prob}
(\rho,u,p)=
\begin{cases}
(0.445,\ \ 0.698,\ \ 3.528),&\quad 0<x<0.5,\\
(0.5,\ \ 0,\ \ 0.571),&\quad 0.5 \leq x<1,
\end{cases}
\end{align}

First, the Sod test case is used to assess the shock-capturing ability of the scheme. The initial condition is evolved for time $t=0.2$, and the solution is shown for $200$ cells. Fig. \ref{fig_sod} shows the density and velocity profiles for both the proposed schemes, and the MP5 scheme. The solutions by all the schemes are in good agreement with the exact solution in addition to the absence of overshoots in regions of discontinuities.

\begin{figure}[H]
\centering
\subfigure[Density]{\includegraphics[width=0.4\textwidth]{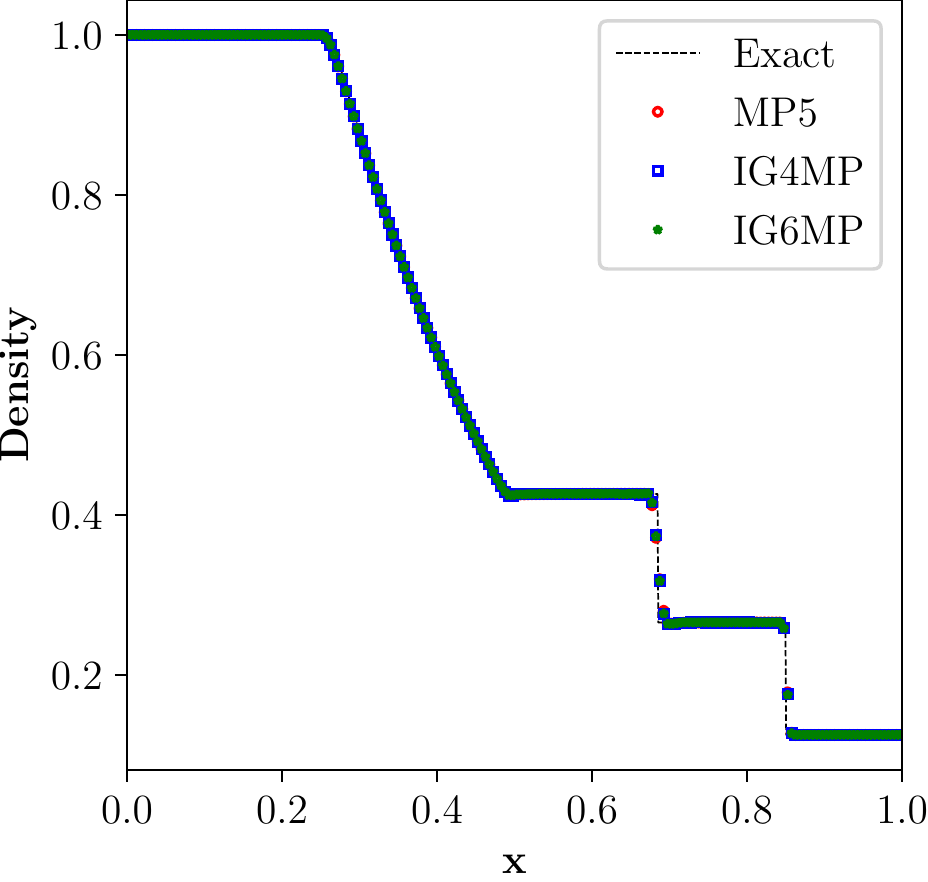}
\label{fig:sod-den}}
\subfigure[Pressure]{\includegraphics[width=0.4\textwidth]{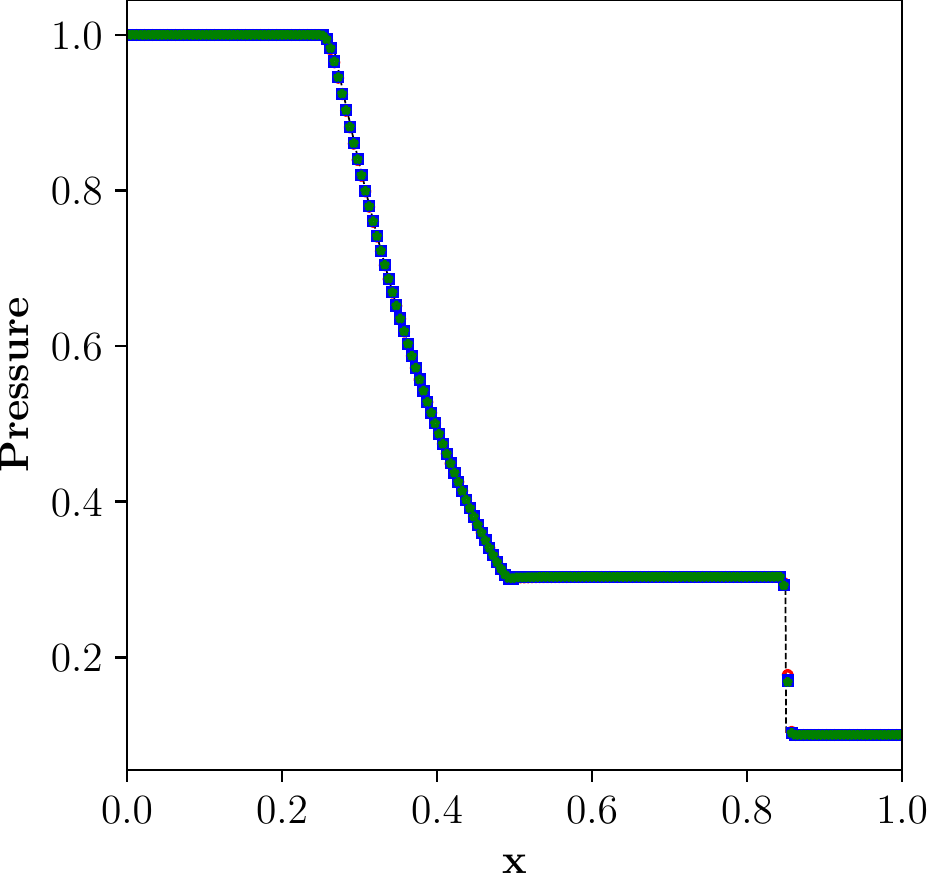}
\label{fig:sod-pres}}
\caption{Numerical solution for Sod problem in Example \ref{sod} for $N = 200$ points, where dashed line: reference solution; green stars: IG6MP; blue squares: IG4MP; red circles: MP5. }
\label{fig_sod}
\end{figure}

Second, we used $200$ cells for the Lax problem, and the solution is obtained at time $t=0.14$. Fig. \ref{fig_lax} presents the density and velocity of the new schemes and the MP5 scheme. The IG4MP scheme resolves well the features of the flow while avoiding oscillations.

\begin{figure}[H]
\centering
\subfigure[Density]{\includegraphics[width=0.4\textwidth]{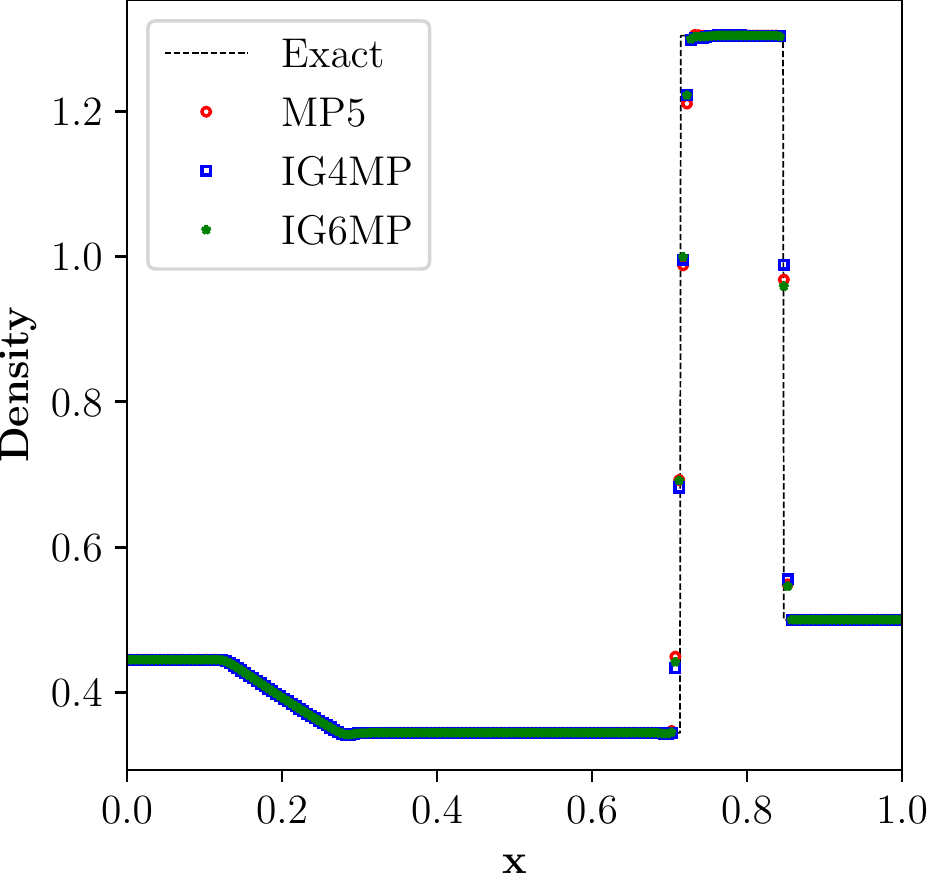}
\label{fig:lax-den}}
\subfigure[Velocity]{\includegraphics[width=0.4\textwidth]{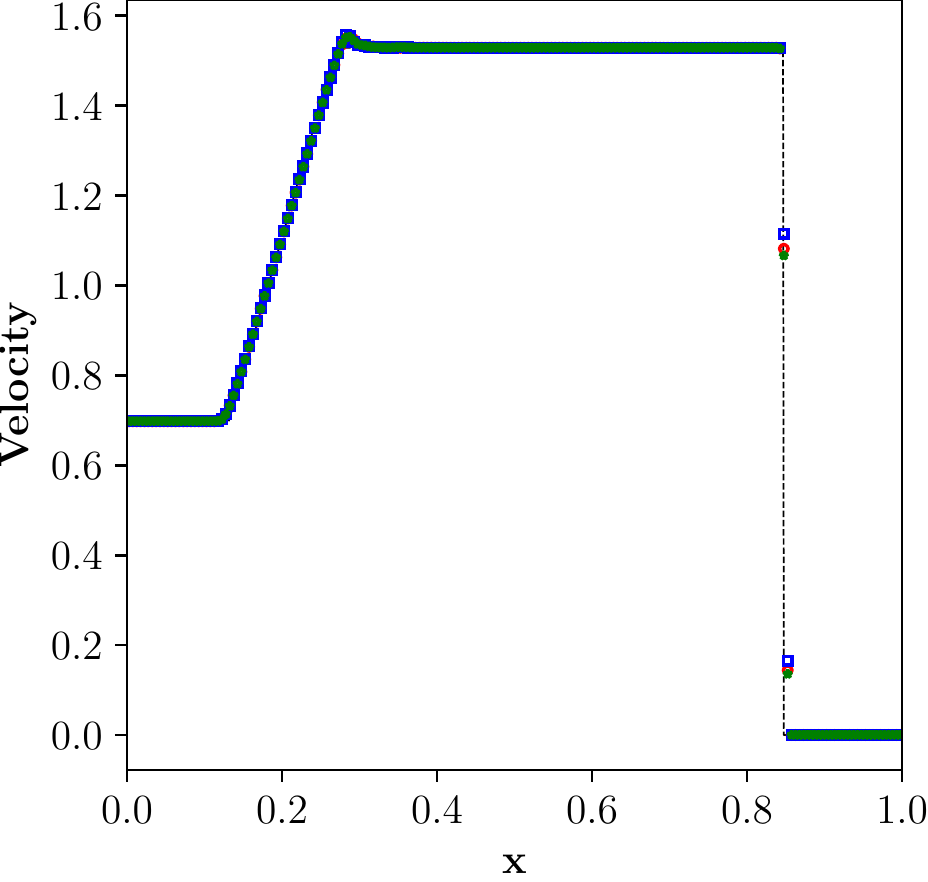}
\label{fig:lax-pres}}
\caption{Numerical solution for Lax problem in Example \ref{sod} for $N=200$ cells, where dashed line: reference solution; green stars: IG6MP; blue squares: IG4MP; red circles: MP5.}
\label{fig_lax}
\end{figure}


\begin{example}\label{Shu-Osher}{Shu-Osher problem}
\end{example}

Third, we consider the Shu-Osher problem \cite{Shu1988}, which is a one-dimensional idealization of shock-turbulence interaction, that simulates the interaction of a right moving shock wave for a given Mach number ($M=3$) superimposed with a perturbed density field. The initial condition for this problem are

\begin{align}\label{shock_den_prob}
(\rho,u,p)=
\begin{cases}
(3.857143,\ \ 2.629369,\ \ 10.3333),&\quad -5<x<-4,\\
(1+0.2\sin(5x),\ \ 0,\ \ 1),&\quad -4 \leq x<5,
\end{cases}
\end{align}

The solution is obtained for a time $t=1.8$ over a grid of $300$ cells. The reference solution is obtained using the WENOZ scheme \cite{Borges2008} on a fine grid of $1600$ cells. The solution for density profile is shown in Fig. \ref{fig:1d-SO}. It is observed that both IGMP schemes perform well in capturing the post-shock oscillations in density. Notably, both the proposed schemes capture the peaks and troughs of the density very well.

\begin{figure}[H]
\centering 
\subfigure[Global profile]{\includegraphics[width=0.4\textwidth]{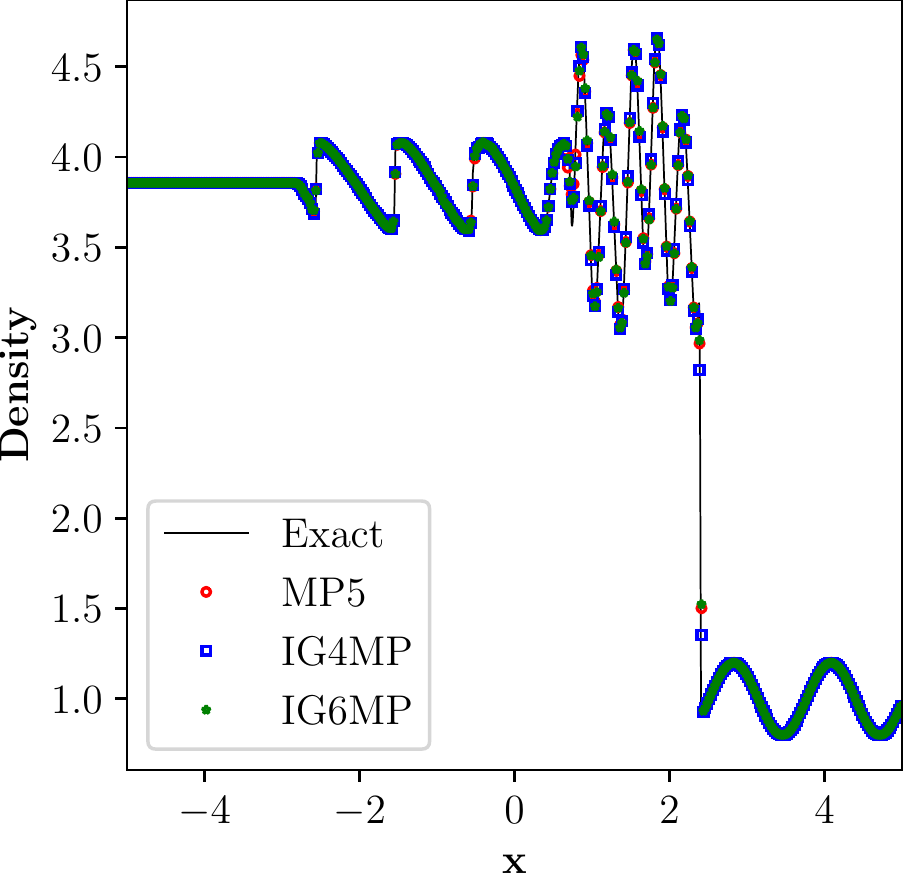}
\label{fig:shu-global}}
\subfigure[Local profile]{\includegraphics[width=0.42\textwidth]{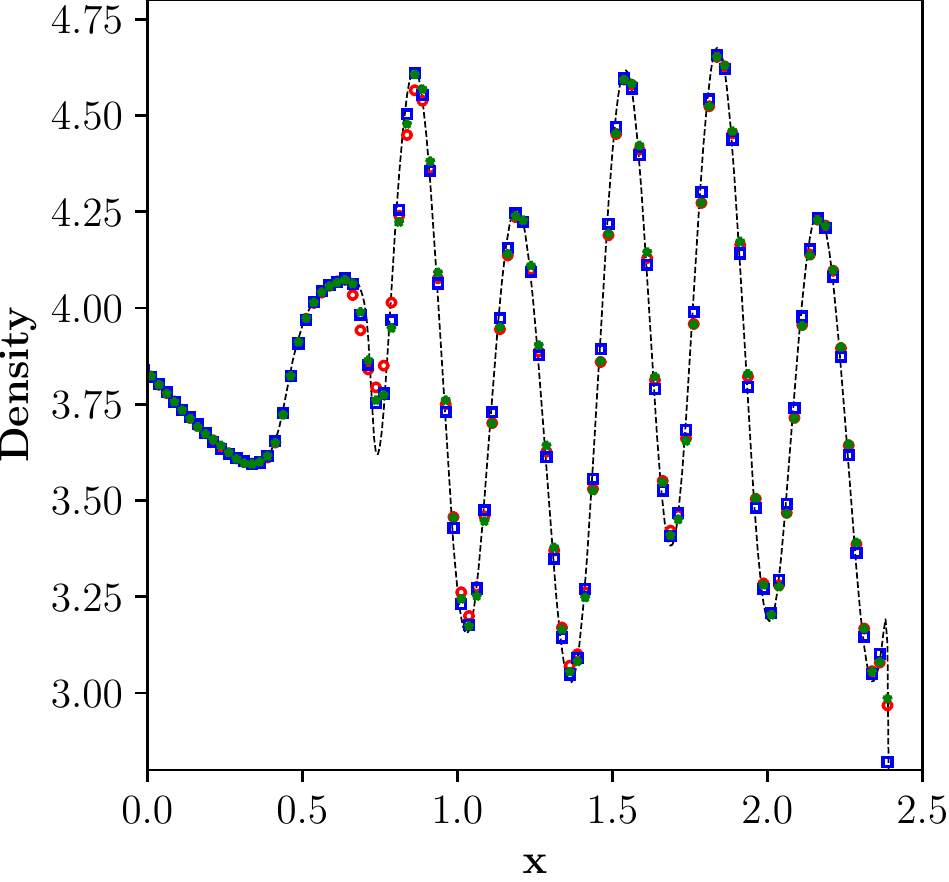}
\label{fig:shu-local}}
\caption{Density profile for Shu-Osher problem, Example \ref{Shu-Osher}, on a grid of $N=300$ cells, where solid or dashed line: reference solution; green stars: IG6MP; blue squares: IG4MP; red circles: MP5.}
\label{fig:1d-SO}
\end{figure}


\begin{example}\label{blast}{Blast wave problem}
\end{example}

In this one-dimensional test case, we consider the blast wave problem by \cite{woodward1984numerical}. This problem simulates the interaction of two blast waves with reflective boundaries on both ends of the domain. The initial condition is given by equation \eqref{bwave}. The simulation is performed for a time of $t=0.038$, on a grid size of $800$ cells. The exact solution is obtained by the WENOZ scheme on a fine grid of $1600$ points. From Fig. \ref{fig_blast}, it can be observed that both the schemes perform well compared to the MP5 solution, particularly in capturing the peak density profile.

 \begin{equation}\label{bwave}
(\rho, u, p)=\left\{\begin{array}{ll}
        (1.0,0.0,1000) ~~~~  0.0<x<0.1, \\
        (1.0,0.0,0.01) ~ ~~~~0.1<x<0.9.\\
        (1.0,0.0,100) ~ ~~~~~ 0.9<x<1.0.
        \end{array}\right.
\end{equation}

\begin{figure}[H]
\centering
\subfigure[Global solution]{\includegraphics[width=0.4\textwidth]{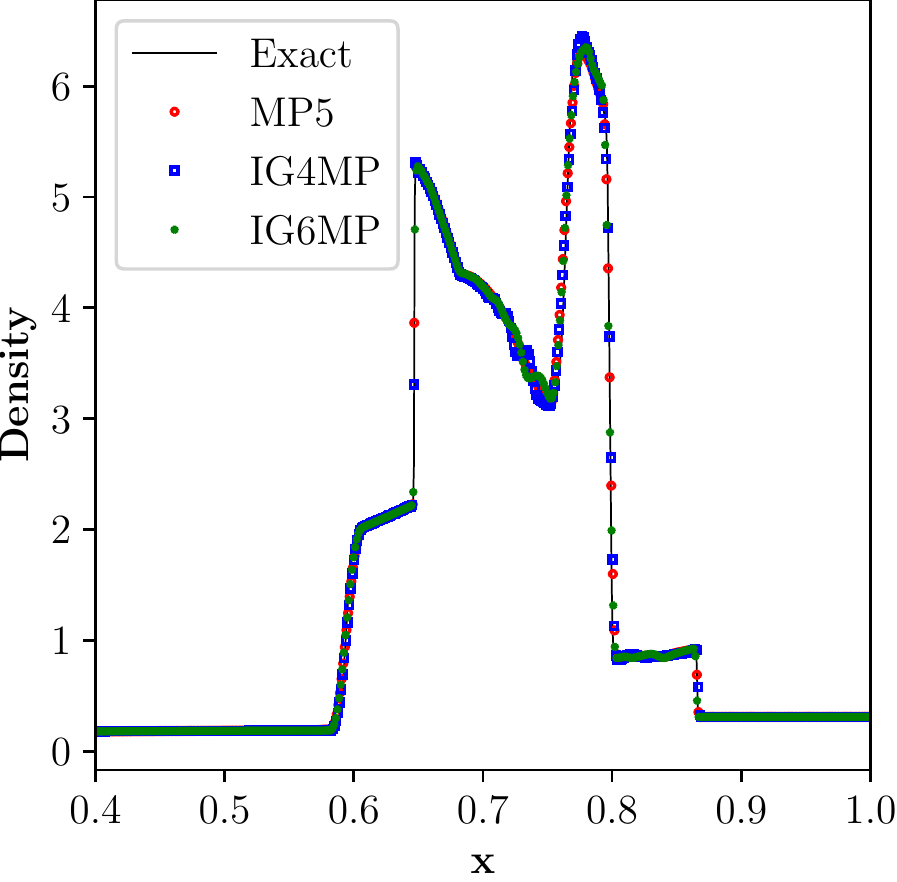}
\label{fig:blast-400}}
\subfigure[Local solution]{\includegraphics[width=0.42\textwidth]{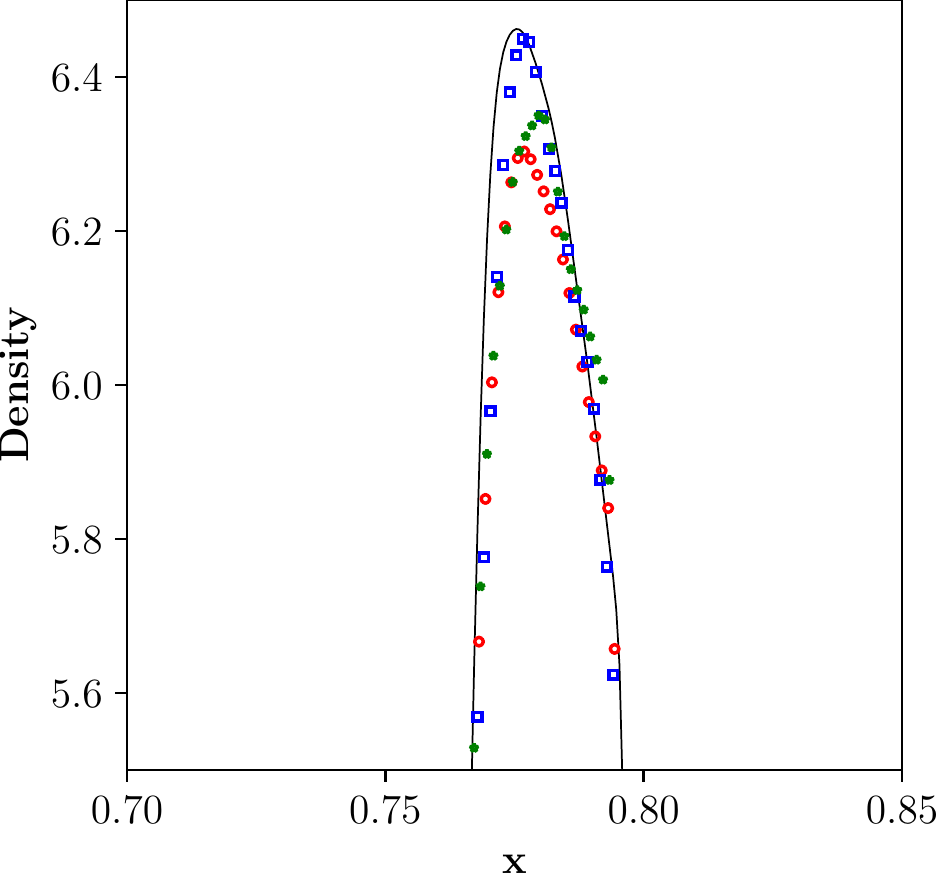}
\label{fig:blast-800}}
\caption{Density profiles obtained by various schemes for Example \ref{blast}. Global and local solution profiles obtained on a grid size of $N=800$ cells are shown in left and right panels, respectively. Solid line: reference solution; green stars: IG6MP; blue squares: IG4MP; red circles: MP5.}
\label{fig_blast}
\end{figure}


\begin{example}\label{Titarev-Toro}{Titarev-Toro problem}
\end{example}

The last one-dimensional test case we consider is the shock-entropy wave problem of Titarev-Toro \cite{titarev2004finite}. In this test, a high-frequency sinusoidal wave interacts with a shock wave. The test case reflects the ability of the scheme to capture the extremely high-frequency waves. The initial conditions are given by equation \eqref{shock_tita} on a domain of $[-5,5]$,

\begin{align}\label{shock_tita}
(\rho,u,p)=
\begin{cases}
(1.515695,\ \ 0.523326,\ \ 1.805),&\quad x<-4.5,\\
(1+0.1\sin(20x\pi),\ \ 0,\ \ 1),&\quad x \geq -4.5.
\end{cases}
\end{align}

 We used $1000$ cells to simulate the problem for time $t=5$, and compared with the reference solution obtained from a WENOZ simulation on a fine grid size of $3000$ cells. The results obtained presented in Fig. \ref{fig_tita} indicate that both the schemes can accurately capture the high-frequency wave. Specifically, we observe that the IGMP schemes capture the linear region significantly better than the MP5 scheme.

\begin{figure}[H]
\centering
\subfigure[Global profile]{\includegraphics[width=0.4\textwidth]{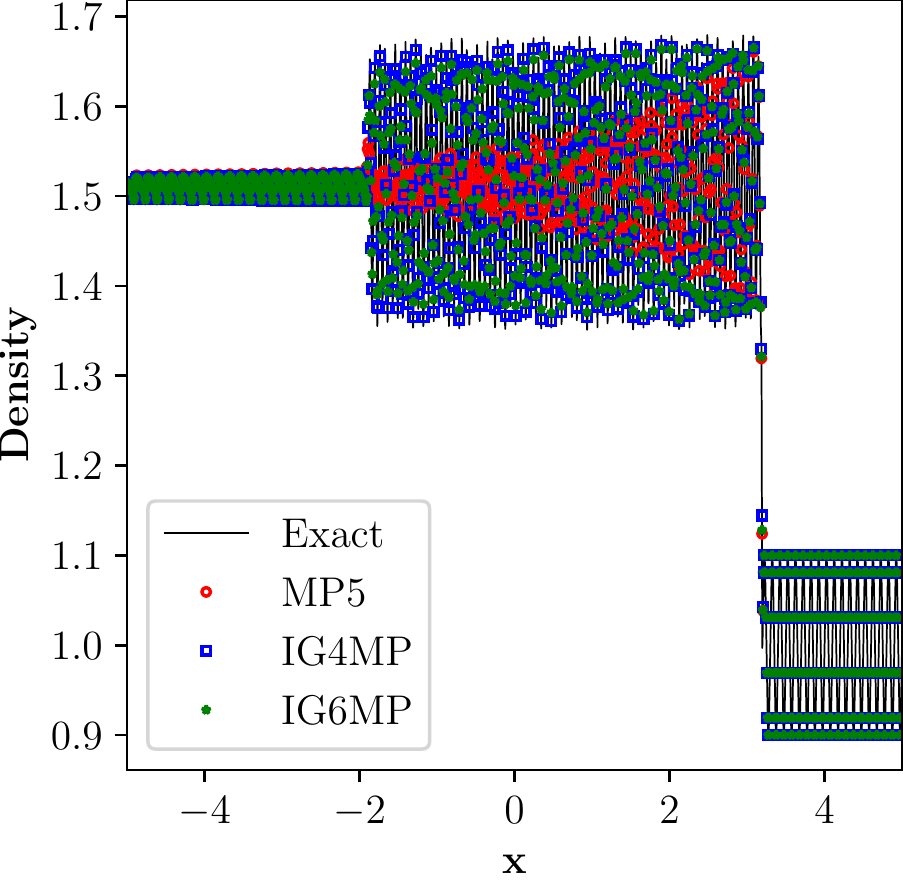}
\label{fig:tita1}}
\subfigure[Local profile]{\includegraphics[width=0.4\textwidth]{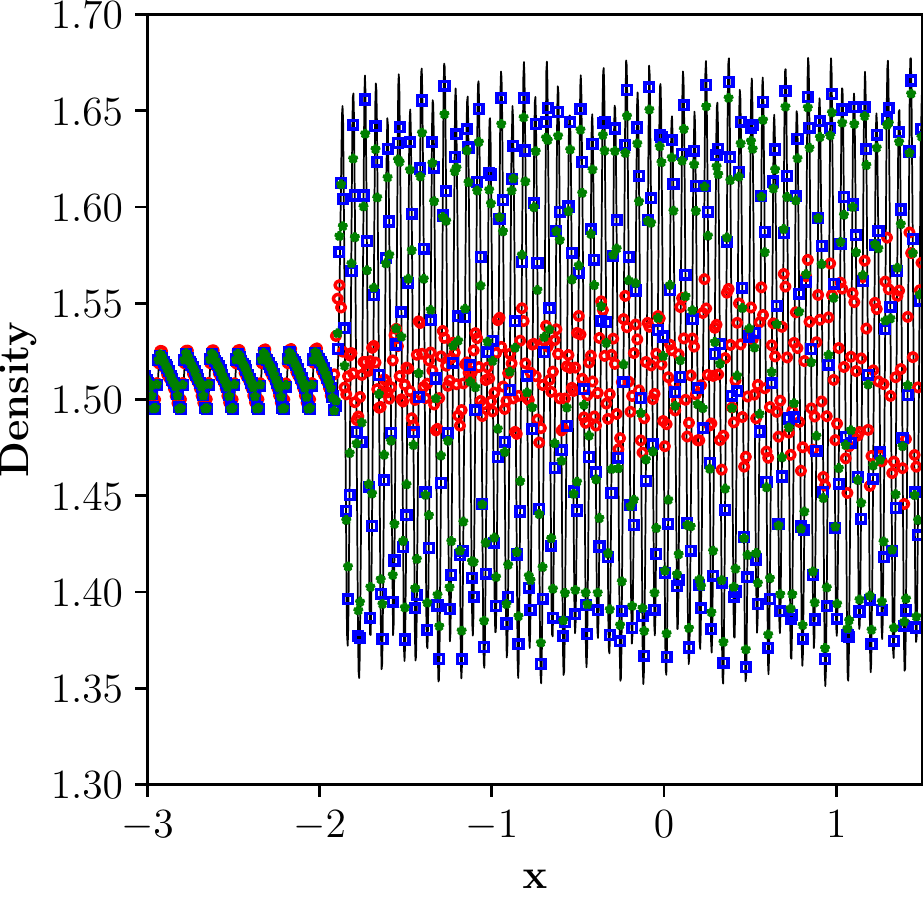}
\label{fig:tita2}}
\caption{Density profiles obtained by various schemes for Example \ref{Titarev-Toro}. Solutions are obtained for grid size of $1000$ cells. Solid line: reference solution; green stars: IG6MP; blue squares: IG4MP; red circles: MP5.}
\label{fig_tita}
\end{figure}

\subsection{Multi-dimensional test cases for Euler equations}

In this section, we present the numerical results of the proposed method for multi-dimensional Euler and Navier-Stokes equations.


\begin{example}\label{shock-entropy}{2D shock-entropy wave test}
\end{example}
In this test case we consider the two-dimensional shock-entropy wave interaction problem proposed in \cite{acker2016improved}. The initial conditions for the test case are as follows,
\begin{align}\label{shock_entropy}
(\rho,u,v,p)=
\begin{cases}
(3.857143, \ \ 2.629369,\ \ 0,\ \ 10.3333),&\quad x\leq-4,\\
(1+0.2\sin(10x \cos\theta+10y\sin\theta),\ \ 0,\ \ 0,\ \ 1),&\quad otherwise,
\end{cases}
\end{align}
with $\theta$ = $\pi/6$ over a domain of $[-5,5]\times [-1,1]$. The initial sine waves make an angle of $\theta$ radians with the $x$ axis. Initial conditions are modified as in \cite{deng2019fifth} with a higher frequency for the initial sine waves compared to that of \cite{acker2016improved} to show the benefits of the proposed method. Mesh size of $400\times 80$ is chosen and the flow is developed for time $t=1.8$. The reference solution is computed on a fine mesh of $1600 \times 320$ by using the WENOZ scheme. Density contour plots shown in Fig. \ref{fig_shock_entropy} indicate that the proposed schemes significantly improve the resolution of the flow structures. The local density profile along $y=0$ is presented in Fig. \ref{fig:SSE-Compare}. The results demonstrate that IG6MP and IG4MP retain the desirable shock-capturing features in the MP5 scheme while capturing the high-frequency region better than the MP5 scheme. 

\begin{figure}[H]
\centering\offinterlineskip
\subfigure[MP5]{\includegraphics[width=0.48\textwidth]{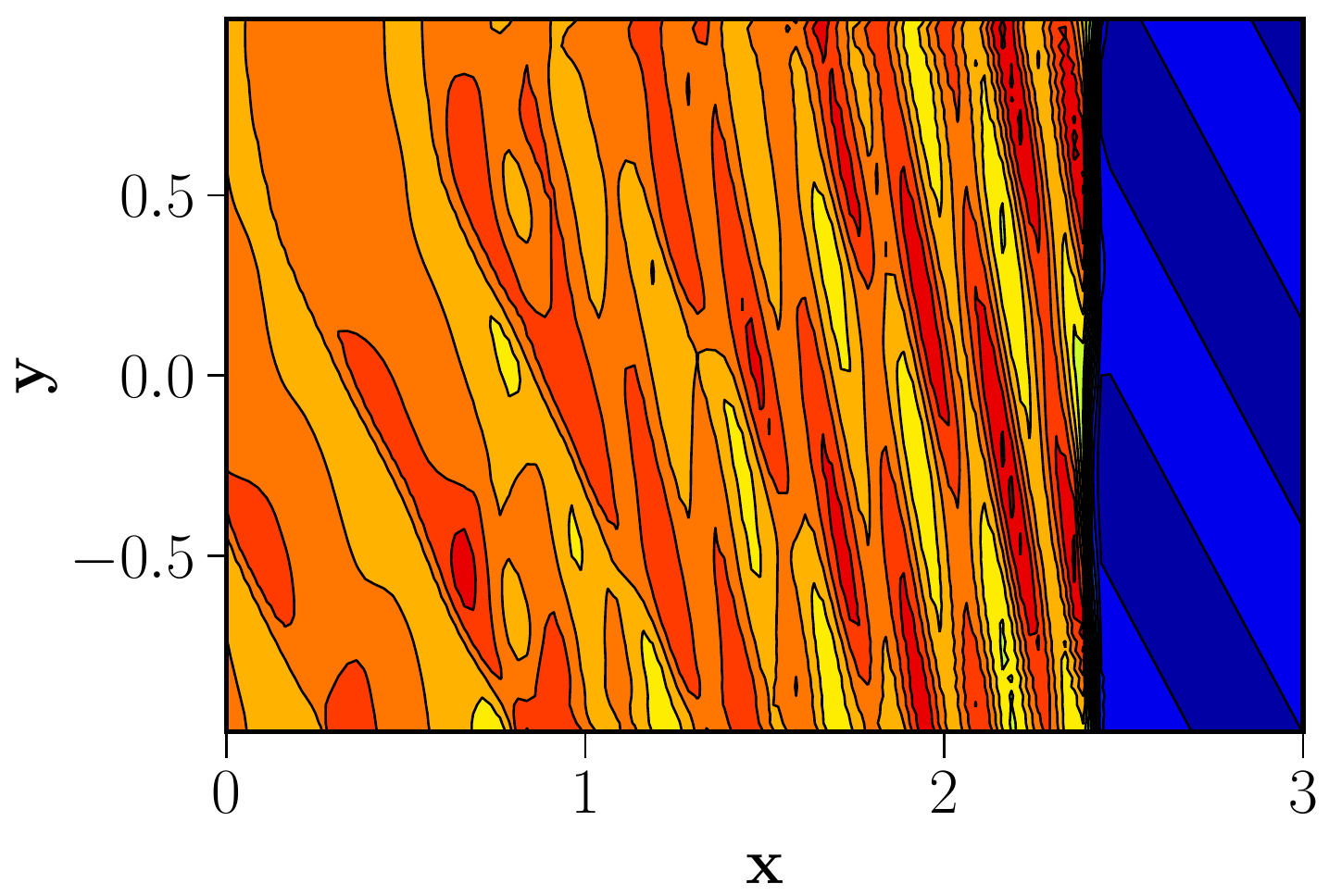}
\label{fig:SE-MP5}}
\subfigure[IG6MP]{\includegraphics[width=0.48\textwidth]{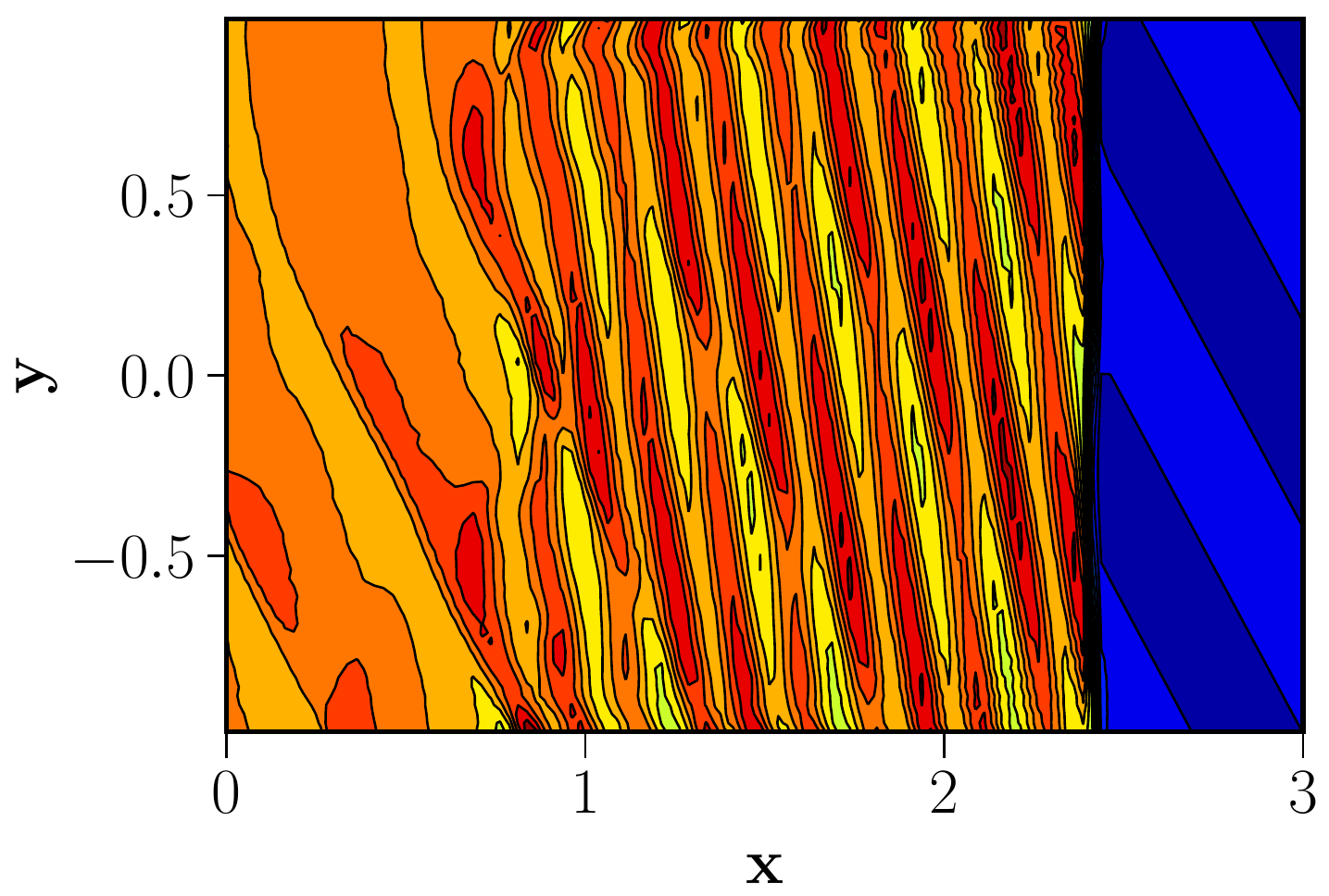}
\label{fig:SE-BVD5}}
\subfigure[IG4MP]{\includegraphics[width=0.48\textwidth]{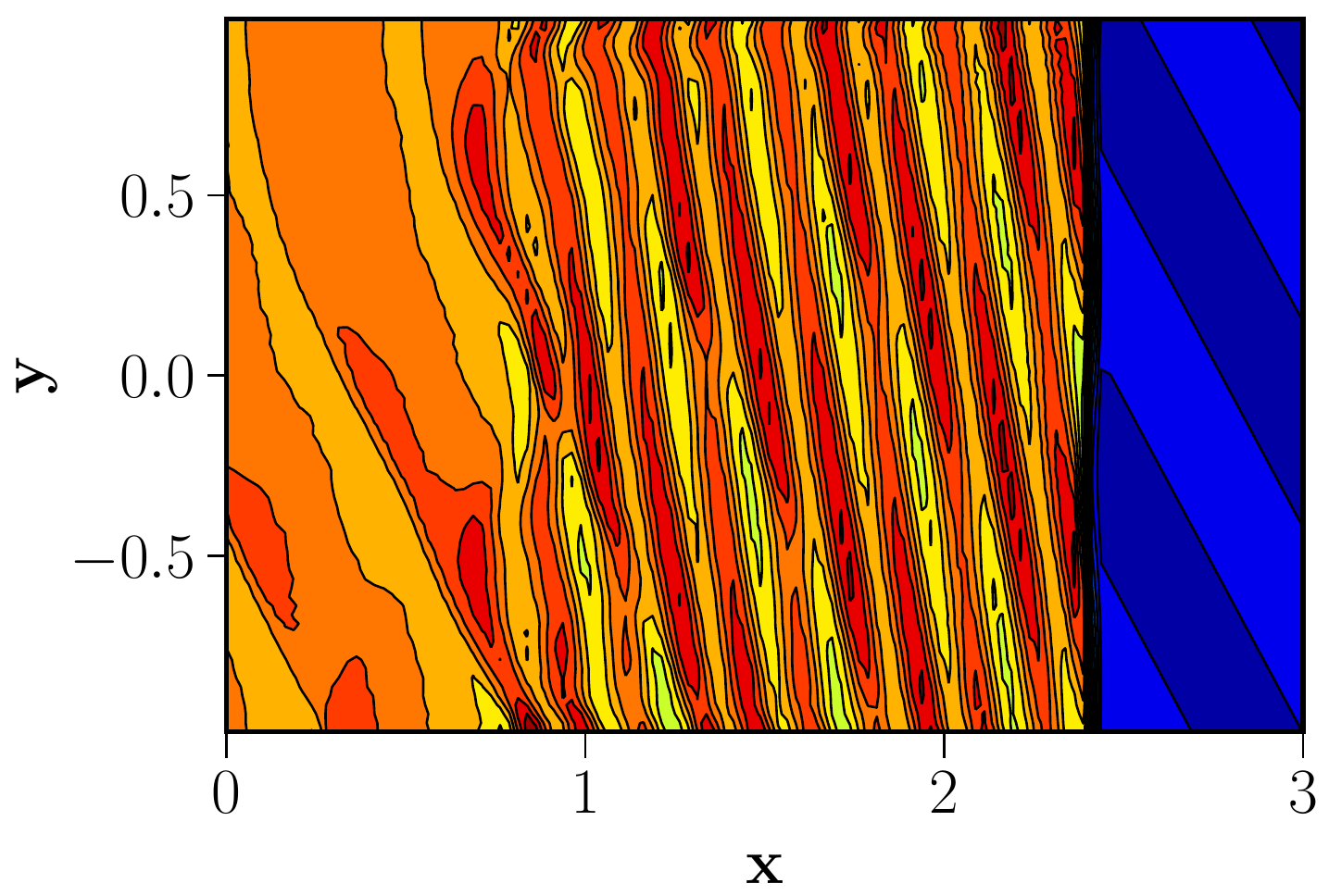}
\label{fig:SE-BVD}}
\subfigure[Local profile]{\includegraphics[width=0.48\textwidth]{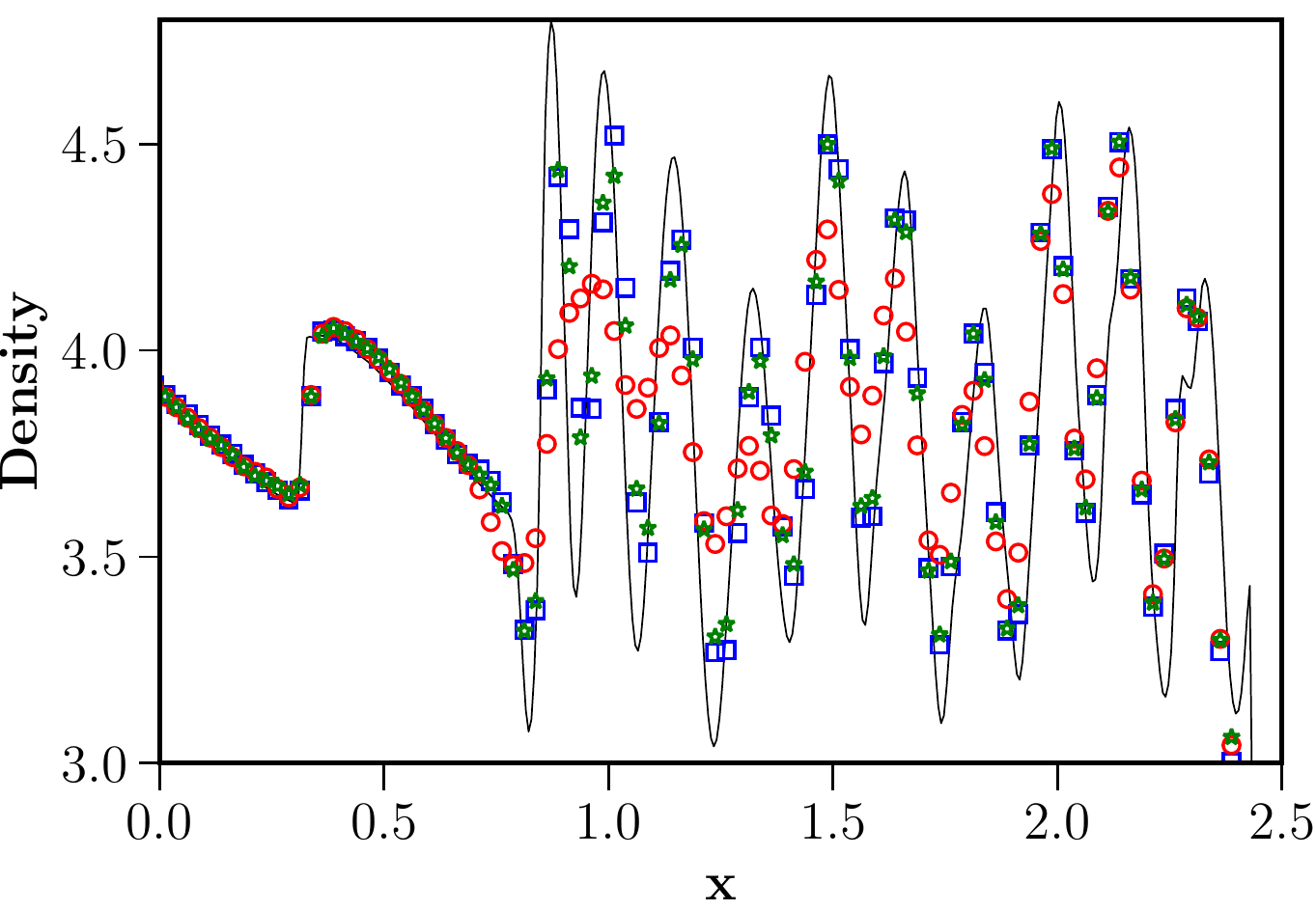}
\label{fig:SSE-Compare}}
\caption{Density contours for the 2D shock-entropy wave test at $t=1.8$, Example \ref{shock-entropy}, for various schemes are shown in Figs. (a), (b) and (c). Fig. (d) shows the local density profile in the region with high-frequency waves for all the schemes. Solid line: reference solution; green stars: IG6MP; blue squares: IG4MP; red circles: MP5.}
\label{fig_shock_entropy}
\end{figure}


\begin{example}\label{ex:rp}{Riemann Problem }
\end{example}

In this test case we consider the Riemann problem of \cite{schulz1993numerical} described as configuration 3. The initial conditions of the problem are given by equation \eqref{riemann_problem} with constant states of the primitive variables along the lines $x=0.8$, and $y=0.8$ in the domain $x,y \in [0,1]$. This produces four shocks at the interfaces of the four quadrants. Also, the small-scale complex structures generated along the slip-lines due to the Kevin-Helmholtz instabilities serve to assess the numerical dissipation of the scheme. Non-reflective boundary conditions are employed on all four boundaries. The test case is run for a time $t=0.8$ on a grid of size $400 \times 400 $. 

 \begin{equation}\label{riemann_problem}
(\rho, u,v, p)=\left\{\begin{array}{ll}
        (1.5, 0, 0, 1.5), ~~~~~~~~~~~~~~~~~~~~~~~~~~~~~~~~\mbox{if} ~~x > 0.8, ~~y > 0.8, \\
        (33/62, 4/\sqrt{11}, 0, 0.3), ~~~~~~~~~~~~~~~~~~\mbox{if} ~~ x\leq 0.8, ~~y > 0.8, \\
        (77/558, 4/\sqrt{11}, 4/\sqrt{11}, 9/310), ~~\mbox{if} ~~ x \leq 0.8, ~~y\leq 0.8, \\
        (33/62, 0, 4/\sqrt{11}, 0.3), ~~~~~~~~~~~~~~~~~~\mbox{if} ~~x > 0.8, ~~ y\leq 0.8.
        \end{array}\right.
\end{equation}

The results are presented for the IG6MP and the IG4MP schemes, and the density contours are presented in Fig. \ref{fig_riemann}. The proposed schemes resolve better rollup behaviour along the slip lines compared to the MP5 scheme. This richness of the fine structures indicate low numerical dissipation features for both the proposed schemes. Specifically, the guitar-like shape of the jet is well captured

\begin{figure}[H]
\begin{onehalfspacing}
\centering\offinterlineskip
\subfigure[MP5]{\includegraphics[width=0.48\textwidth]{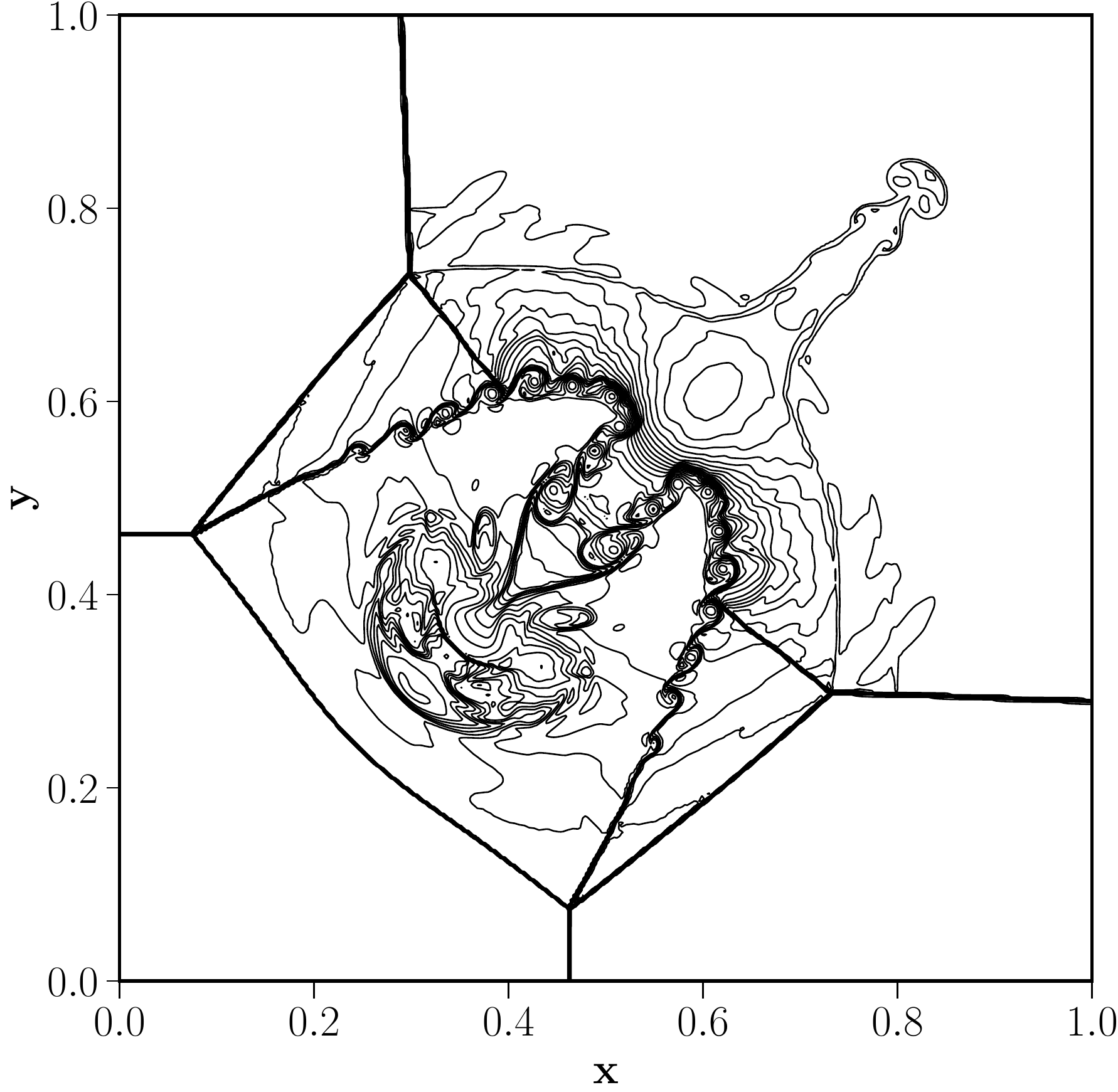}
\label{fig:MP5_RM}}
\subfigure[IG6MP]{\includegraphics[width=0.48\textwidth]{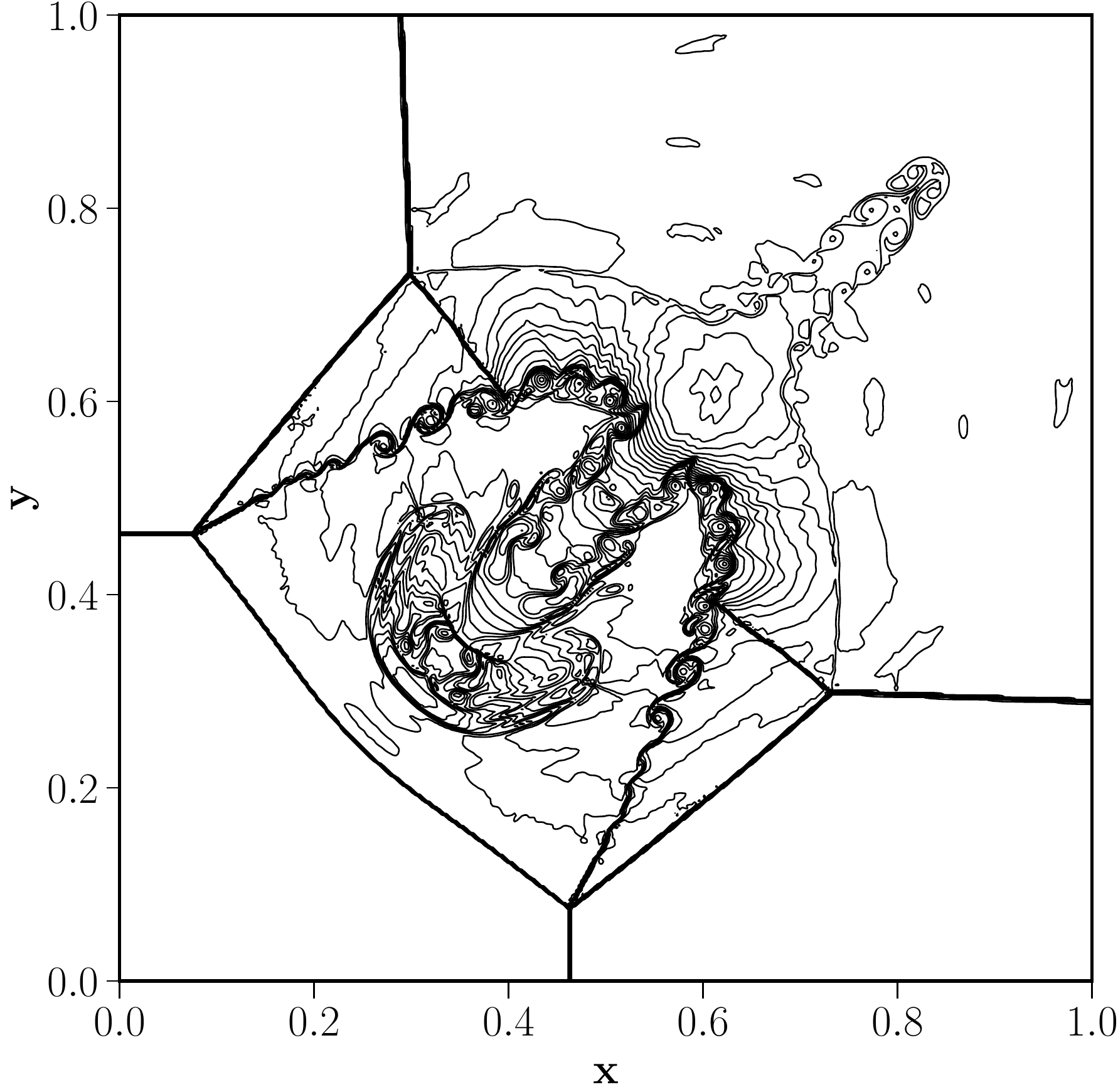}
\label{fig:IG6MC_RM}}
\subfigure[IG4MP]{\includegraphics[width=0.48\textwidth]{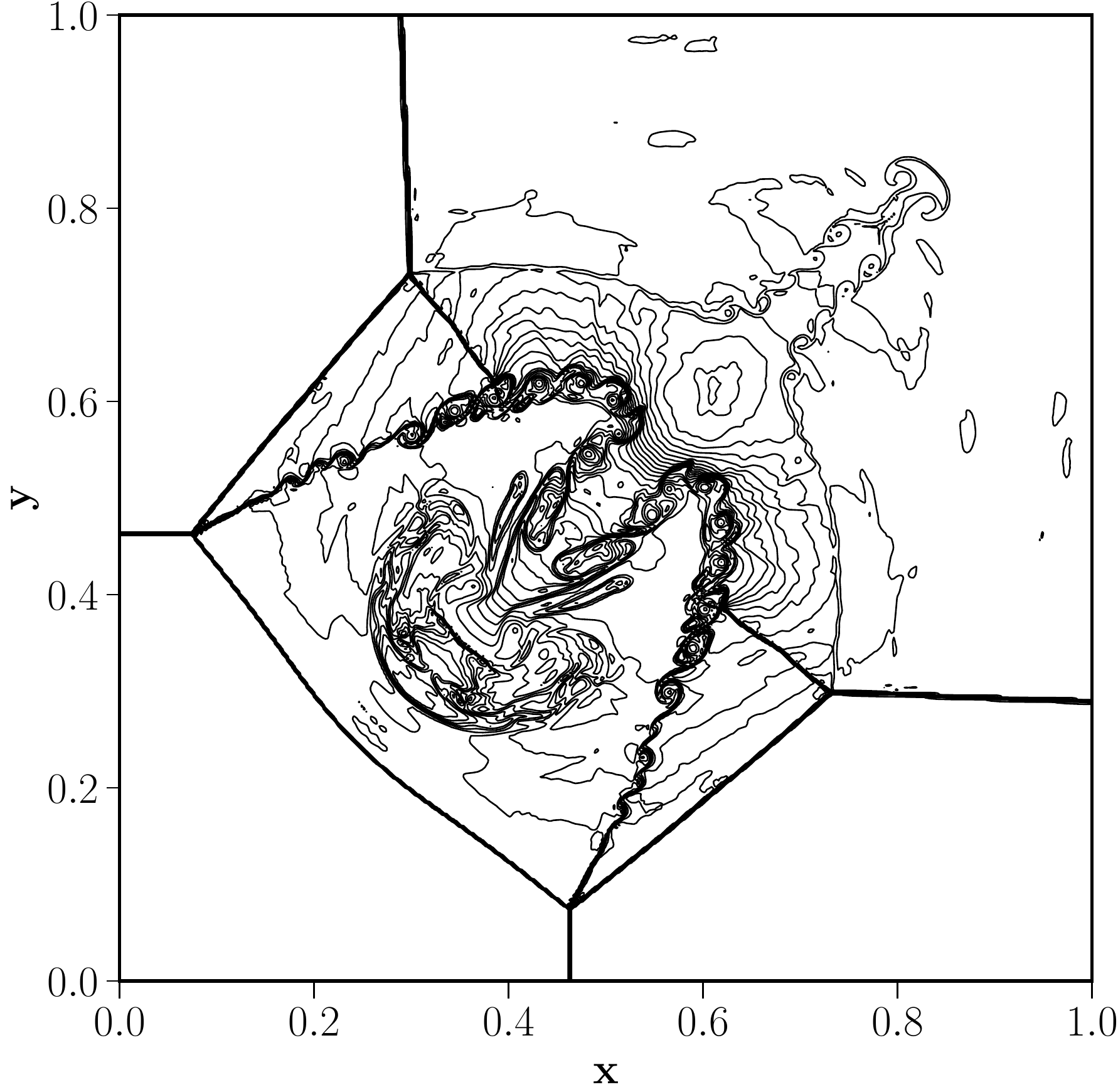}
\label{fig:IG4MP_RM}}
\caption{Density contours of the Riemann problem in Example \ref{ex:rp} for different schemes on a grid size of 400 $\times$ 400.}
\label{fig_riemann}
\end{onehalfspacing}
\end{figure}


\begin{example}\label{ex:rt} {{\color{black}Rayleigh-Taylor instability}}
\end{example}

Rayleigh-Taylor instability occurs at the interface between fluids with different densities when acceleration is directed from the heavy fluid to the light one. In this test case, two initial gas layers with different densities are subjected to unit magnitude's gravity, where the resulting acceleration is directed towards the lighter fluid. A small disturbance of the contact line triggers the instability. This problem has been extensively studied using high order shock-capturing schemes in the literature, see, e.g. \cite{Shi2002}, with the following initial conditions,

 \begin{equation}
\begin{aligned}
(\rho,u,v,p)=
\begin{cases}
(2.0,\ 0,\ -0.025\sqrt{\frac{5p}{3\rho}\cos(8\pi x)},\ 2y+1.00),&\quad 0\leq y< 0.5,\\
(1.0,\ 0,\ -0.025\sqrt{\frac{5p}{3\rho}\cos(8\pi x)},\ 1y+3/2),&\quad 0.5\leq y\leq 1.0,
\end{cases}
\end{aligned}
\label{eu2D_RT}
\end{equation}
over the computational domain $[0, 1/4]\times [0,1]$. Reflective boundary conditions are imposed on the right and left boundaries via ghost cells. The flow conditions are set to $\rho=1$, $p=2.5$, and $u=v=$0 on top boundary and $\rho=2$, $p=1.0$, and $u=v=0$ on bottom boundary with the specific heat ratio, $\gamma$, of $5/3$. The source term $S=(0,0,\rho, \rho v)$ is added to the Euler equations. We performed simulations on a uniform mesh of resolution $120 \times 480$ and the computations are conducted until $t = 1.95$.  

\begin{figure}[H]
\begin{onehalfspacing}
\centering\offinterlineskip
\subfigure[MP5]{%
\includegraphics[width=0.20\textwidth]{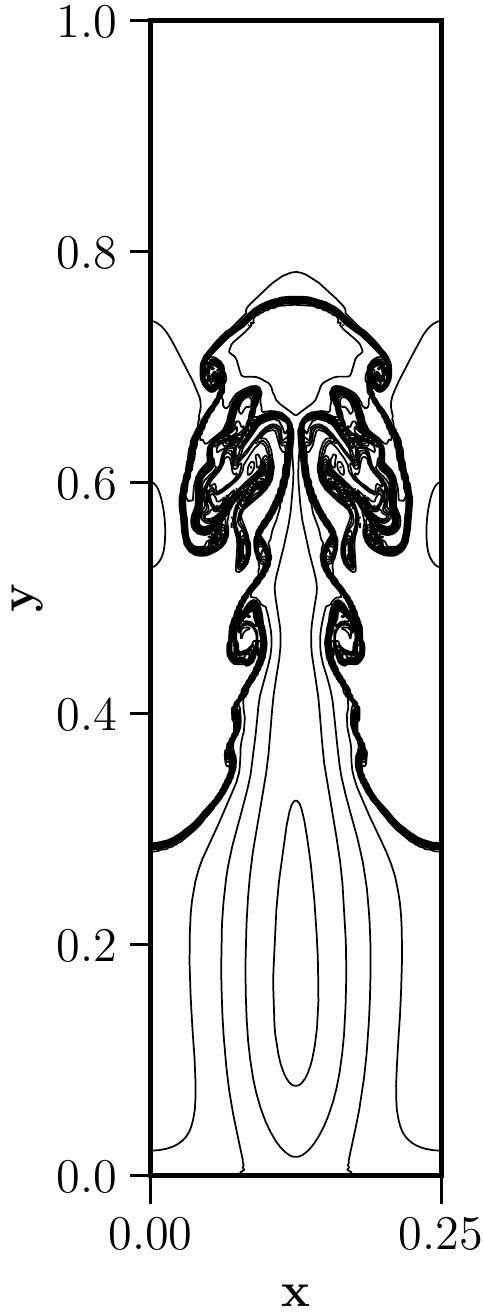}
\label{fig:RT_MP51}}
\subfigure[IG6MP]{%
\includegraphics[width=0.20\textwidth]{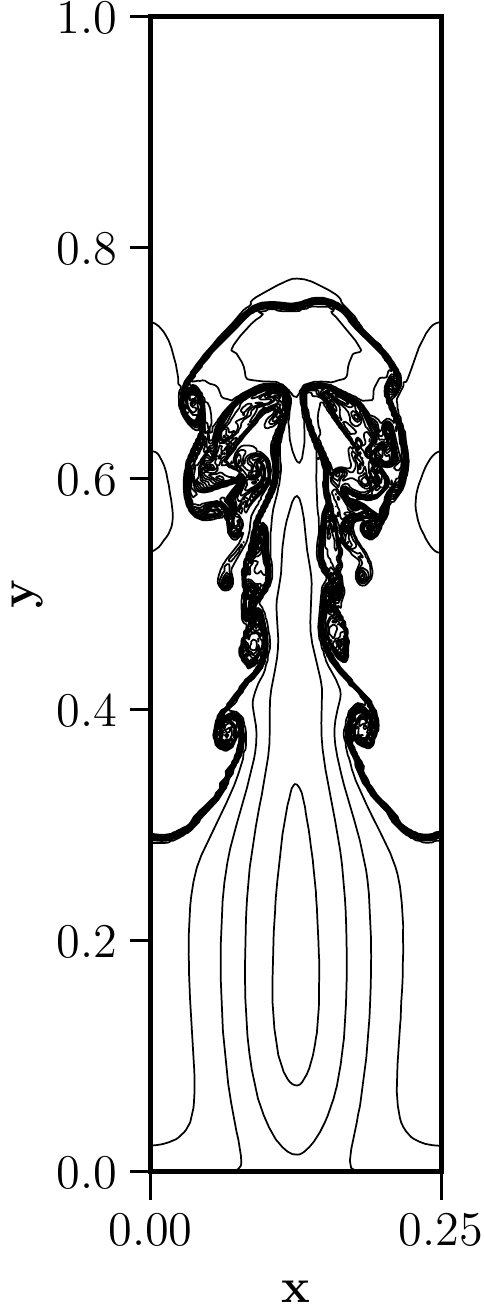}
\label{fig:RT_IG6T}}
\subfigure[IG4MP]{%
\includegraphics[width=0.20\textwidth]{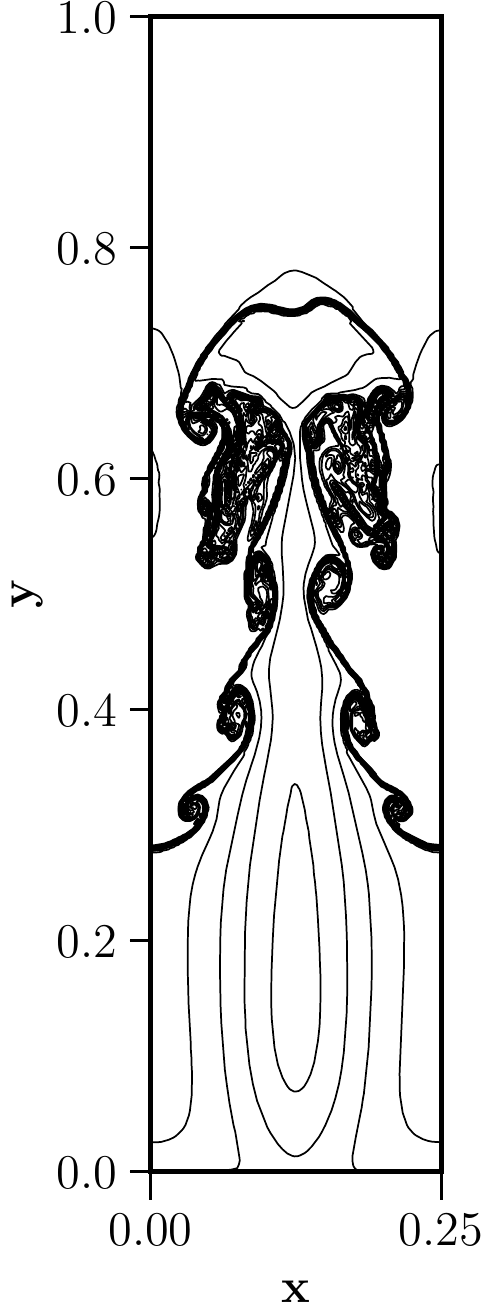}
\label{fig:RT_IG4T}}

\caption{Comparison of density contours obtained by different numerical schemes for the test case in Example \ref{ex:rt} on a grid size of 480 $\times$ 120.}
\label{fig:2d-RT1}
\end{onehalfspacing}
\end{figure}
Fig. \ref{fig:2d-RT1} shows the density distribution of the Rayleigh-Taylor instability problem. We can observe that the impilcit gradient schemes produced more small vortices in the shear layer, indicating that they have better resolution to capture small scale features of the flow.

\begin{example}\label{ex:kh}{Kelvin Helmholtz instability}
\end{example}
The Kelvin–Helmholtz (KH) instability occurs where there is a velocity difference across the interface between two fluids, and it plays an important role in the evolution of the mixing layer and the transition to turbulence. Consider the following initial conditions over a periodic domain of [0, 1] × [0, 1],

\begin{equation}
\begin{array}{l}
p=2.5, \quad \rho(x, y)=\left\{\begin{array}{l}
2, \text { If } 0.25<y \leq 0.75 \\
1, \text { else, }
\end{array}\right. \\
u(x, y)=\left\{\begin{array}{l}
0.5, \quad \text { If } 0.25<y \leq 0.75 \\
-0.5, \text { else },
\end{array}\right. \\
v(x, y)=0.1 \sin (4 \pi x)\left\{\exp \left[-\frac{(y-0.75)^{2}}{2 \sigma^{2}}\right]+\exp \left[-\frac{(y-0.25)^{2}}{2 \sigma^{2}}\right]\right\}, \text{where}\ \sigma = 0.05/\sqrt {2}.
\end{array}
\end{equation}
The computational domain is discretized with  512 cells in each direction, and the final time is taken to be $t$ = 0.8. The numerical solutions are computed using the IGMP schemes are depicted in Fig. 19. We can observe that the IGMP schemes capture complex structures and small-scale vortices when compared to the MP5 scheme.

\begin{figure}[H]
\centering\offinterlineskip
\subfigure[MP5]{\includegraphics[width=0.34\textheight]{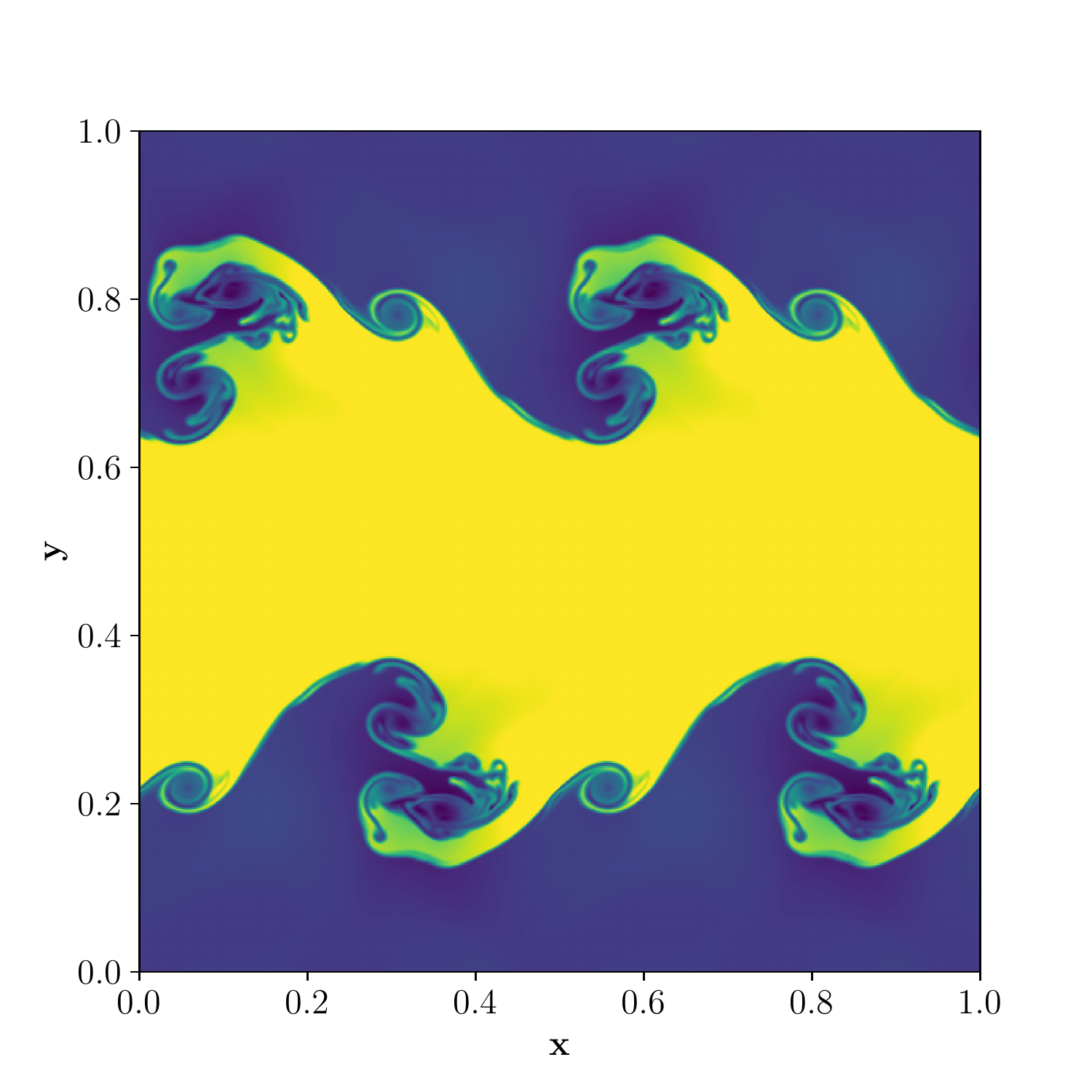}
\label{fig:MP5_KH}}
\subfigure[IG6MP]{\includegraphics[width=0.34\textheight]{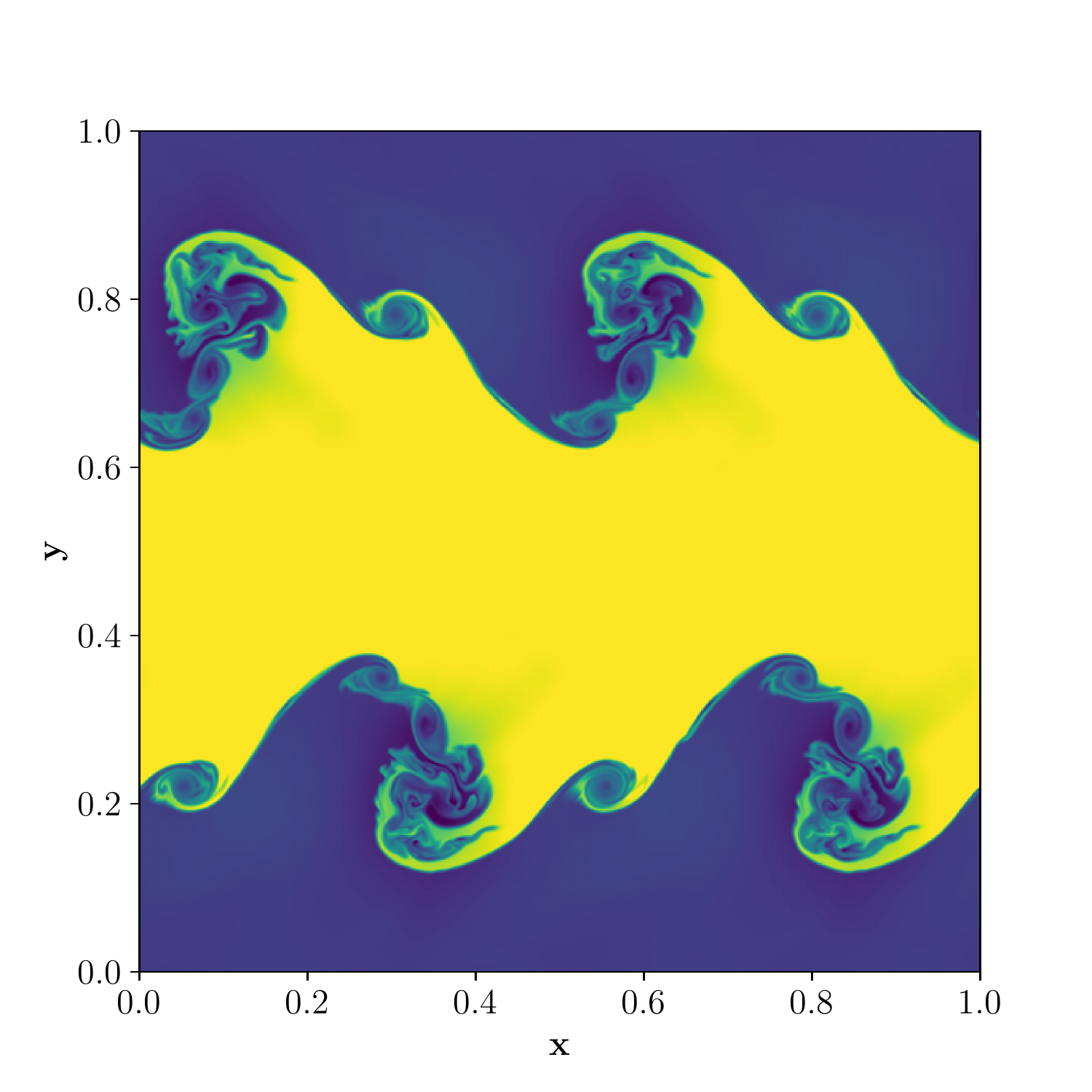}
\label{fig:IG6_KH}}
\subfigure[IG4MP]{\includegraphics[width=0.34\textheight]{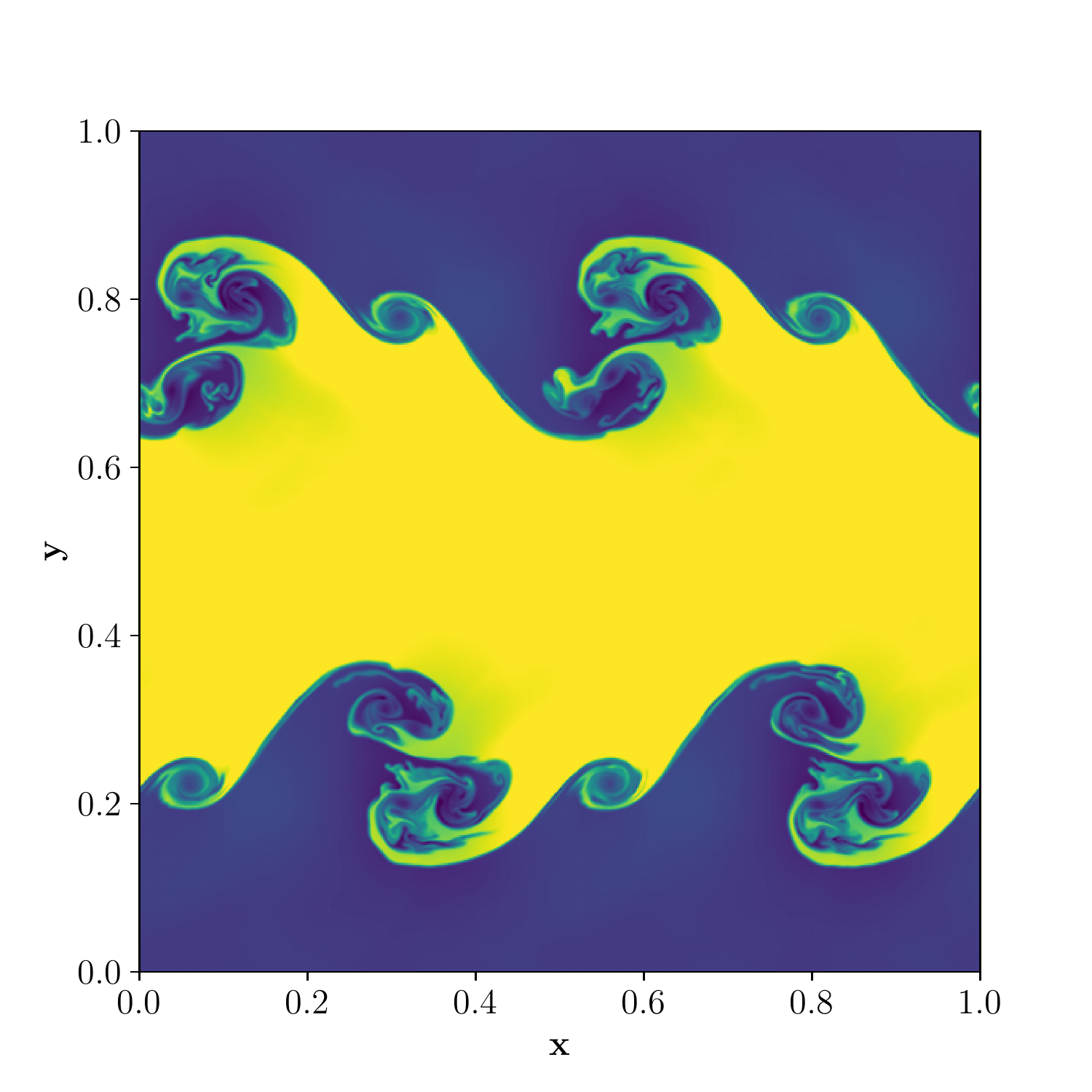}
\label{fig:IG4_KH}}
\caption{Numerical results of Kelvin-Helmholtz instability test case, Example \ref{ex:kh}, by different schemes on a grid size of 512 $\times$ 512.}
\label{fig_KH}
\end{figure}

\begin{example}\label{ex:dmr}{Double Mach reflection}
\end{example}

Next, we consider the double-Mach reflection problem proposed by \cite{woodward1984numerical}, and is typically used to test resolution capability of numerical schemes. In this test case, an unsteady planar shock-wave of Mach 10 impinges on an inclined surface of 30 degrees to the horizontal axis. This inclined surface is simplified by tilting the shock-wave to avoid modelling the oblique no-slip physical wall boundary. The near-wall jet structure and the vortex structures appearing from the contact discontinuity that emerges from the triple-point indicate the proposed scheme's numerical dissipation.

Post-shock flow conditions are set at the left boundary, and zero gradient conditions are applied at the right boundary. At the bottom boundary, reflecting boundary conditions are applied in $x \in [1/6,3]$ and the post-shock conditions in $x \in \left[ 0, 1/6 \right]$. Furthermore, the exact solution of the moving shock is imposed at the upper at $y=1$ and is time-dependent. The computational domain is taken as $x \in [0,3], y \in [0,1]$ and the simulation is performed for time $t=0.2$ on a grid of $768 \times 256$ cells.

\begin{equation}
\begin{aligned}
(\rho,u,v,p)=
\begin{cases}
&(8,\ 8.25 \cos 30^\circ,\ -8.25 \sin 30^\circ,\
116.5),\quad x<1/6+\frac{y}{\tan 60^\circ},\\
&(1.4,\ 0,\ 0,\
1),\quad\quad\quad\quad\quad\quad\quad\quad\quad\quad\quad\ x>
1/6+\frac{y}{\tan 60^\circ}.
\end{cases}
\end{aligned}
\label{eu2D_mach}
\end{equation}

Observations made from Figs. \ref{fig_doublemach} and  \ref{fig_doublemach-full} indicate that the IG4MP and the IG6MP schemes have better resolution of the KH instabilities. Notably the resolution of the shear layers along the slip lines and the near-wall jet region are well resolved. It can be noted that the current IG4MP scheme is slightly better in resolving the shear layer along the slip line than the IG6MP scheme. 
\begin{figure}[H]
\centering\offinterlineskip
\subfigure[MP5]{\includegraphics[width=0.24\textheight]{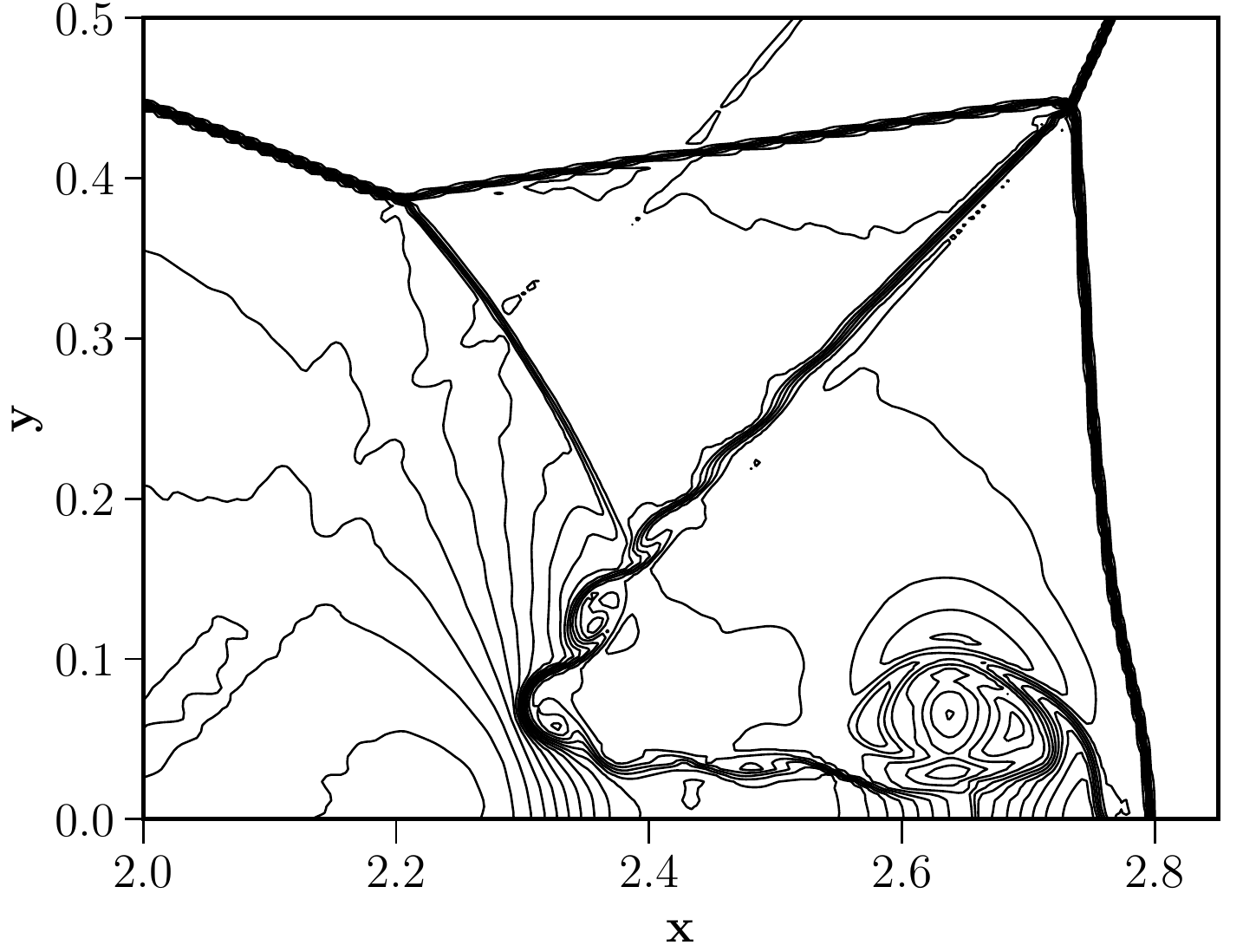}
\label{fig:MP5_DB}}
\subfigure[IG6MP]{\includegraphics[width=0.24\textheight]{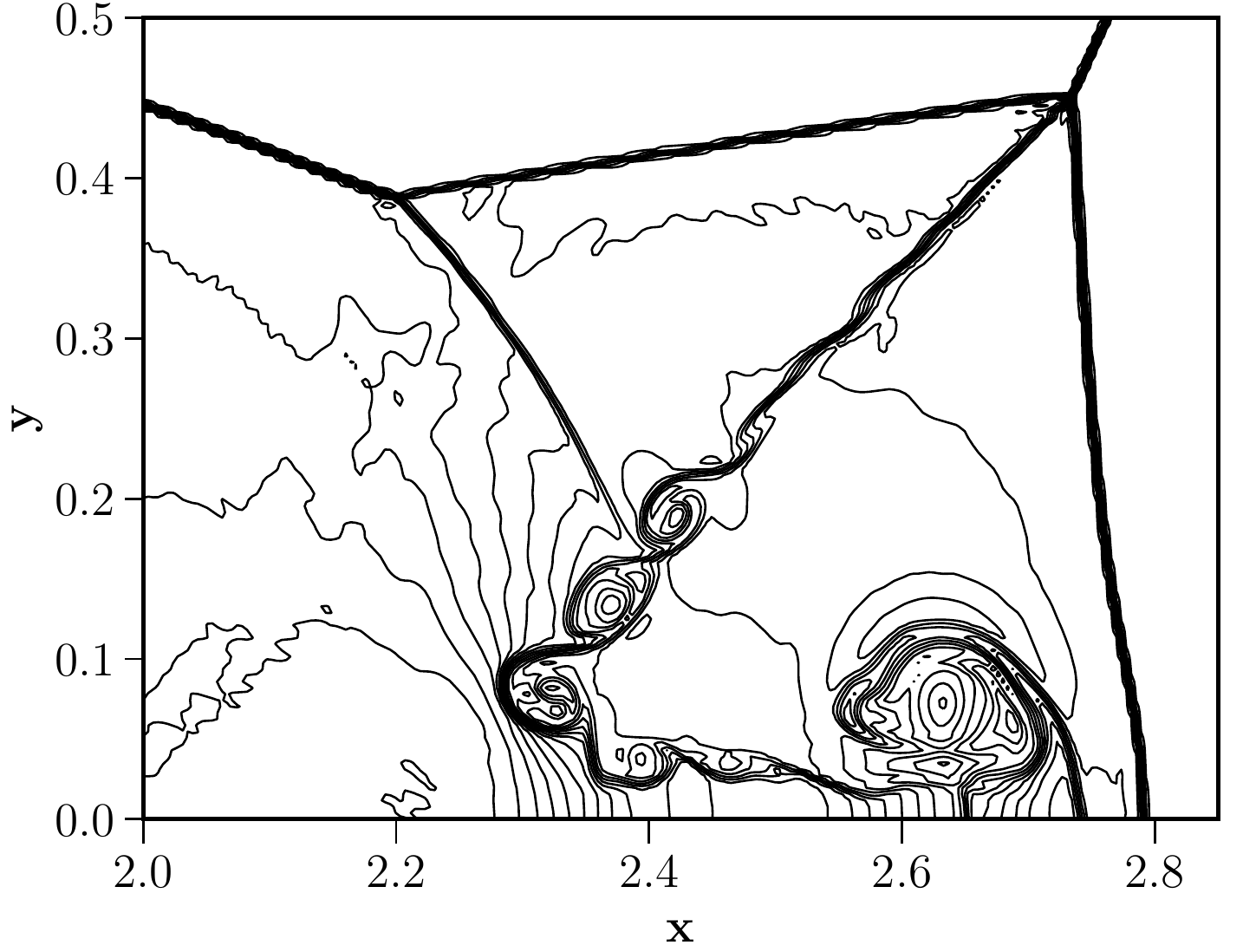}
\label{fig:HOCUS6_DB}}
\subfigure[IG4MP]{\includegraphics[width=0.24\textheight]{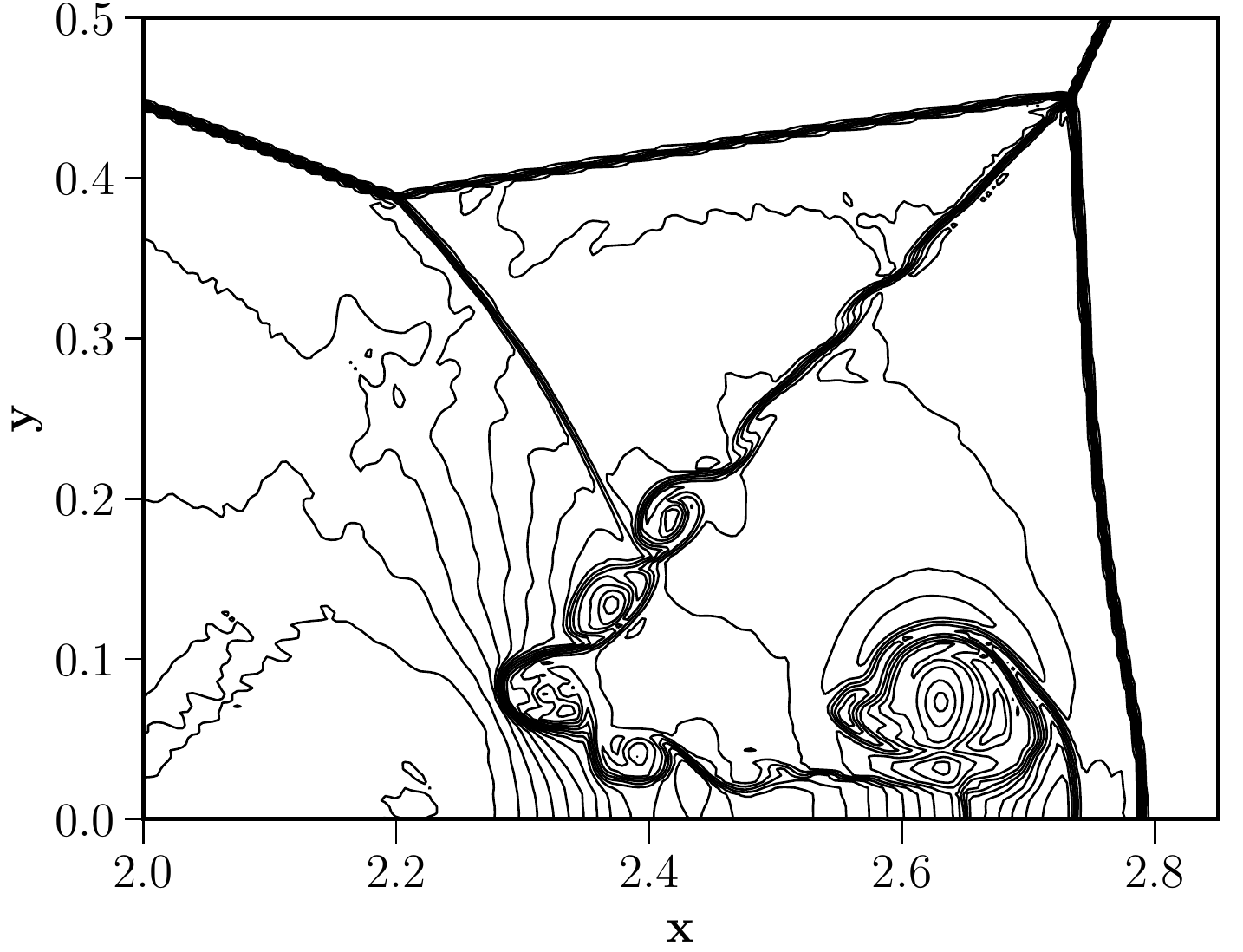}
\label{fig:MP5_DB1}}
\caption{Density contours in the blown-up region around the Mach stem for Example \ref{ex:dmr} on a grid size of 768 $\times$ 256.}
\label{fig_doublemach}
\end{figure}
\begin{figure}[H]
\centering\offinterlineskip
\subfigure[MP5]{\includegraphics[width=0.8\textwidth]{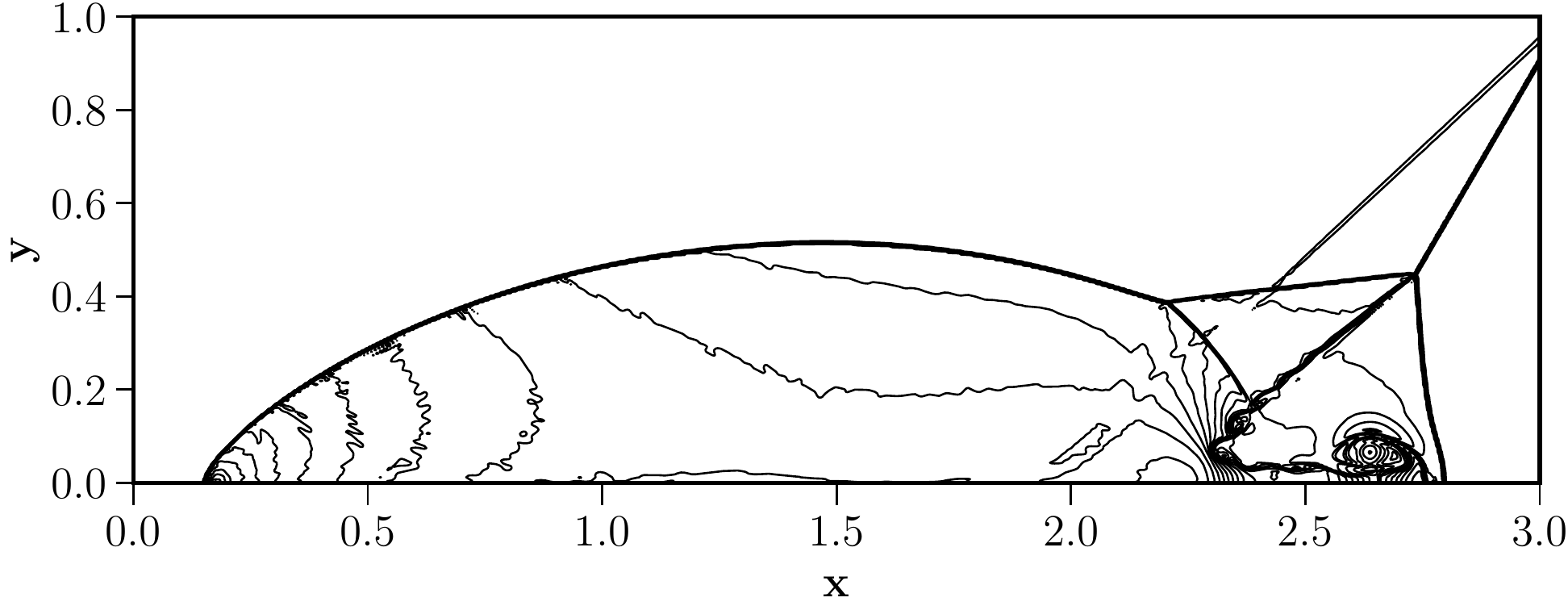}
\label{fig:MP5_DB12}}
\subfigure[IG6MP]{\includegraphics[width=0.8\textwidth]{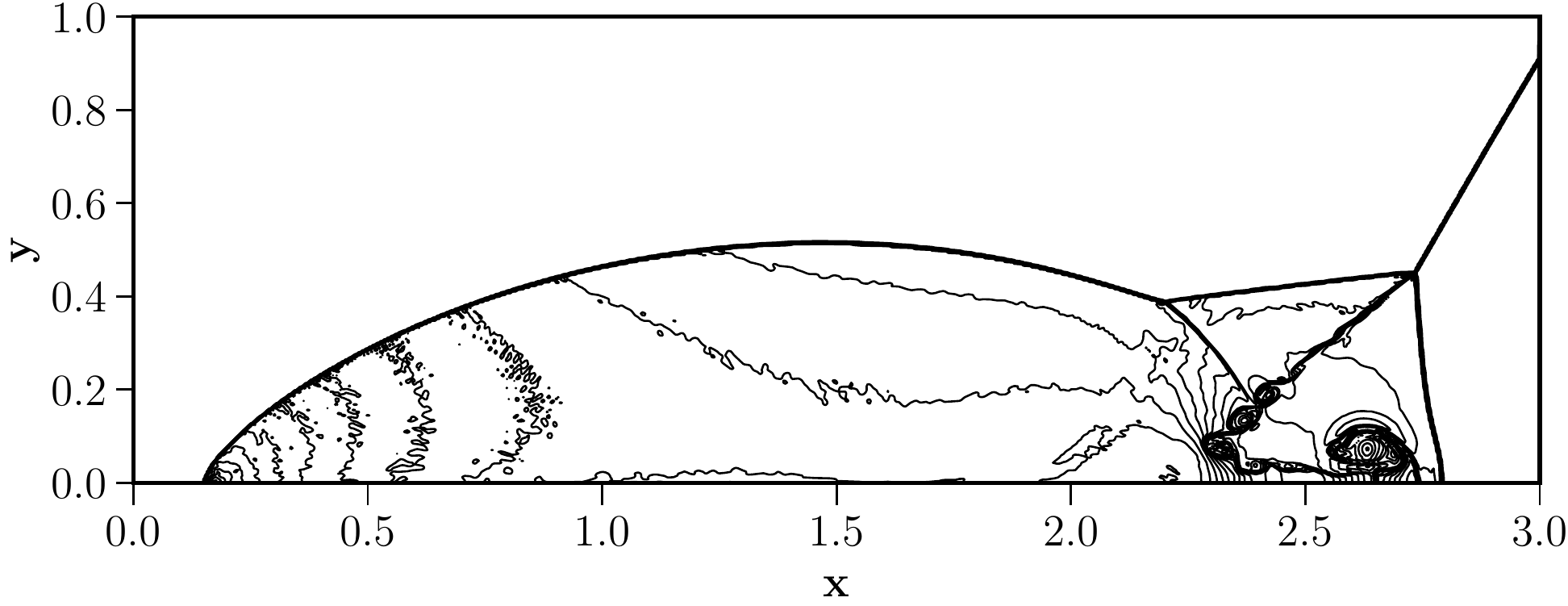}
\label{fig:IG6MP_DB12}}
\subfigure[IG4MP]{\includegraphics[width=0.8\textwidth]{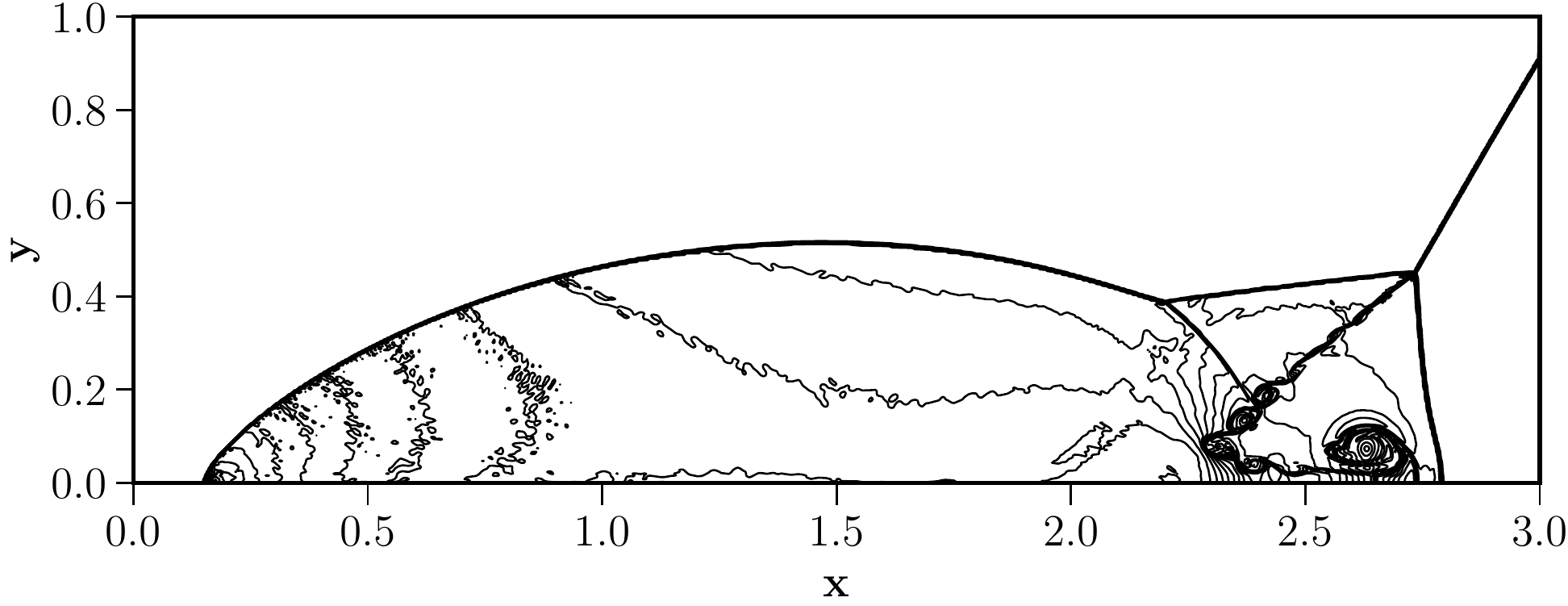}
\label{fig:IG4MP_DB12}}
\caption{Forty equi-spaced density contours of Example \ref{ex:dmr} for different schemes.}
\label{fig_doublemach-full}
\end{figure}


\begin{example}\label{ex:TGV}{Inviscid Taylor-Green Vortex}
\end{example}

Next, we consider the three-dimensional inviscid Taylor-Green vortex problem, with initial conditions given by equation \eqref{itgv} and are applied to the domain of size $x,y,z \in [0,2\pi)$. Periodic boundary conditions are applied for all boundaries. The ratio of the specific heats of the gas is taken as $\gamma=5/3$. The simulations are performed for a time $t=10$ on a grid size of $64 \times 64 \times 64$. 

\begin{equation}\label{itgv}
\begin{pmatrix}
\rho \\
u \\
v \\
w \\
p \\
\end{pmatrix}
=
\begin{pmatrix}
1 \\
\sin{x} \cos{y} \cos{z} \\
-\cos{x} \sin{y} \cos{z} \\
0 \\
100 + \frac{\left( \cos{(2z)} + 2 \right) \left( \cos{(2x)} + \cos{(2y)} \right) - 2}{16}
\end{pmatrix}
\end{equation}

This flow problem is essentially incompressible as the mean pressure is chosen to be very large. The Taylor-Green vortex is the simplest problem for analysing the nonlinear transfer of kinetic energy among the different scales of the flow. It contains several physical processes that are key to understanding turbulence. The vortices in the initial flow stretch and produce smaller-scale features with time. This problem can be used as a test to examine the scale-separation ability of different schemes to under-resolved flow. We compare the ability of different schemes to preserve kinetic energy and also the growth of enstrophy in time, i.e., the sum of vorticity of all the vortex structures, indicating the schemes ability to preserve as many structures as possible. The enstrophy can be described as the integral of the square of the vorticity that can be computed as the integral of the magnitude of vorticity, $\overrightarrow{\omega}$, over the whole domain, 

\begin{equation}
Enstrophy = \sum_{cells}{} \left \| \overrightarrow{\omega} \right \|
\end{equation}
Fig. \ref{fig_TGV} shows the normalised kinetic energy and normalised enstrophy with respect to the initial values for different schemes. We can observe that the WENOZ is the worst among all the schemes in preserving the kinetic energy over time. The IGMP schemes preserve the kinetic energy better than the MP5 scheme, with IG4MP slightly better than IG6MP, which shows the BVD algorithm's effectiveness. We also conducted simulations using the unlimited linear schemes IG4, IG6 and C5 to show the dissipation properties of the novel reconstruction schemes proposed in the paper. It can be seen that the kinetic energy is well preserved by the IG4 scheme in comparison with C5 and IG6, which corroborates our theoretical observations from the Fourier analysis in Fig. \ref{fig_disp}. Another advantage of the present IG schemes is the re-use of the velocity gradients in computing the enstrophy, which is the same used for the reconstruction of the interface states given by Equations \ref{eqn:ddx}, for IG4 and IG6 schemes, respectively. From the enstrophy plot, it is observed that the present schemes outperform the original MP5 and the WENOZ schemes significantly. The IGMP schemes also outperform the HOCUS5 scheme, presented in Ref. \cite{chamarthi2021high}, a combination of MP5 and C5 schemes. Even though the BVD algorithm effectively captures discontinuities, there is still a significant difference between the kinetic energy and enstrophy values computed by the linear IG and nonlinear IGMP schemes, which can be improved in the future.

\begin{figure}[H]
\centering
\subfigure[Kinetic energy]{\includegraphics[width=0.48\textwidth]{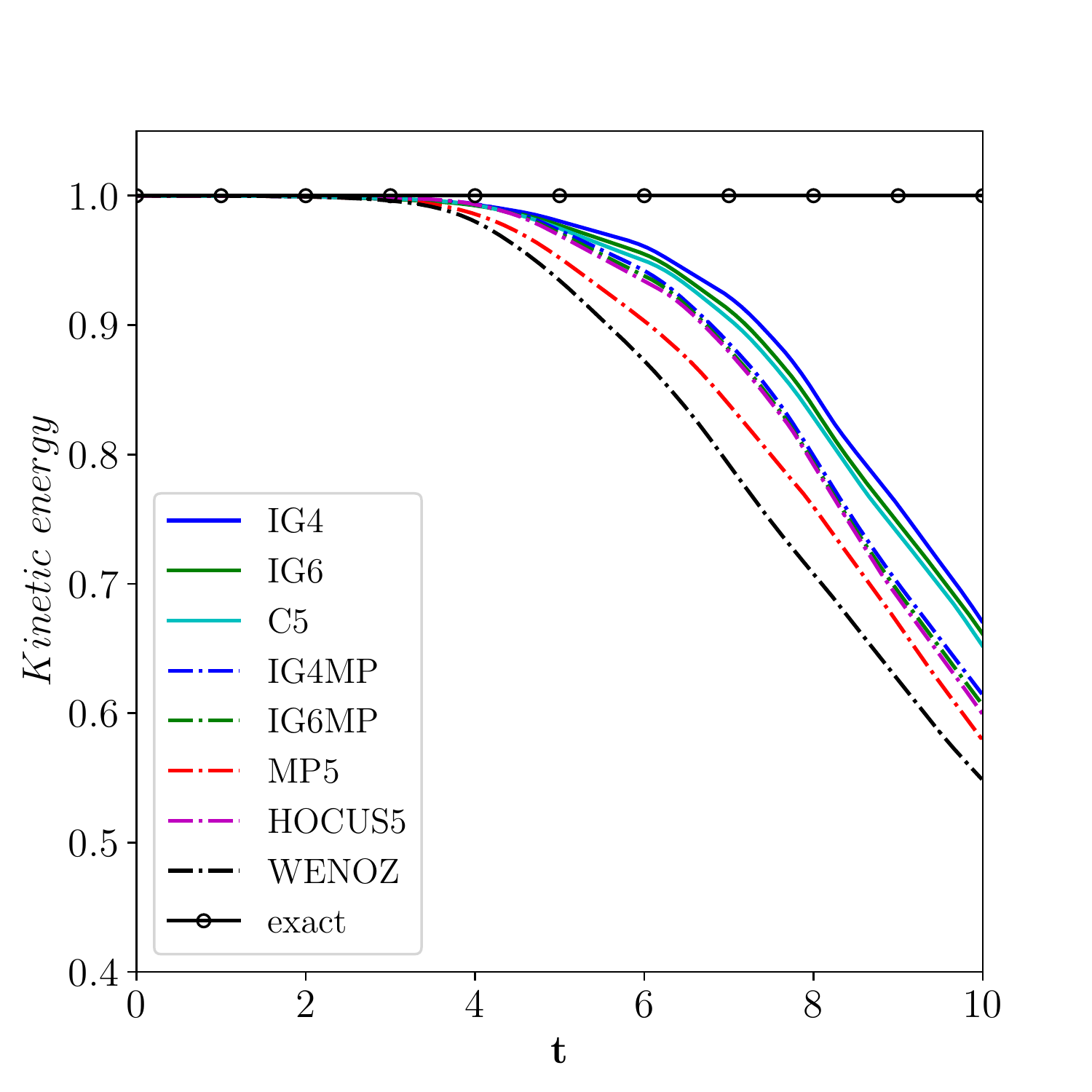}
\label{fig:TGV_KE}}
\subfigure[Enstrophy]{\includegraphics[width=0.48\textwidth]{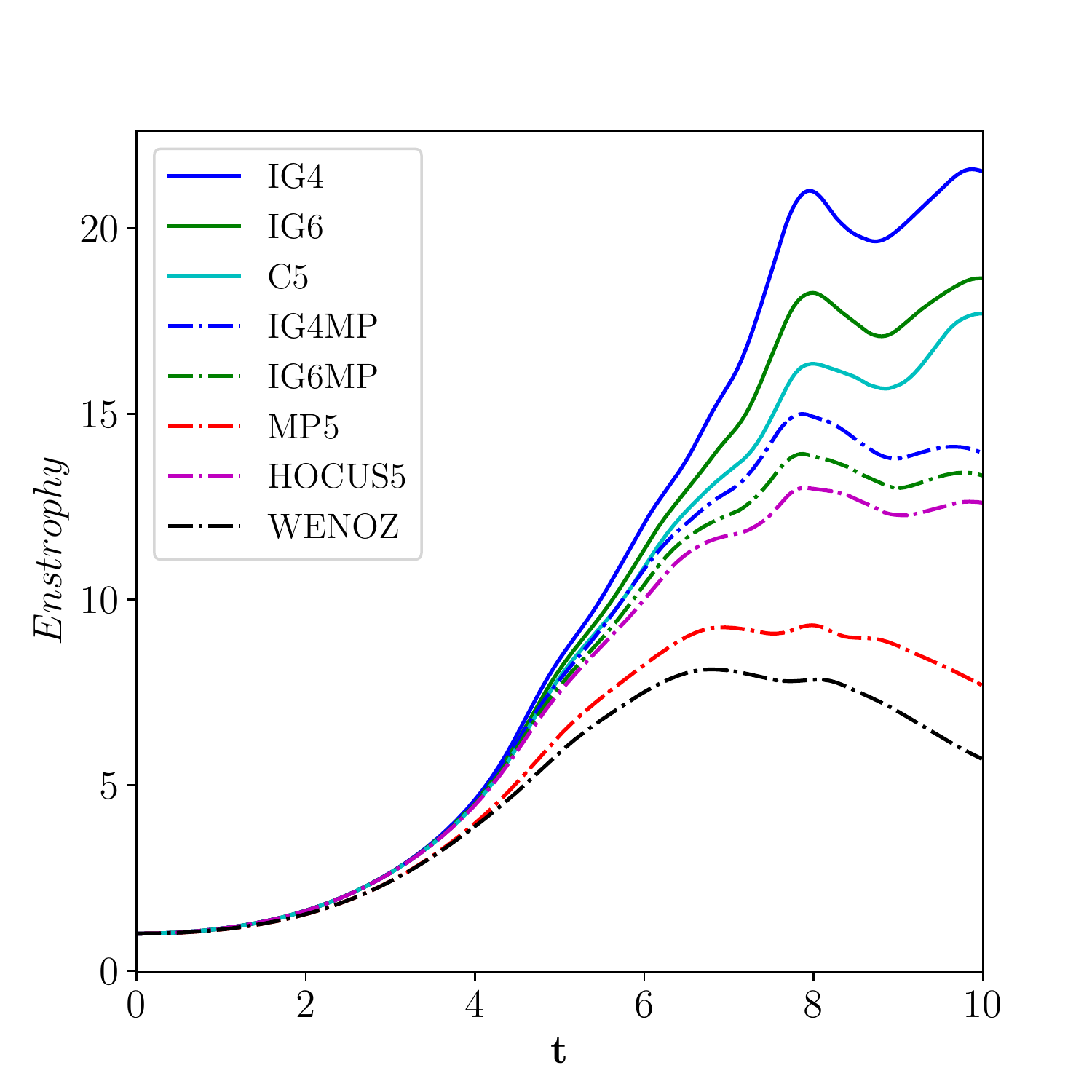}
\label{fig:TGV_ens}}
\caption{Normalised kinetic energy and enstrophy for different schemes presented in Example \ref{ex:TGV} on grid size of $64^3$. Solid line with circles: exact solution; solid blue line: IG4; solid green line: IG6; solid cyan line: C5; dashed blue line: IG4MP; dashed green line: IG6MP; dashed red line: MP5; dashed magenta line: HOCUS5; dashed black line: WENOZ.}
\label{fig_TGV}
\end{figure}

\begin{example}\label{euler-accuracy}{Accuracy of the proposed schemes}
\end{example}
First, we present the order of accuracy (OOA) of the new schemes by convecting the initial profile given by equation \eqref{accu-euler} in the domain $x,y \in [-1, 1]$. The solution is obtained at time $t=2$, and the time-step is varied as a function of the grid size as $\Delta t = \text{CFL} \Delta x^{2.0}$. 
\begin{align}\label{accu-euler}
(\rho,u,v,p)= (1+0.5 \sin(x+y),\ \ 1.0,\ \, 1.0 \ \ 1.0)
\end{align}

The $L_2$ norm of the error between the exact and the obtained solution is used to compute the OOA. The OOA's are presented in Table \ref{tab:ooa} indicates that the new schemes are fourth-order accurate in space. Hu et al. \cite{hu2012dispersion} optimized the linear schemes for favourable spectral properties satisfying the dispersion–dissipation relation and but such optimization lead to order degeneration. Even though the proposed schemes are only fourth-order accurate, they have superior dispersion and dissipation properties which is evident from the above test cases.

\begin{table}[H]
 \centering
  \footnotesize
    \begin{tabular}{crcrcrc}
    \hline
      N    & MP5   & OOA & \multicolumn{1}{c}{IG6MP} & OOA & \multicolumn{1}{c}{IG4MP} & OOA \\
          \hline
    10    & 6.79E-03 & -     & \multicolumn{1}{c}{5.98E-04} & -     & \multicolumn{1}{c}{4.65E-04} & - \\
          \hline
    20    & 2.24E-04 & 4.92  & \multicolumn{1}{c}{4.59E-05} & 3.71  & \multicolumn{1}{c}{4.37E-05} & 3.41 \\
          \hline
    40    & 7.06E-06 & 4.98  & \multicolumn{1}{c}{2.54E-06} & 4.18  & \multicolumn{1}{c}{2.30E-06} & 4.25 \\
          \hline
    80    & 2.21E-07 & 5.00  &  \multicolumn{1}{c}{1.77E-07} & 3.84  &  \multicolumn{1}{c}{1.74E-07} & 3.72 \\
                    \hline
    \end{tabular}%
 \caption{$L_2$ errors and numerical order of accuracy for test case in Example \ref{euler-accuracy}. $N$ is the number of cells in the domain}
 \label{tab:ooa}%
\end{table}%


\subsection{Two-dimensional test cases for Navier-Stokes equations}
\begin{example}\label{ex:ldc}{Lid driven cavity}
\end{example}

Next, we validate the solution of the classical lid-driven cavity problem with the present schemes.  Even though the test case is incompressible flow, it serves as a validation for the proposed schemes. The two-dimensional square domain of the problem is set to unit size. The no-slip boundary conditions are applied to all boundaries, and additionally, the top wall moves with a fixed velocity in the positive x-direction ($u_{top}$ =1). It is challenging to achieve a steady-state using unsteady compressible flow methods at the low-Mach number. The steady-state convergence was considered to be achieved when the calculated velocity fields remain unchanged for 400 consecutive time steps \cite{houim2011low}. The steady-state solution is obtained for two Reynolds numbers ($Re = \frac{\rho u_{top}}{\mu}$) of 400, 1000, on a grid size of $128 \times 128$. The streamline plots and the contour plots of the x-velocity for different Reynolds numbers are presented in Fig. \ref{fig_ldc_a} and for the IG4MP scheme. Also, as shown in Fig. \ref{fig_ldc_b} the normalised velocity profiles along the mid-sections of the domain are observed to be in good agreement with Ghia et al. \cite{ghia1982high} for both the IG4MP and IG6MP schemes.

\begin{figure}[H]
\begin{onehalfspacing}
\centering\offinterlineskip
\subfigure[$Re=400$]{\includegraphics[width=0.4\textwidth]{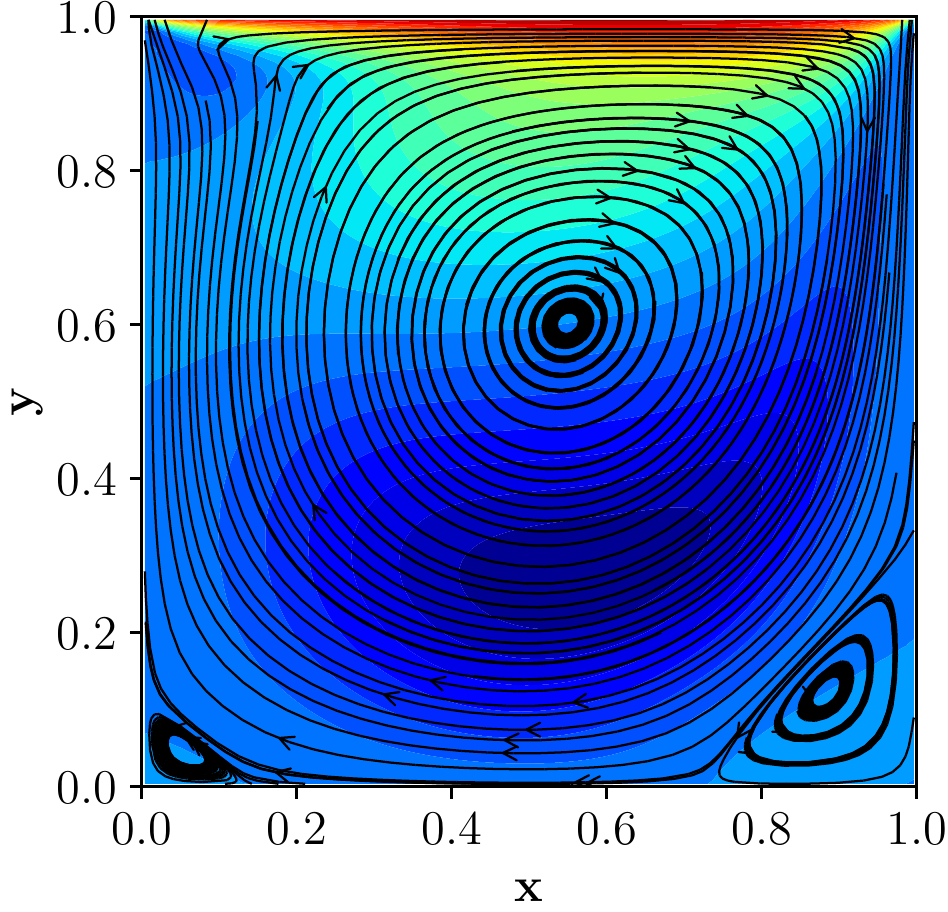}
\label{fig:ldc1}}
\subfigure[$Re=1000$]{\includegraphics[width=0.4\textwidth]{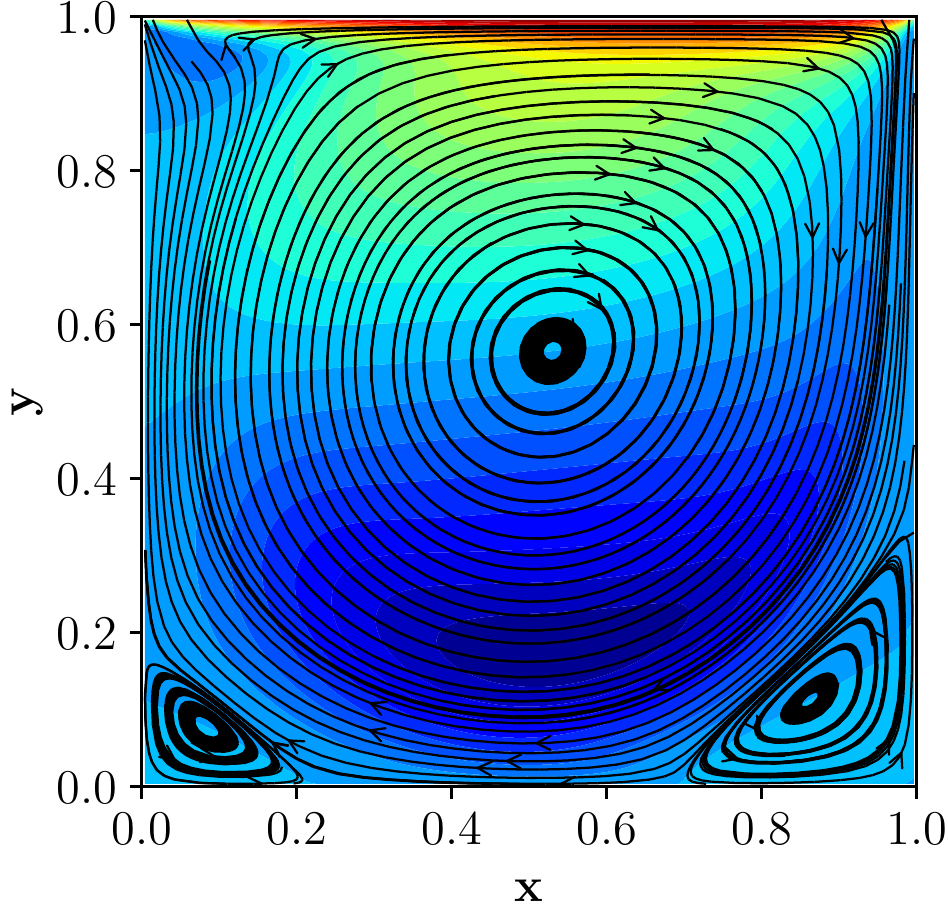}
\label{fig:ldc2}}
\caption{Streamline and x-velocity contour plots of the lid driven cavity problem for different Reynolds numbers using the IG4MP scheme.}
\label{fig_ldc_a}
\end{onehalfspacing}
\end{figure}

\begin{figure}[H]
\centering\offinterlineskip
\subfigure[Re=400]{\includegraphics[width=0.4\textwidth]{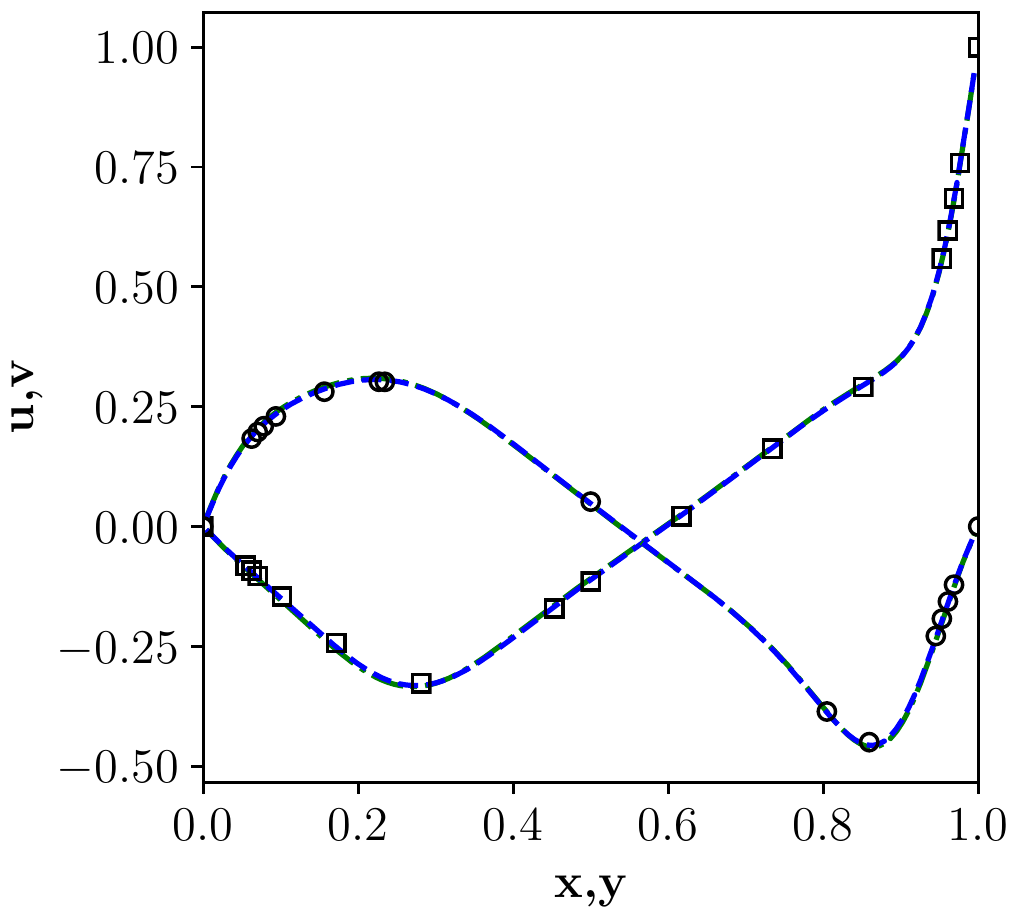}
\label{fig:ldc3}}
\subfigure[Re=1000]{\includegraphics[width=0.4\textwidth]{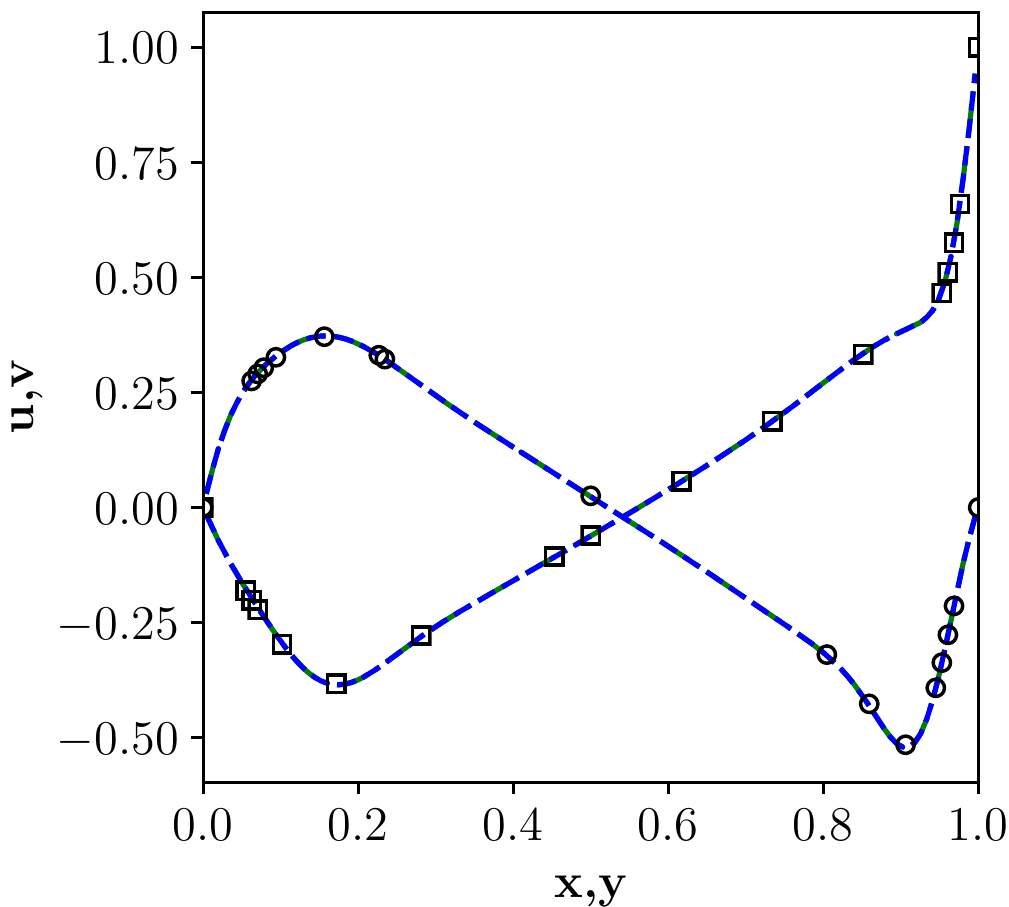}
\label{fig:ldc4}}
\caption{Normalised velocity profiles along the mid-sections of the domain, for different Reynolds numbers, for Example \ref{ex:ldc} on a grid size of 128 $\times$ 128. Black circles: u velocity for $Re$ = 400 and 100; black squares: v velocity for $Re$ = 400 and 1000; green dashed line: IG6MP; blue dotted line: IG4MP.}
\label{fig_ldc_b}
\end{figure}


\begin{example}\label{ex:vs}{Viscous Shock tube}
\end{example}

This last test case for the single component flows is the viscous shock-tube problem in a square shock tube of domain of unit length and half unit height and insulated walls taken from Ref.\cite{daru2009numerical}. In this problem, the propagation of the incident shock wave and contact discontinuity leads to development of a thin boundary layer at the bottom wall. After its reflection at the right wall, the shock wave interacts with this boundary layer. These interactions result in a complex vortex system, separation region, and a typical lambda-shaped shock pattern. The flow contains various length scale structures, making it an ideal test case for testing high-resolution schemes. The initial conditions are:
\begin{equation}\label{vst}
\begin{aligned}
\left( {\rho , u,v, p} \right) = \left\{ \begin{array}{l}
\left( {120, 0 ,0,120/\gamma } \right),  \quad 0 < x < 0.5,\\
\left( {1.2 , 0 ,0, 1.2/\gamma } \right),  \quad 0.5 \le x < 1,
\end{array} \right.
\end{aligned}
\end{equation}
A shock wave is initialized at the mid-location in the domain, which interacts with a developing boundary layer as reflected from the slip wall at the right boundary. Also, the contact discontinuity and the reflected shock wave interact with each other. The domain for this test case is taken as $x \in [0,1], y \in [0,0.5]$. The left boundary is set as an outflow boundary. The initial conditions are given by equation \eqref{vst}, with the ratio of specific heats of $\gamma = 7/5$. The flow is simulated for time $t=1$, keeping the Mach number of the shock wave at 2.37. The problem is solved for two different Reynolds numbers $Re=500$ and $Re=1000$ on a grid size of $750 \times 375$ for the $Re=500$ case, and $1280 \times 640$, for the $Re=1000$ case, respectively. For the MP5 scheme the gradients of the velocity components and temperature are computed using the generic fourth order explicit formula as in Ref. \cite{Fu2016}, given by Equation (\ref{eqn:4e}) and denoted as MP5 - 4E , and for the IGMP schemes the gradients are as explained in Section \ref{sec-3.2}. The observations from the simulations are as follows:

\begin{itemize}
\item For the $Re=500$ case, the density profiles obtained by various schemes are shown in Fig. \ref{fig_viscousshock_500}. The numerical results obtained by IGMP schemes are almost the same and also similar to the converged results in Ref.\cite{daru2009numerical} despite using coarse mesh. The result obtained by the MP5 scheme shows a significant difference from the IGMP scheme considering the height and shape of the primary vortex. Density distributions along the bottom wall of the shock tube, obtained at $t$ = 1 for IGMP schemes, is compared with the converged results of Daru and Tenaud \cite{daru2009numerical} (See their Fig. 2 ). It can be seen that the present method gives similar results on a much coarser grid resolution than theirs. 

\begin{figure}[H]
\begin{onehalfspacing}
\centering\offinterlineskip
\subfigure[MP5 - 4E]{\includegraphics[width=0.48\textwidth]{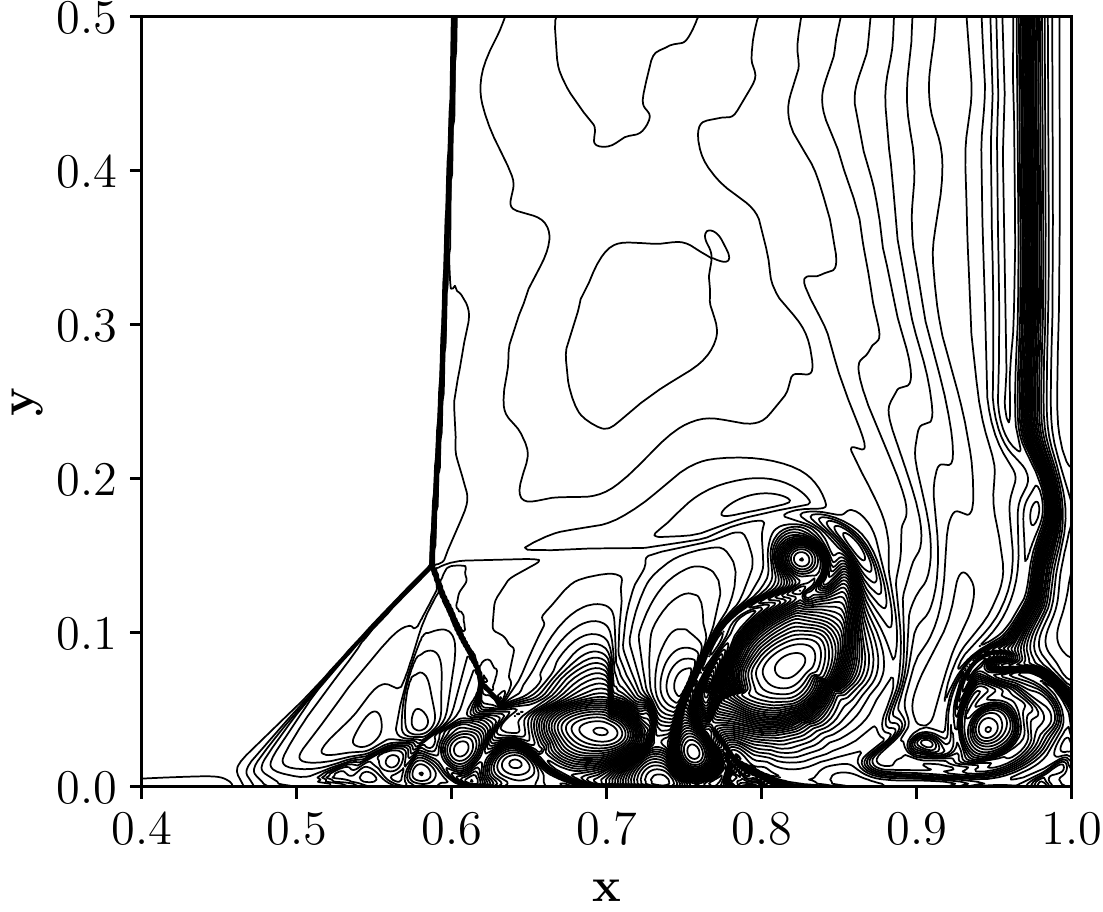}
\label{fig:MP5_VS4E}}
\subfigure[IG6MP]{\includegraphics[width=0.48\textwidth]{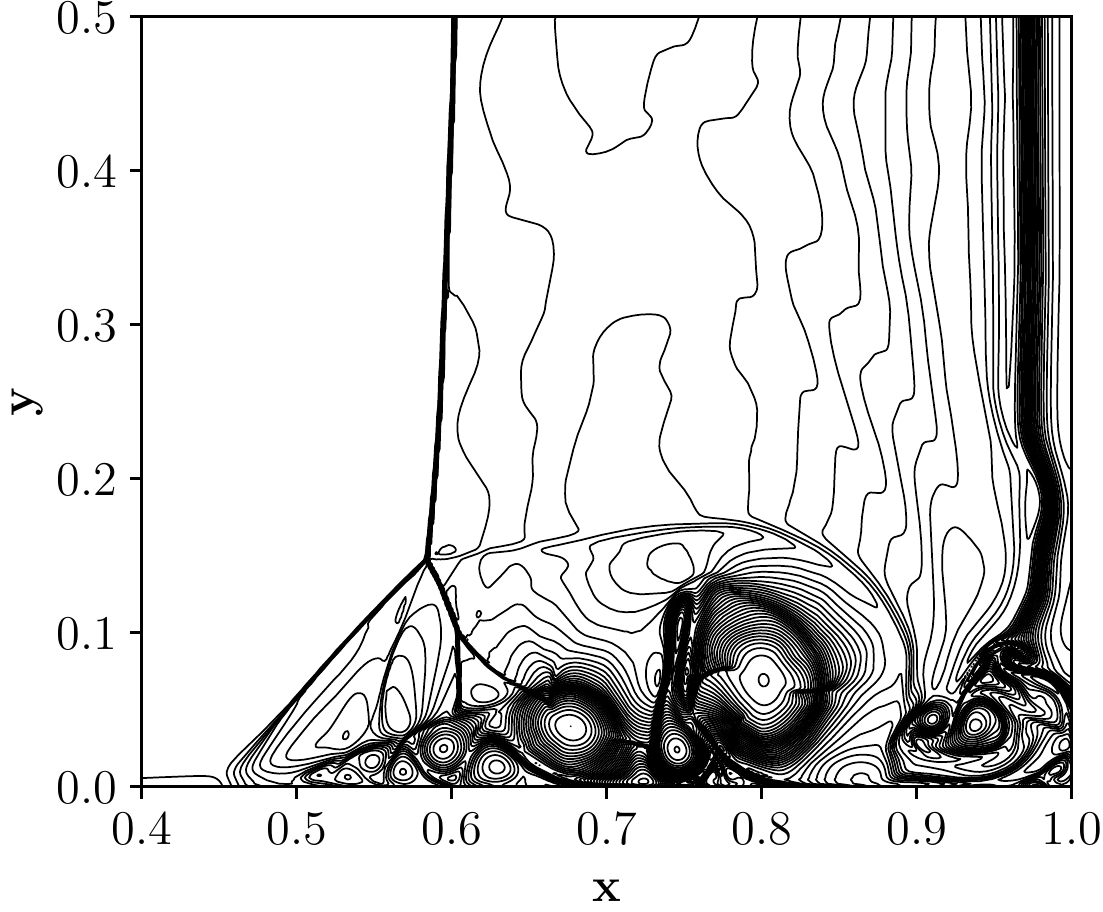}
\label{fig:VST_IG6MP}}
\subfigure[IG4MP]{\includegraphics[width=0.48\textwidth]{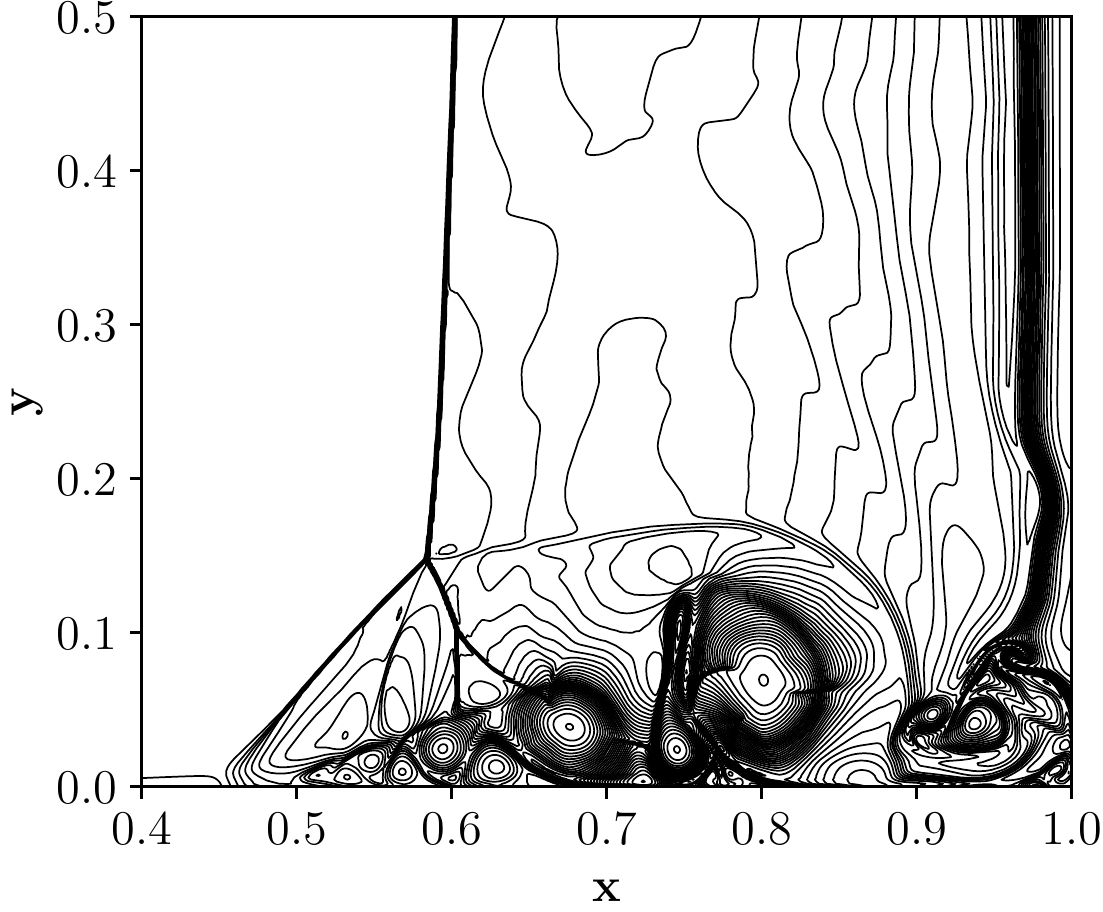}
\label{fig:VST_IG4MP}}
\subfigure[Density profile along the wall, $Re$= 500.]{\includegraphics[width=0.48\textwidth]{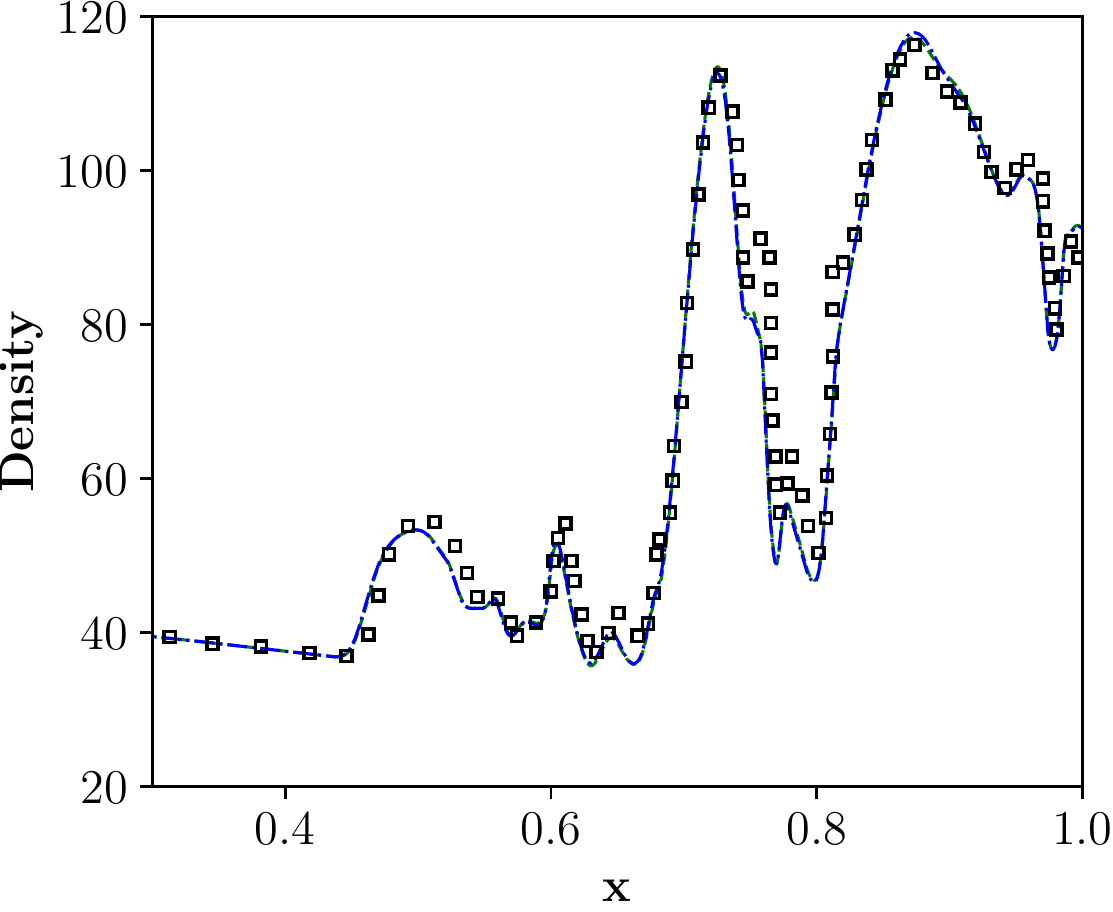}
\label{fig:HOCUS6_VS12}}
\caption{Density contours for different schemes for the Example \ref{ex:vs} and $Re=500$ for various schemes. Density along the wall for both the proposed schemes are compared to the reference solution represented by black squares. The green line represents the IG6MP scheme, and the blue line represents the IG4MP scheme, respectively.}
\label{fig_viscousshock_500}
\end{onehalfspacing}
\end{figure}

\item Fig. \ref{fig_viscousshock_1000} shows the density profiles for the $Re$ = 1000. When the Reynolds number is increased to 1000, the flow structures are much more complicated. The boundary layer separates at several points, giving rise to the development of very complicated vortex structures and the interactions between vortices and shock waves. Converged results for this test case can be found in Ref.\cite{Fu2016} (see their Fig. 21) and \cite{zhou2018grid} (their Fig. 6e) which are carried out on a grid resolution of 5120 $\times$ 2560 and 5000 $\times$ 2500, respectively. The MP scheme shows a significant difference from the converged results, with the primary vortex being completely distorted.
On the other hand, the present schemes captures the lambda-shaped shocks accurately along with the small scale features. The rotating vortices on the lower right corner fit very well with the converged results along with the shape of the primary vortex, especially for the IG4MP scheme. Density distributions along the bottom wall for $Re$ =1000 agree very well with those of Zhou et al.\cite{zhou2018grid} (their Fig. 7b). These results indicate that the present schemes can compute the multi-scale flows with shock waves in high resolution, even on coarse grids.

\begin{figure}[H]
\begin{onehalfspacing}
\centering\offinterlineskip
\subfigure[MP5 - 4E]{\includegraphics[width=0.48\textwidth]{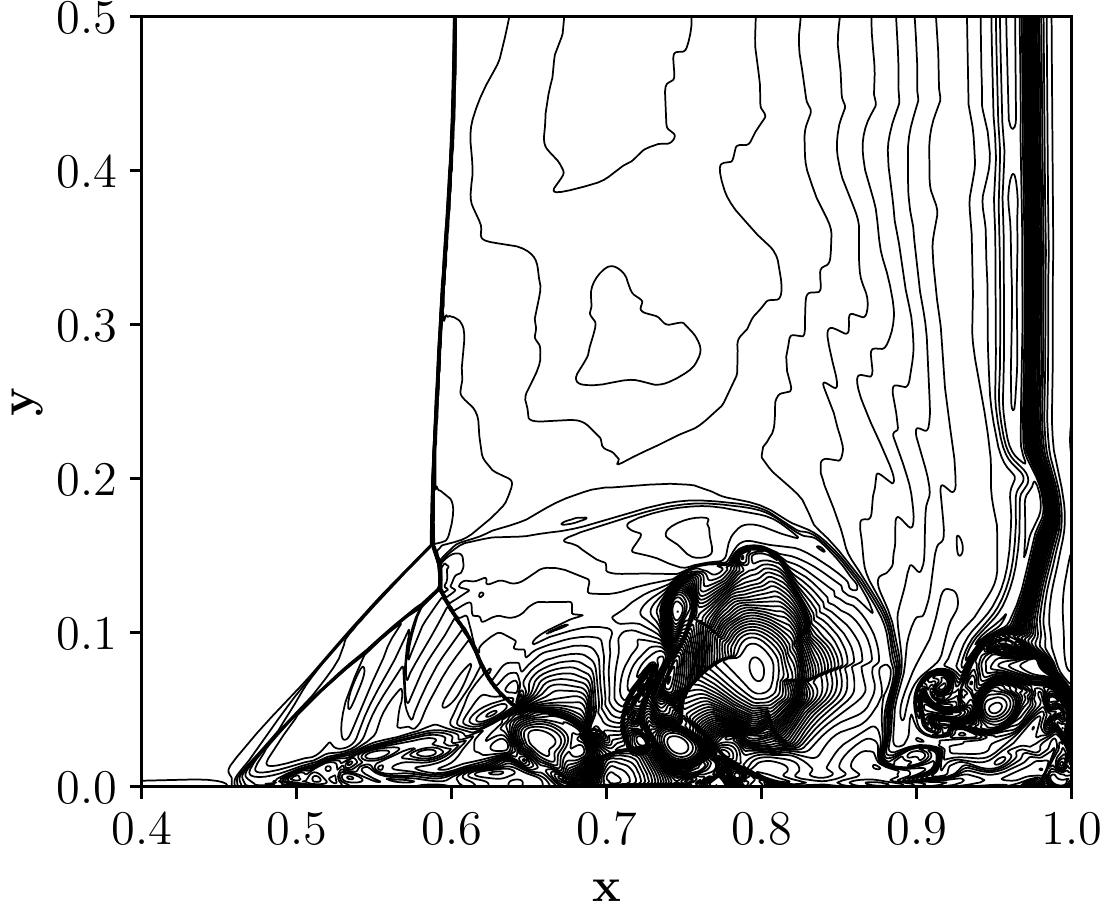}
\label{fig:MP5_VS123}}
\subfigure[IG6MP]{\includegraphics[width=0.48\textwidth]{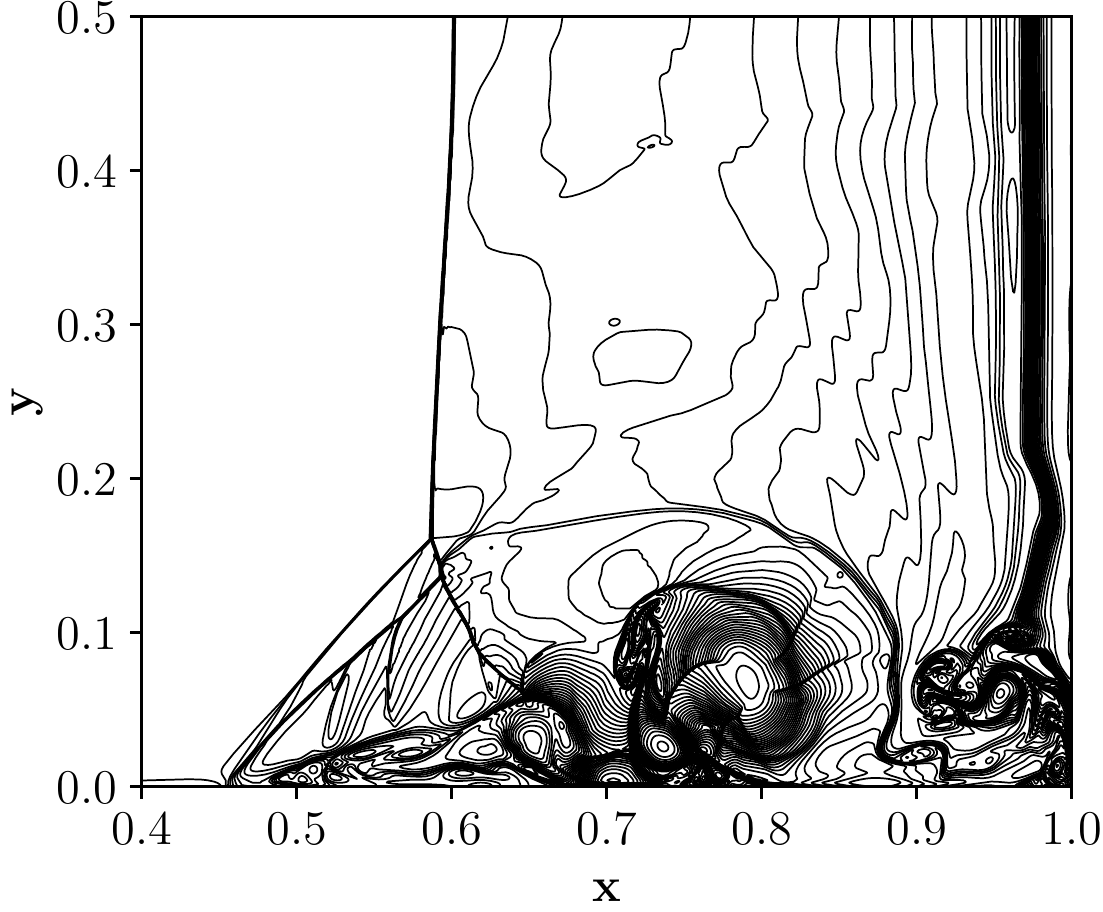}
\label{fig:M}}
\subfigure[IG4MP]{\includegraphics[width=0.48\textwidth]{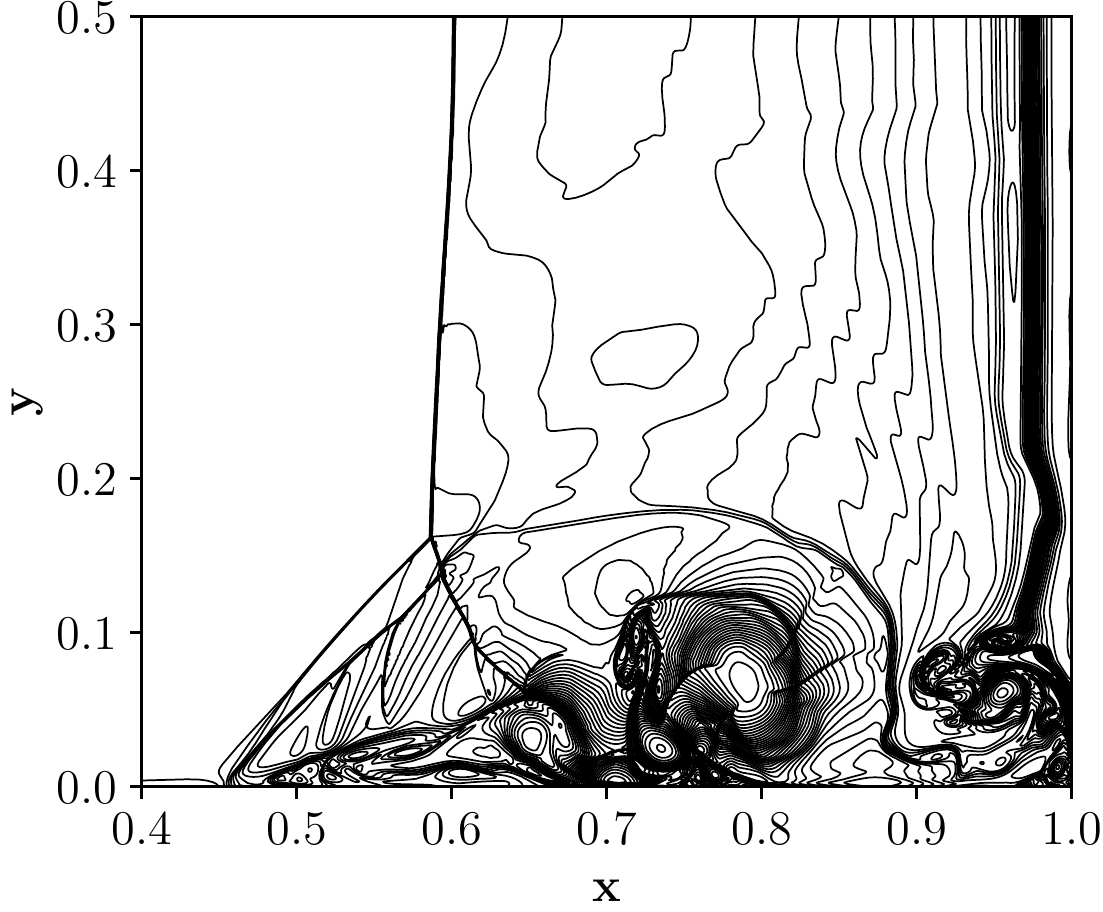}
\label{fig:H}}
\subfigure[Density profile along the wall, $Re$=1000.]{\includegraphics[width=0.48\textwidth]{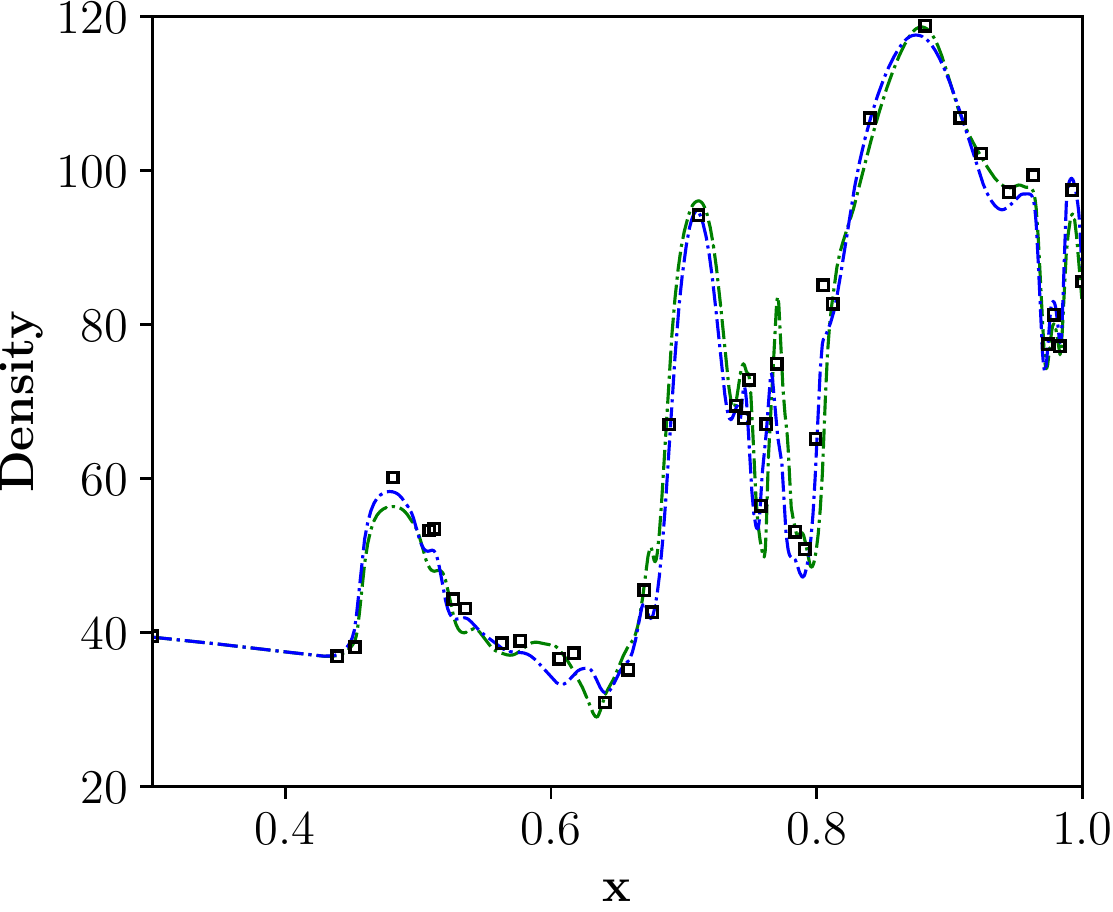}
\label{fig:MP5_VS2}}
\caption{Density contours for different schemes for the Example \ref{ex:vs} and $Re=1000$ for various schemes. Density along the wall for both the proposed schemes are compared to the reference solution represented by black squares. The green line represents the IG6MP scheme, and the blue line represents the IG4MP scheme, respectively.}
\label{fig_viscousshock_1000}
\end{onehalfspacing}
\end{figure}

\item The fourth-order explicit approximation for gradient evaluation in viscous terms for the MP5 scheme is based on \cite{Fu2016}. Daru and Tenaud \cite{daru2009numerical} indicated that the accuracy of viscous terms had negligible on their results. However, we noticed that if the velocity and temperature gradients are computed with implicit gradients, even the results obtained by the MP5 scheme do improve significantly, as shown in Fig \ref{fig_viscousshock_1000_mp}. It is evident that the accuracy of the viscous terms also affects the resolution of the results for this test case.

\begin{figure}[H]
\begin{onehalfspacing}
\centering\offinterlineskip
\subfigure[Re - 500]{\includegraphics[width=0.48\textwidth]{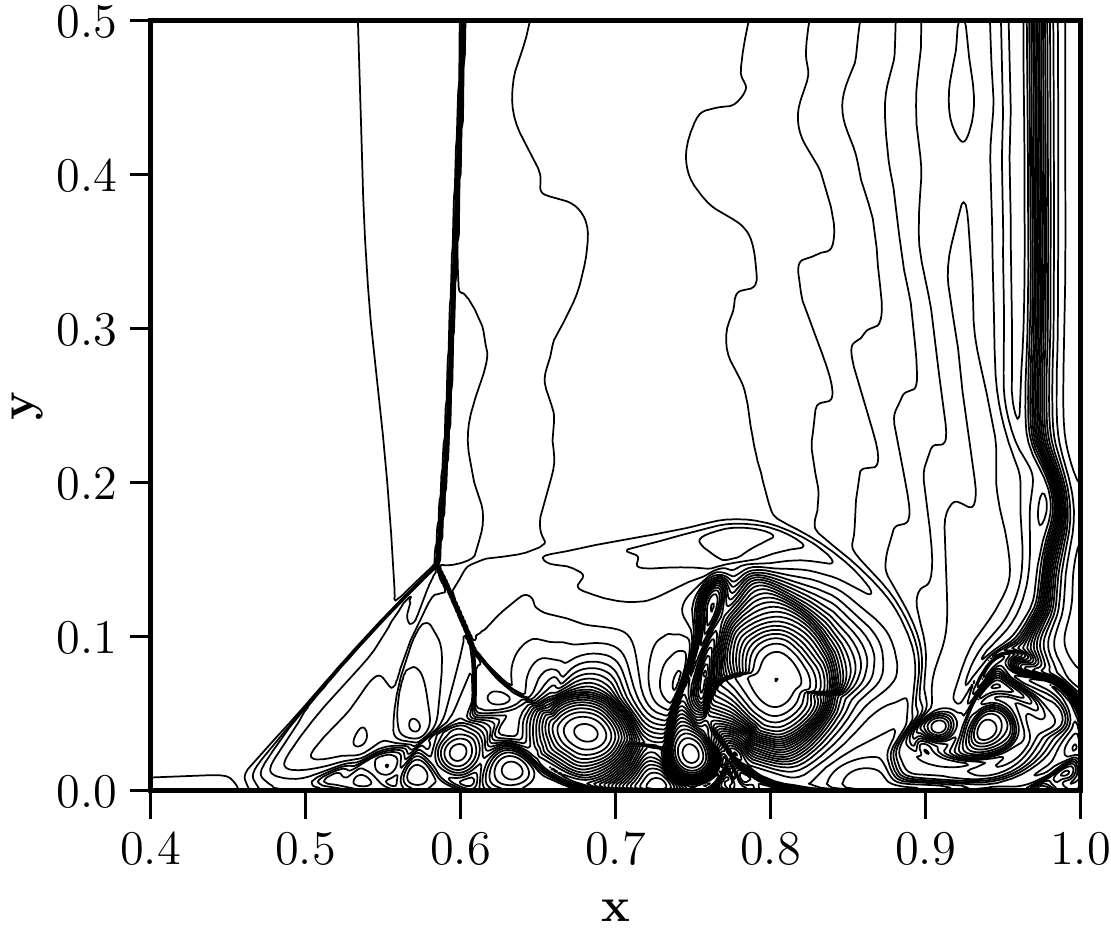}
\label{fig:MP5_4T_500}}
\subfigure[Re - 1000]{\includegraphics[width=0.48\textwidth]{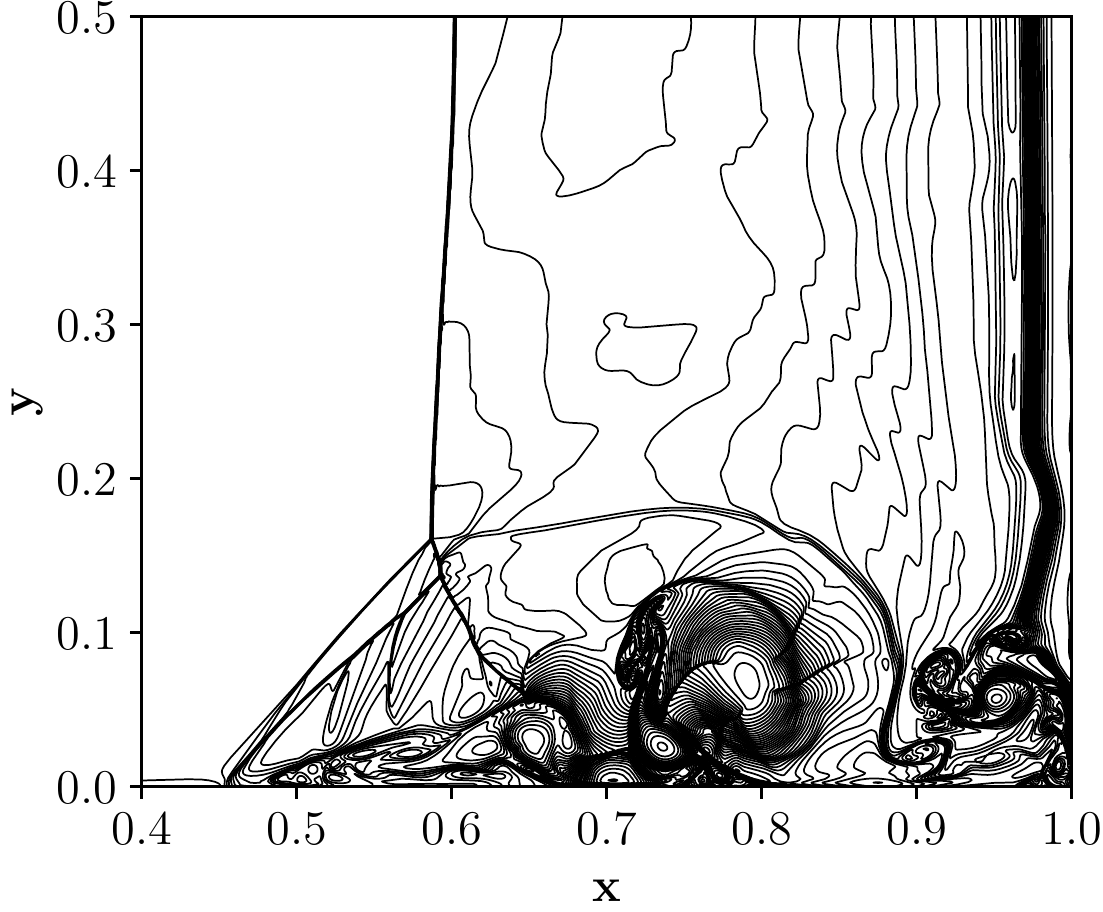}
\label{fig:MP5_4T_1000}}
\caption{Density contours for the Example \ref{ex:vs} for $Re=500$ and $Re=1000$ using implicit gradients for the viscous terms and MP5 scheme for the  convective terms.}
\label{fig_viscousshock_1000_mp}
\end{onehalfspacing}
\end{figure}

\end{itemize}


\subsection{Numerical test cases for multi-species flows}
This section presents the numerical results of the proposed algorithm for one- and two-dimensional inviscid and viscous multi-component flows.
\begin{example}\label{ex:multi-sod}{Multi-species shock tube}
\end{example}

The first one-dimensional test case is the two-fluid modified shock tube of Abgrall and Karni \cite{abgrall2001computations}. The initial conditions of the test case are as follows:
\begin{equation}
\left(\alpha_{1} \rho_{1}, \alpha_{2} \rho_{2}, u, p, \alpha_{1}, \gamma \right)=\left\{\begin{array}{ll}
\left(1 , 0, 0, 1, 1, 1.4 \right) & \text { for } x<0 \\
\left(0,0.125,0,0.1, 0, 1.6 \right) & \text { for }  x \geq 0 .
\end{array}\right.
\end{equation}
Simulations are carried out on a uniformly spaced grid with 100 cells on the spatial domain $-0.5 \leq x \leq 0.5$ with a constant CFL of 0.1 and the final time is $t$=0.2. Fig. \ref{fig_multisod} shows the density, pressure and volume fraction profiles obtained by various schemes. Both the IG6MP and the IG4MP schemes capture the shock wave and the material interface without any oscillations. Solution profile depicted by black triangles in Fig. \ref{fig:multi_sod-den_zoom} is obtained by using the IG6MP scheme without the volume fraction boundedness and positivity of phasic densities condition given by Equation (\ref{eqn:multi-positivity}). 
\begin{figure}[H]
\centering
\subfigure[Global density]{\includegraphics[width=0.4\textwidth]{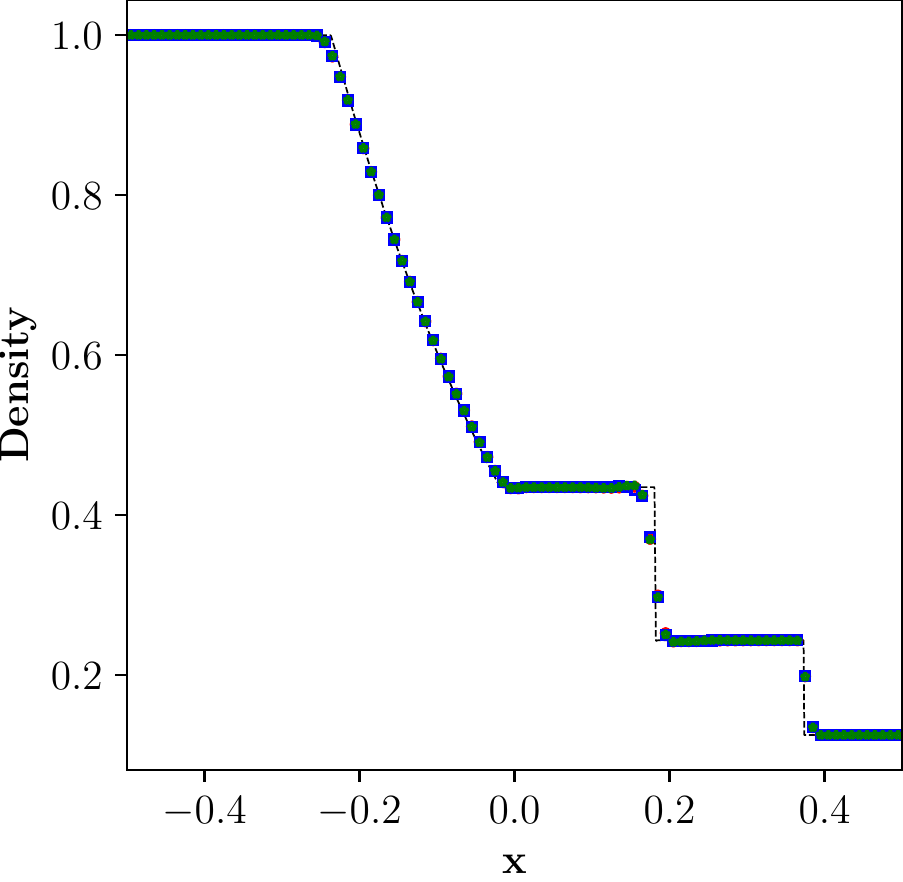}
\label{fig:multi_sod-den}}
\subfigure[Local density ]{\includegraphics[width=0.4\textwidth]{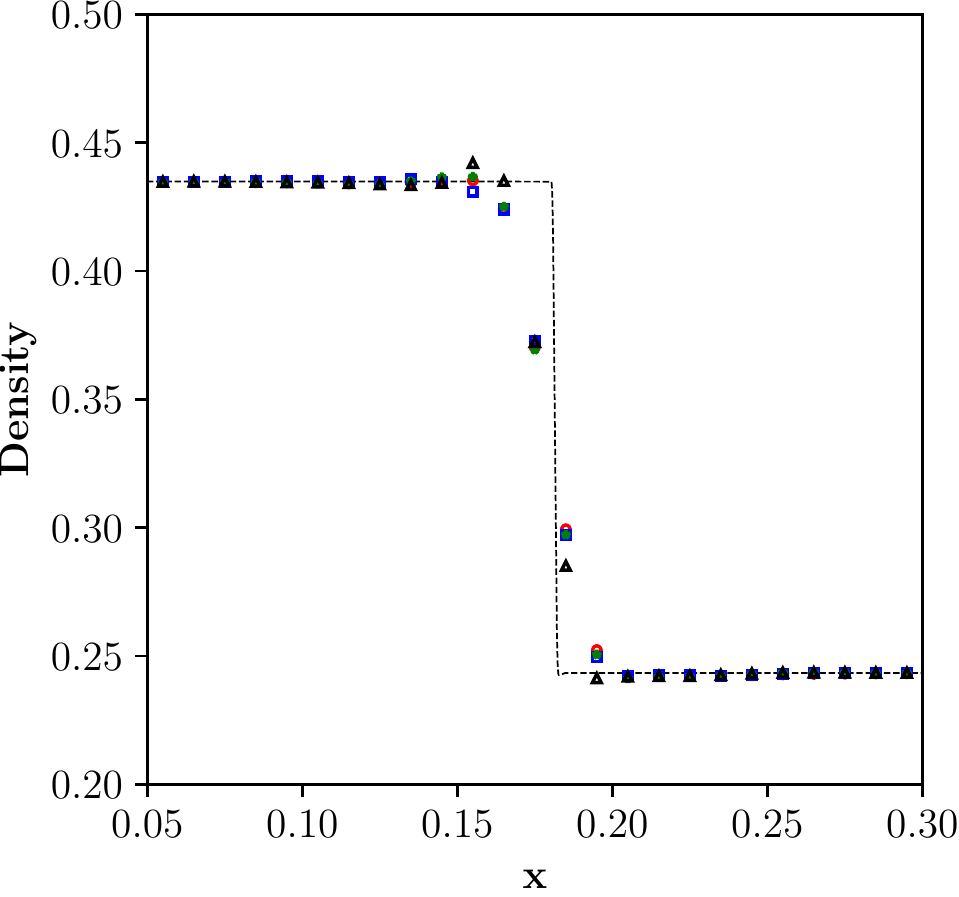}
\label{fig:multi_sod-den_zoom}}
\subfigure[Pressure]{\includegraphics[width=0.4\textwidth]{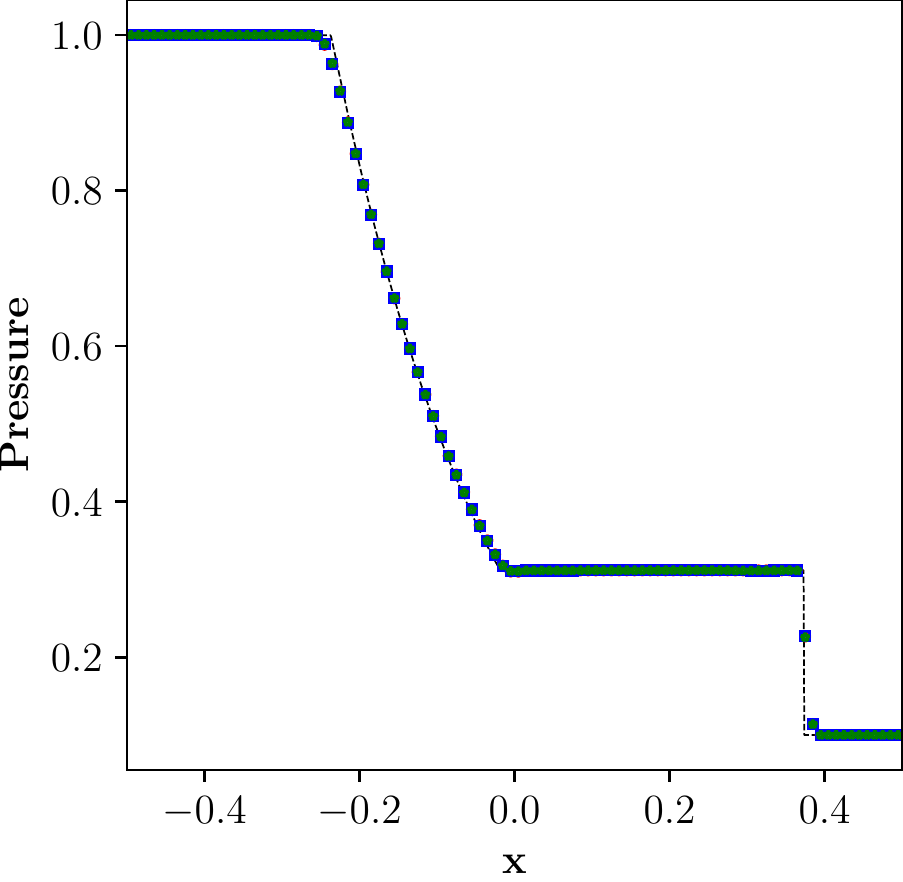}
\label{fig:multi_sod-pres}}
\subfigure[Volume fraction]{\includegraphics[width=0.4\textwidth]{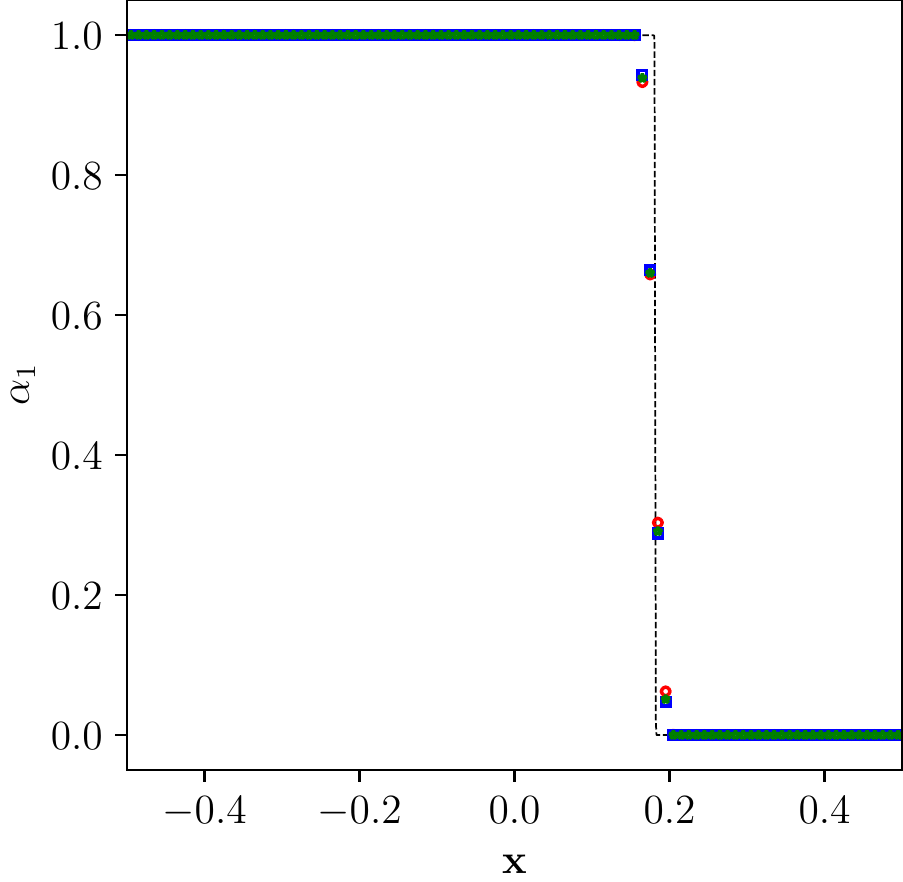}
\label{fig:multi_sod-volf}}
\caption{Numerical solution for multi-species shock tube problem in Example \ref{ex:multi-sod} on a grid size of $N=100$. Dashed line: reference solution; green stars: IG6MP; blue squares: IG4MP; red circles: MP5.}
\label{fig_multisod}
\end{figure}


\begin{example}\label{ex:isolated}{Isolated interface advection}
\end{example}

The next one-dimensional test case is the inviscid isolated interface advection introduced by Johnsen and Colonius \cite{johnsen2006implementation}. This problem is used to test a numerical schemes' capabilities in maintaining the velocity and pressure equilibrium, i.e. preventing spurious oscillations at material interfaces. The initial conditions are as follows with the values of velocity and pressure constant across the material interface,
\begin{equation}
\left(\alpha_{1} \rho_{1}, \alpha_{2} \rho_{2}, u, p, \alpha_{1}, \gamma \right)=\left\{\begin{array}{ll}
\left(10.0 , 0.0,0.5,1/1.4, 1, 1.6\right) & 0.25 \leq x<0.75 \\
\left(1.0 , 0.0,0.5,1/1.4, 1, 1.4\right) & x<0.25 \text { or } x \geq 0.75
\end{array}\right.
\end{equation}

Simulations are performed on a uniformly spaced grid with 50 cells and the spatial domain is $0.0 \leq x \leq 1.0$. Periodic boundary conditions are applied at both boundaries. Simulation is carried out with a constant time step of 0.05 for precisely one period, i.e. until the material interface returns to its original position. Fig. \ref{fig:multi_isolated-den}, shows the density profiles obtained the all the schemes and it can be see that these solutions are free of numerical oscillations. Fig. \ref{fig:multi_isolated-pres}, shows that the errors in pressure are close to zero for all the IGMP schemes in comparison with the exact solution.
\begin{figure}[H]
\centering
\subfigure[Density]{\includegraphics[width=0.4\textwidth]{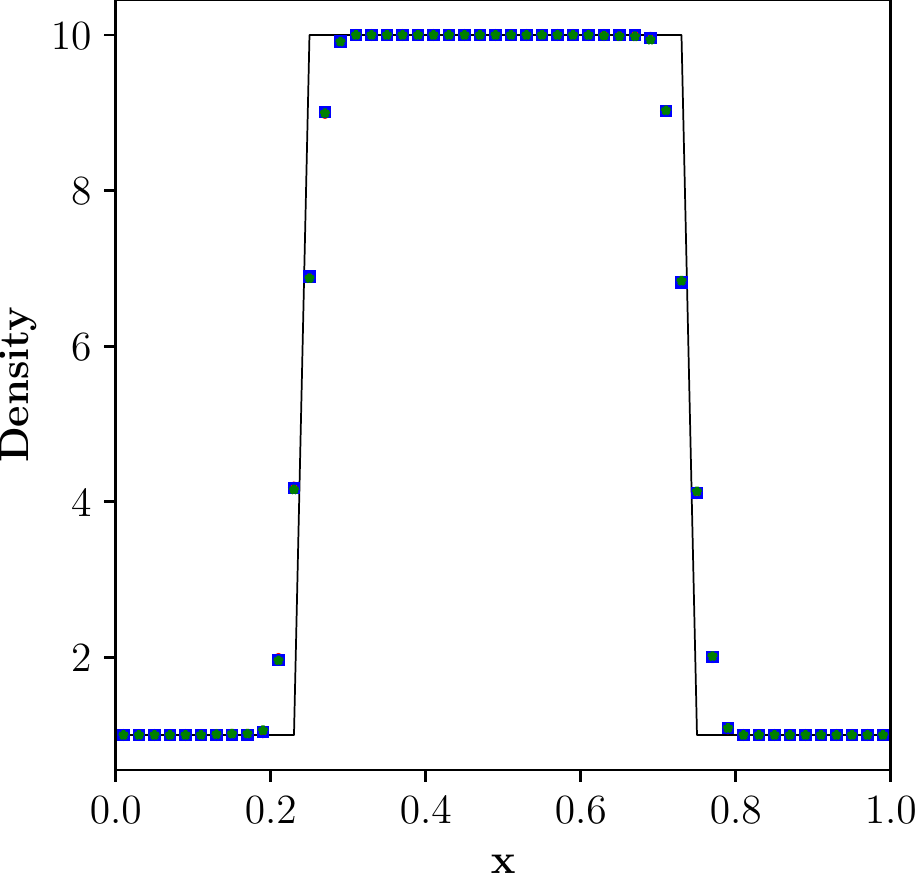}
\label{fig:multi_isolated-den}}
\subfigure[Pressure]{\includegraphics[width=0.4\textwidth]{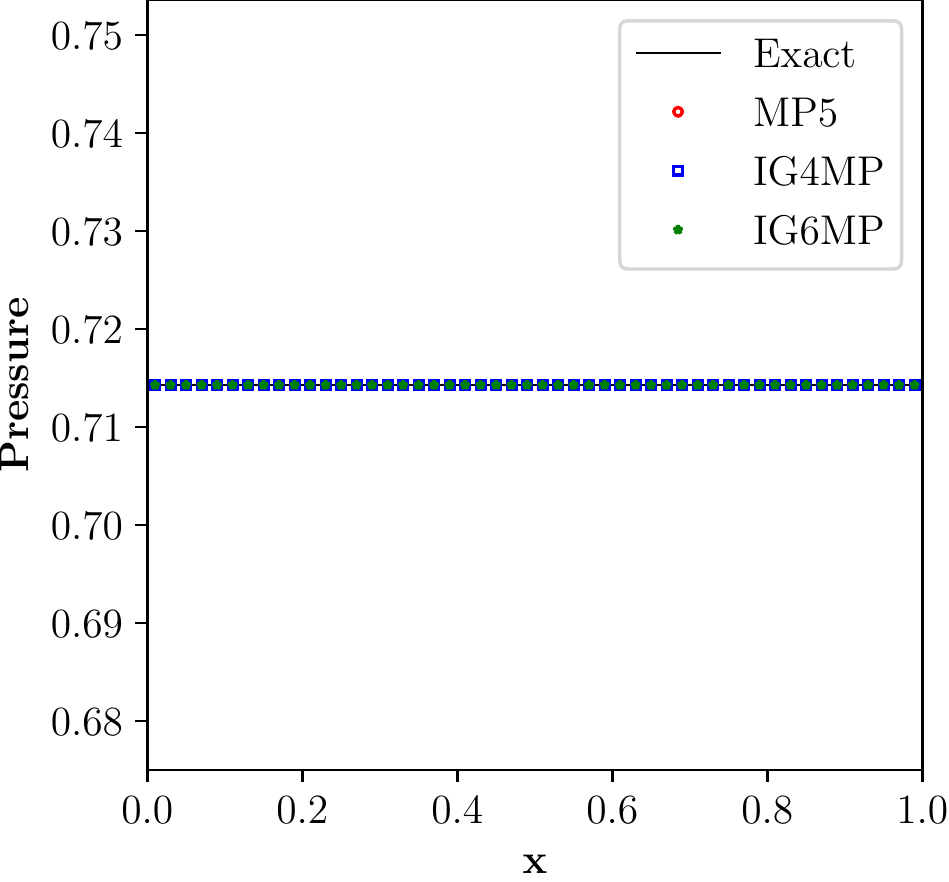}
\label{fig:multi_isolated-pres}}
\caption{Numerical solution for isolated material interface advection problem in Example \ref{ex:isolated}  on a grid size of $N=50$. Solid line: reference solution; green stars: IG6MP; blue squares: IG4MP; red circles: MP5.}
\label{fig_multi_isolated}
\end{figure}


\begin{example}\label{ex:curtain}{Inviscid shock-curtain interaction}
\end{example}
The interaction between a shockwave and a material interface is a challenging problem for numerical schemes.  The one-dimensional inviscid shock-curtain interaction problem introduced by Abgrall \cite{abgrall1996prevent} is considered in this final one-dimensional test case to study such behaviour. The initial conditions are as follow:
\begin{equation}
\left(\alpha_{1} \rho_{1}, \alpha_{2} \rho_{2}, u, p, \alpha_{1}, \gamma \right)=\left\{\begin{array}{ll}
\left(1.3765 , 0.0,0.3948,1.57, 1, 1.4\right), & 0 \leq x<0.25 \\
\left(1.0 , 0.0,0.0,1.0, 1, 1.4\right), & 0.25 \leq x<0.4 \text { or } 0.6 \leq x<1 \\
\left(0.0 , 0.138,0,1.0, 0, 1.67\right), & 0.4 \leq x<0.6
\end{array}\right.
\end{equation}
The shock wave travels in the air and moves to the right to interact with a helium bubble/curtain in the region $0.4 \leq x \leq 0.6$. Simulations are carried out on a grid size of 200 cell until $t$ = 0.3. Figs. \ref{fig:multi_curtain-den} and \ref{fig:multi_curtain-vel} show the density and velocity profiles obtained for all the schemes. IGMP schemes captured the discontinuities without any overshoots or undershoots in the solutions. The differences between IGMP schemes and MP scheme is negligible for all the one-dimensional test cases, and BVD algorithm effectively captures both shocks and material interfaces.

\begin{figure}[H]
\centering
\subfigure[Density]{\includegraphics[width=0.4\textwidth]{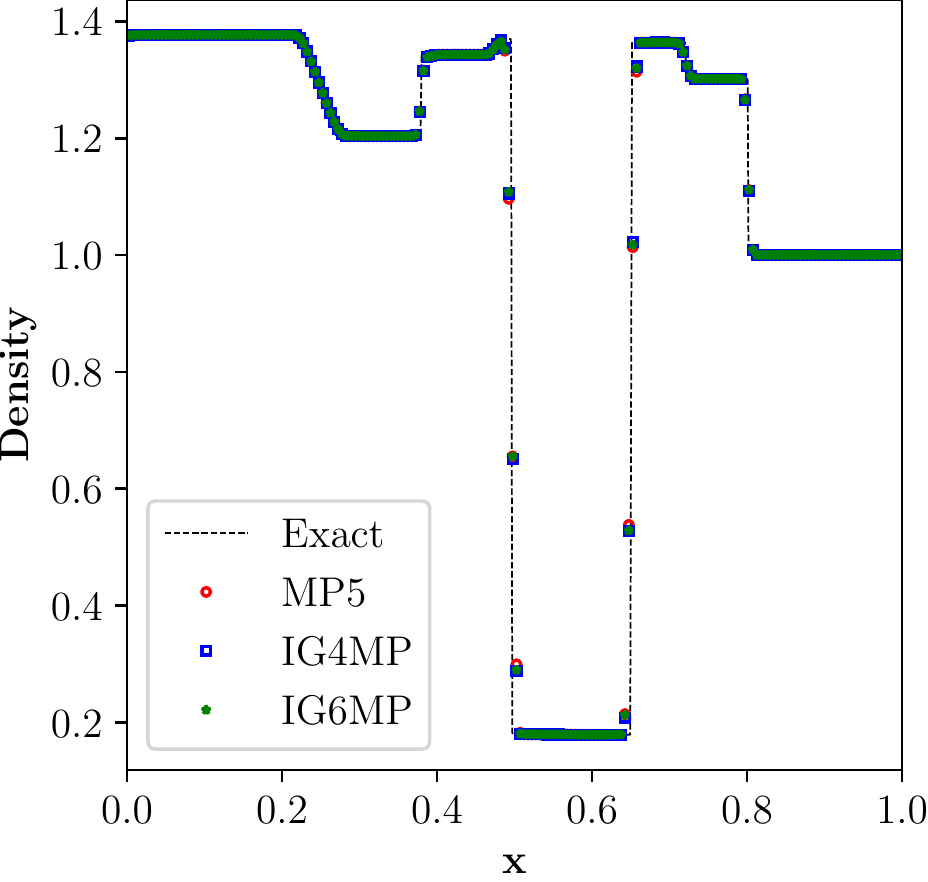}
\label{fig:multi_curtain-den}}
\subfigure[Velocity]{\includegraphics[width=0.4\textwidth]{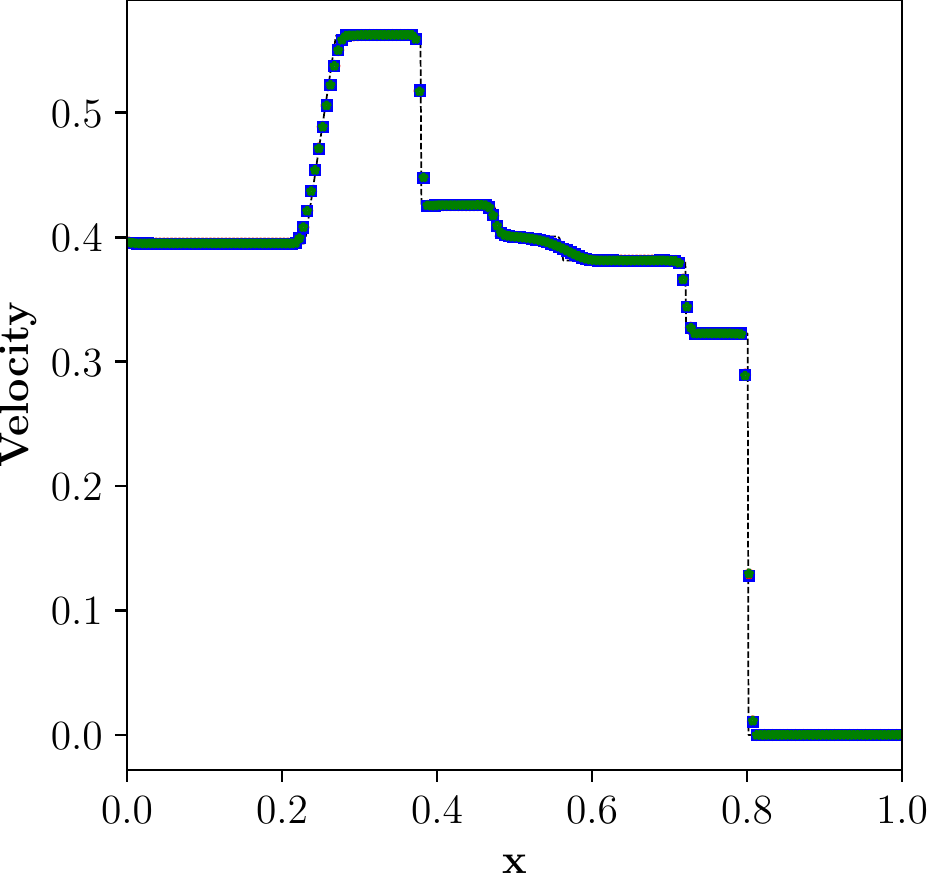}
\label{fig:multi_curtain-vel}}
\caption{Numerical solution for multi-species shock-curtain interaction problem in Example \ref{ex:curtain}  on a grid size of $N=100$. Dashed line: reference solution; green stars: IG6MP; blue squares: IG4MP; red circles: MP5.}
\label{fig_multi_curtain}
\end{figure}

\begin{example}\label{ex:He-bubble}{Two dimensional shock-cylinder interaction}
\end{example}

In this two-dimensional test case, the shock-cylinder interaction where a Helium bubble with density lighter than the air surrounding it interacts with a planar shock wave. This problem has been studied by Hass and Sturtevant  \cite{haas1987interaction} experimentally, and many researchers carried out computational studies for the same \cite{coralic2014finite, kawai2011high, Nonomura2012, Wong2017}. After the impact of the shock wave on the helium cylinder, the small scale vortices generated due to the baroclinic effects along the He-air interface can demonstrate the capabilities of the numerical scheme. Initial conditions of the test case are similar to that of Refs. \cite{kawai2011high, Wong2017} and are given by:

\begin{equation}
(\rho, u, v, p, \gamma)=\left\{\begin{array}{ll}
(1.3764,-0.3336,0,0,1.5698 / 1.4,1.4), & \text { for post-shock air, } \\
(0.1819,0,0,1 / 1.4,1.648), & \text { for helium cylinder, }\\
(1,0,0,1 / 1.4,1.4), & \text { for pre-shock air }.
\end{array}\right.
\end{equation}

The computational domain extends from $0.0 \leq x \leq 6.5D$ and $0.0 \leq y \leq 1.78D$. Initially, the Helium bubble is located at $[3.5$D$, 0.89$D$]$, where $D$ denotes the bubble's diameter and is taken as $D$ = 1 and a left moving normal shock of M 1.22 is placed at $x$=4.5$D$. All simulations are carried out with a CFL of 0.1, and slip wall boundary condition is imposed on top and bottom boundaries. The right boundary condition is extrapolated from the inside, and the left boundary is imposed as outflow. Two different grid-resolutions, 1300 $\times$ 356, and 2600 $\times$ 712, are employed for the simulation and the results of the normalized density gradient magnitude $\phi = $exp$(|\nabla \rho|/|\nabla \rho|_{max} )$ obtained by various schemes are shown in Figs. \ref{fig_bubble_He_fine} -\ref{fig_bubble_He_finest}. The present results are in good agreement with the computations reported in Ref. \cite{wang2020consistent}. It can be seen from these figures that there are no noticeable spurious oscillations in both the IGMP schemes. With increased grid resolution, the interface has become thinner, and more secondary instabilities of Kelvin-Helmholtz can be seen in the IG4MP scheme in the later stages, $t$ = 6.9, of the simulations. 

\begin{figure}[H]
\begin{onehalfspacing}
\centering\offinterlineskip
\subfigure[t = 3.25]{\includegraphics[width=0.28\textwidth]{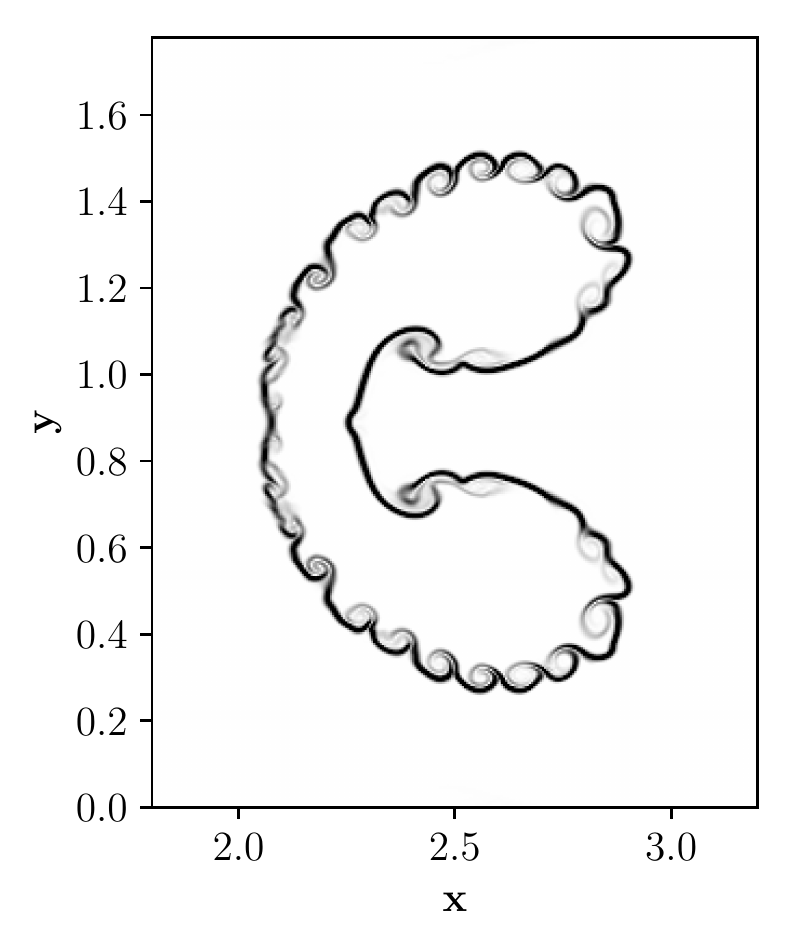}
\label{fig:MP5_325_fine}}
\subfigure[t = 4.90]{\includegraphics[width=0.28\textwidth]{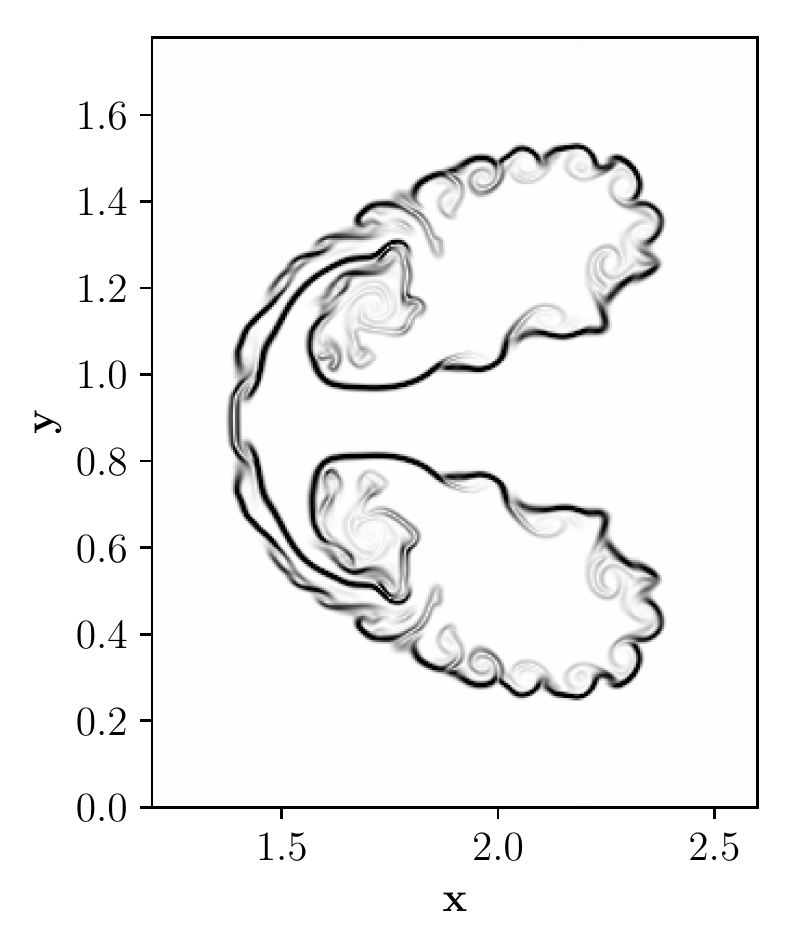}
\label{fig:MP5_490_fine}}
\subfigure[t = 6.90]{\includegraphics[width=0.28\textwidth]{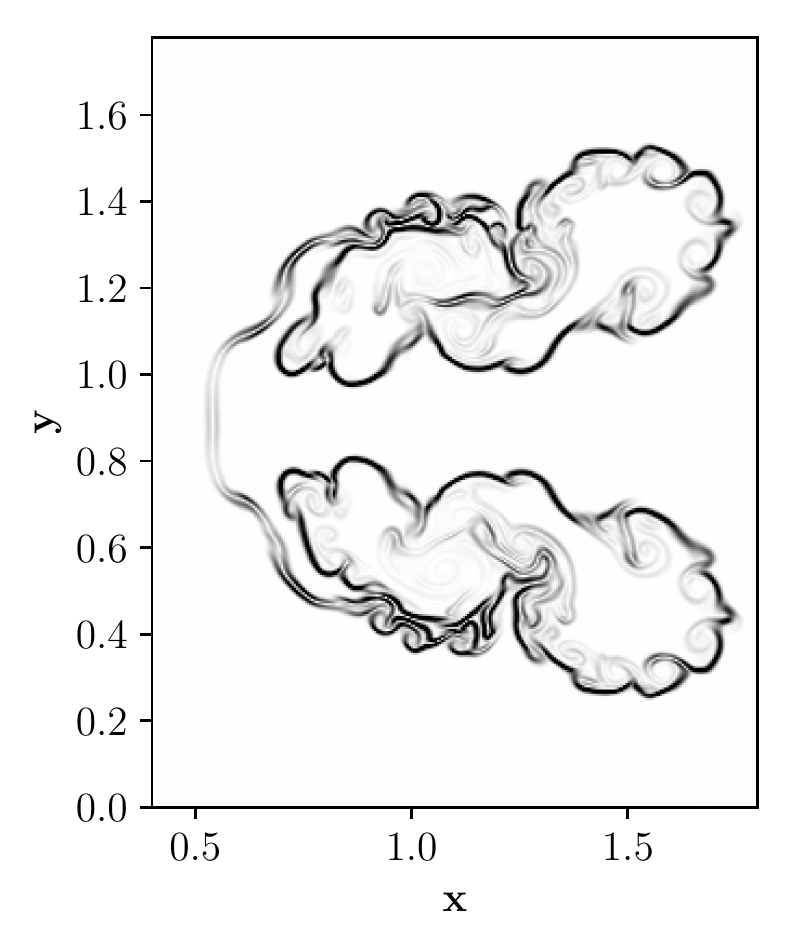}
\label{fig:MP5_695_fine}}
\subfigure[t = 3.25]{\includegraphics[width=0.28\textwidth]{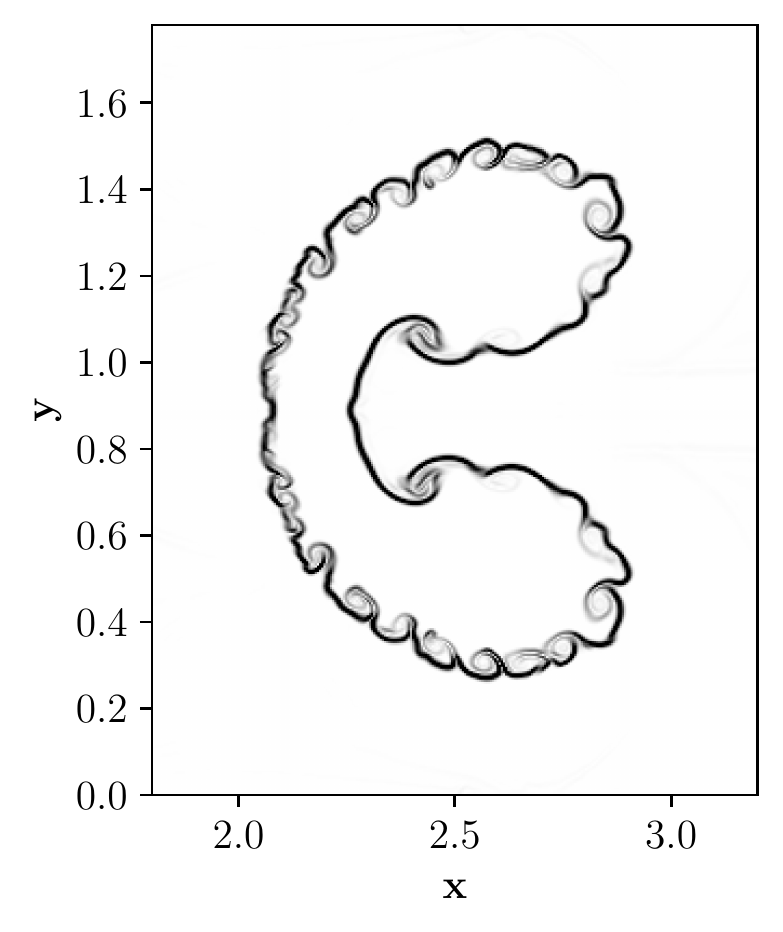}
\label{fig:IG6_325_fine}}
\subfigure[t = 4.90]{\includegraphics[width=0.28\textwidth]{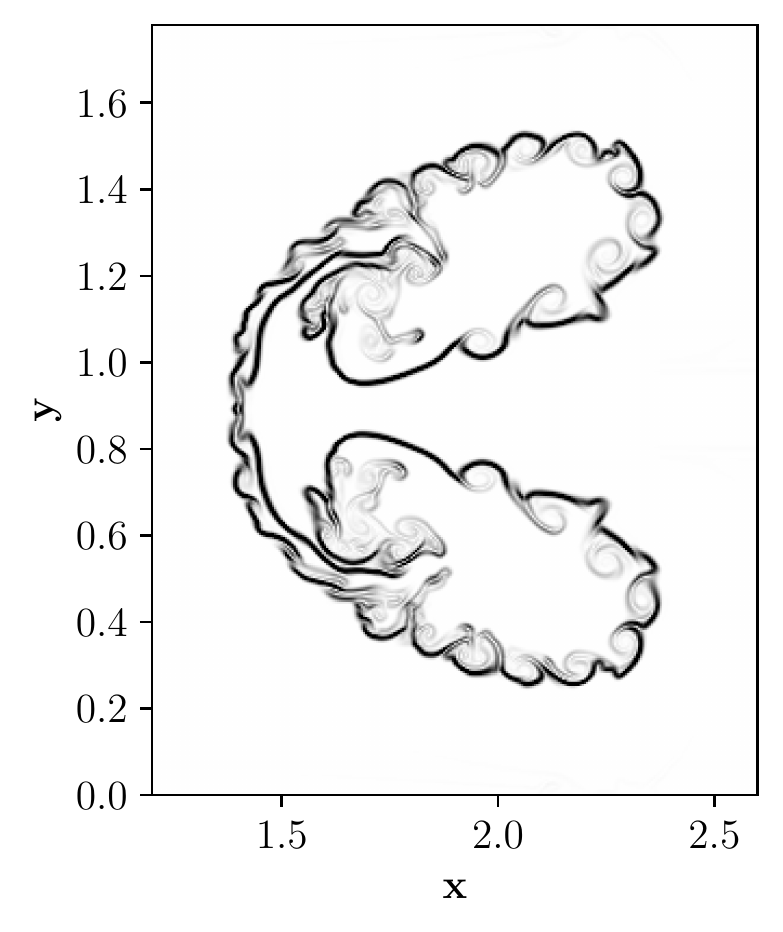}
\label{fig:IG6_490_fine}}
\subfigure[t = 6.90]{\includegraphics[width=0.28\textwidth]{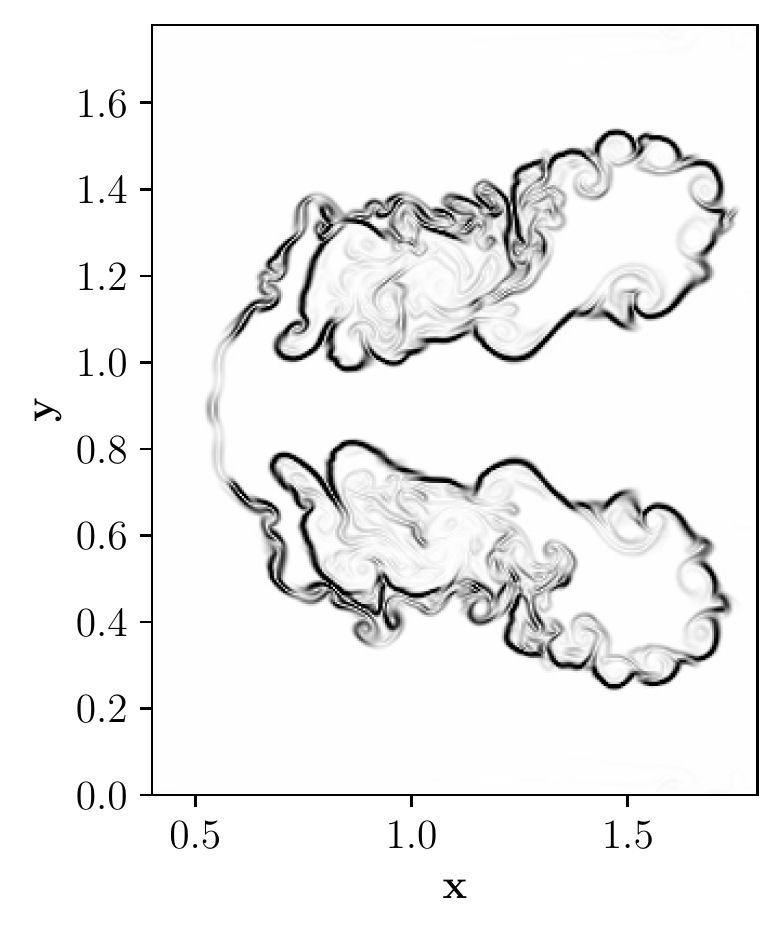}
\label{fig:IG6_695_fine}}
\subfigure[t = 3.25]{\includegraphics[width=0.28\textwidth]{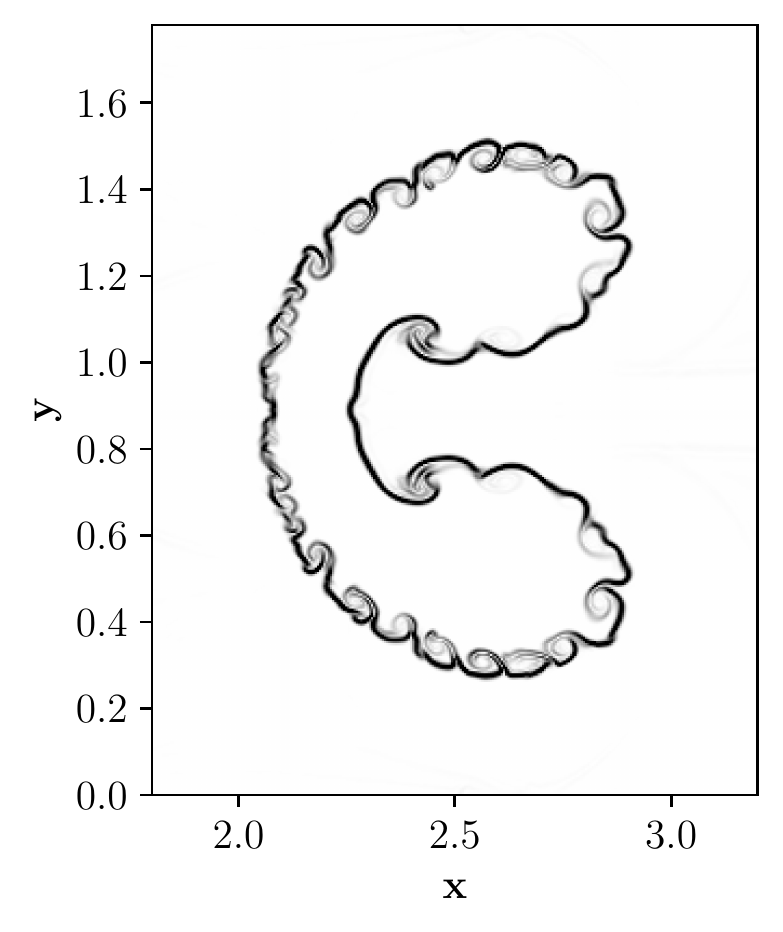}
\label{fig:IG4_325_fine}}
\subfigure[t = 4.90]{\includegraphics[width=0.28\textwidth]{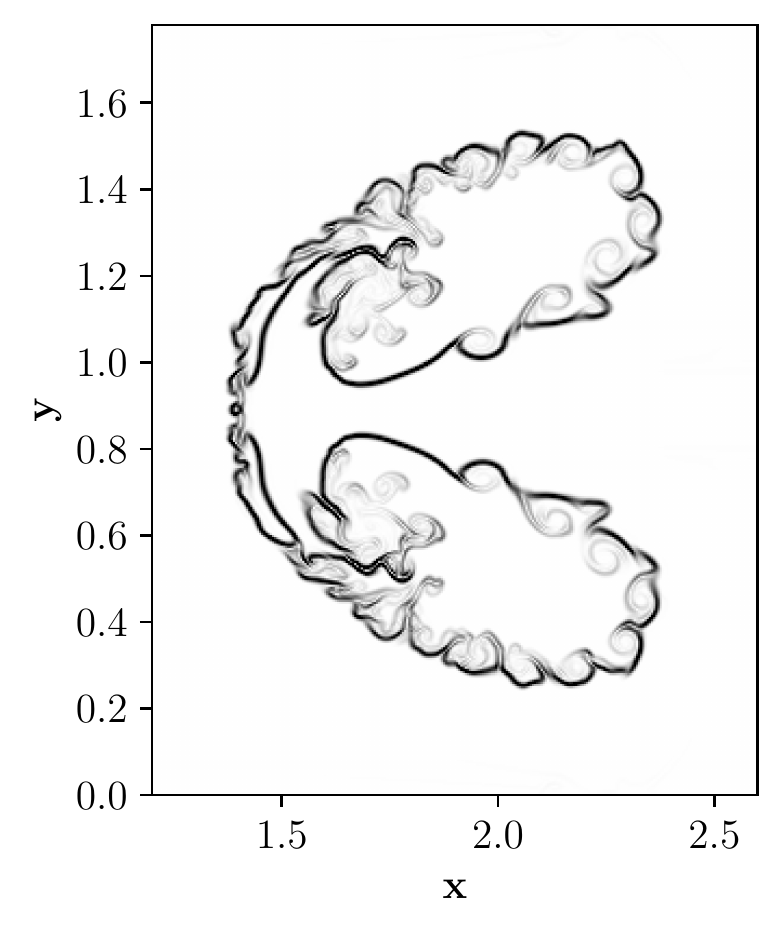}
\label{fig:IG4_490_fine}}
\subfigure[t = 6.90]{\includegraphics[width=0.28\textwidth]{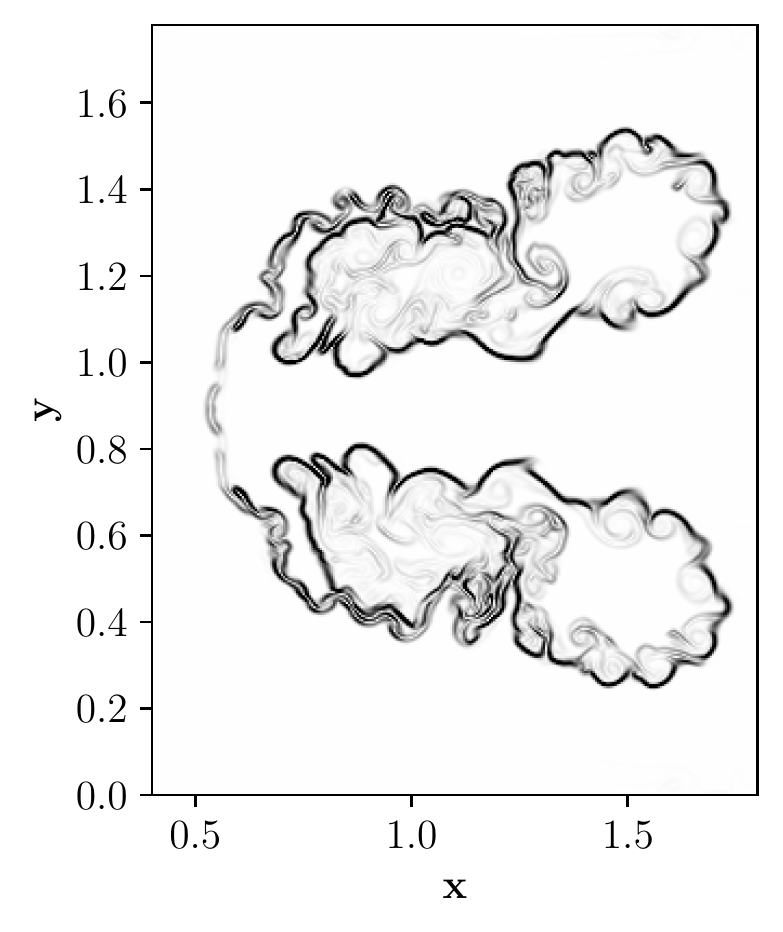}
\label{fig:IG4_695_fine}}
\caption{Comparison of normalized density gradient magnitude, $\phi$, contours for two-dimensional shock-bubble interaction problem in Example \ref{ex:He-bubble} on a grid resolution of 1300 $\times$ 356. Contours are from 1 to 1.7 at different times $t$ using different schemes. Top row: MP5; Middle row: IG6MP and bottom row: IG4MP.}
\label{fig_bubble_He_fine}
\end{onehalfspacing}
\end{figure}

\begin{figure}[H]
\begin{onehalfspacing}
\centering\offinterlineskip
\subfigure[t = 3.25]{\includegraphics[width=0.28\textwidth]{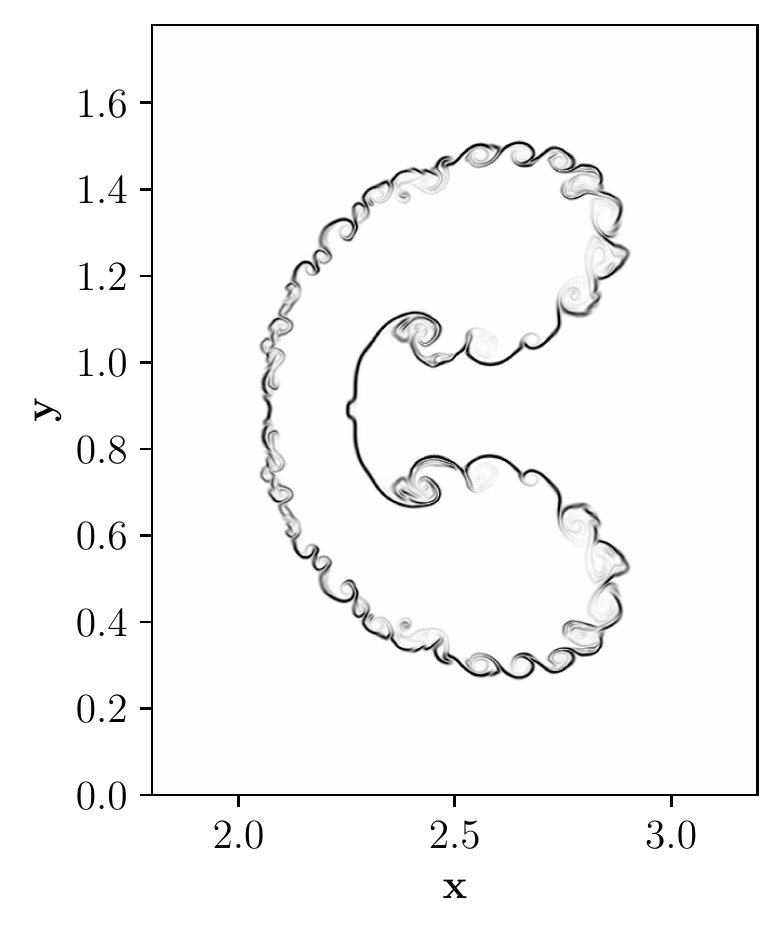}
\label{fig:MP5_325_finest}}
\subfigure[t = 4.90]{\includegraphics[width=0.28\textwidth]{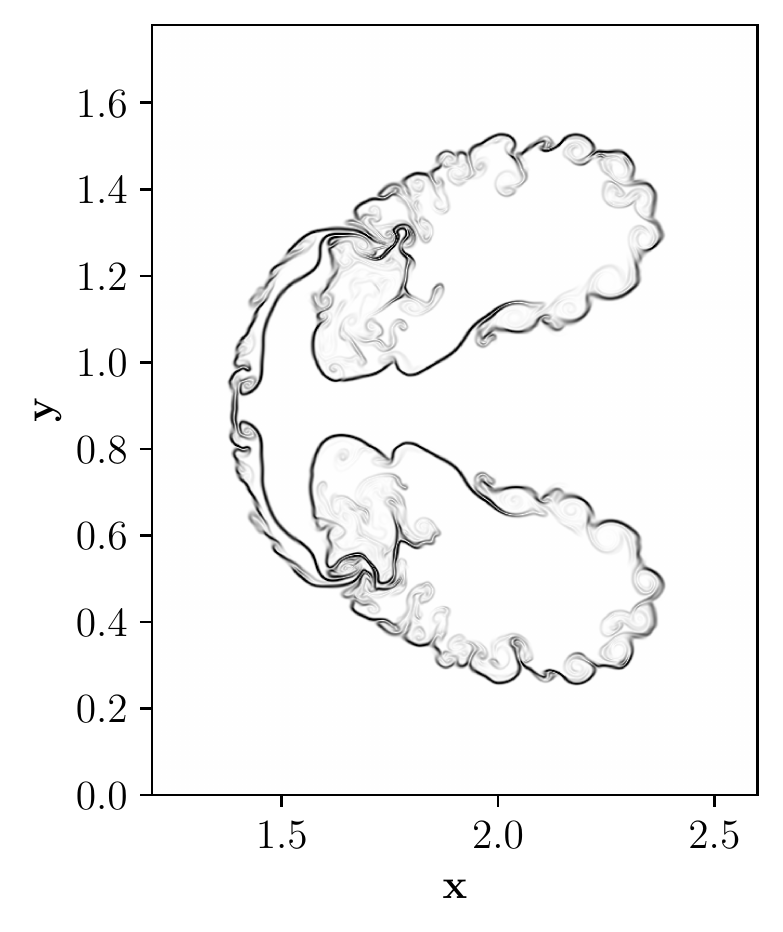}
\label{fig:MP5_490_finest}}
\subfigure[t = 6.90]{\includegraphics[width=0.28\textwidth]{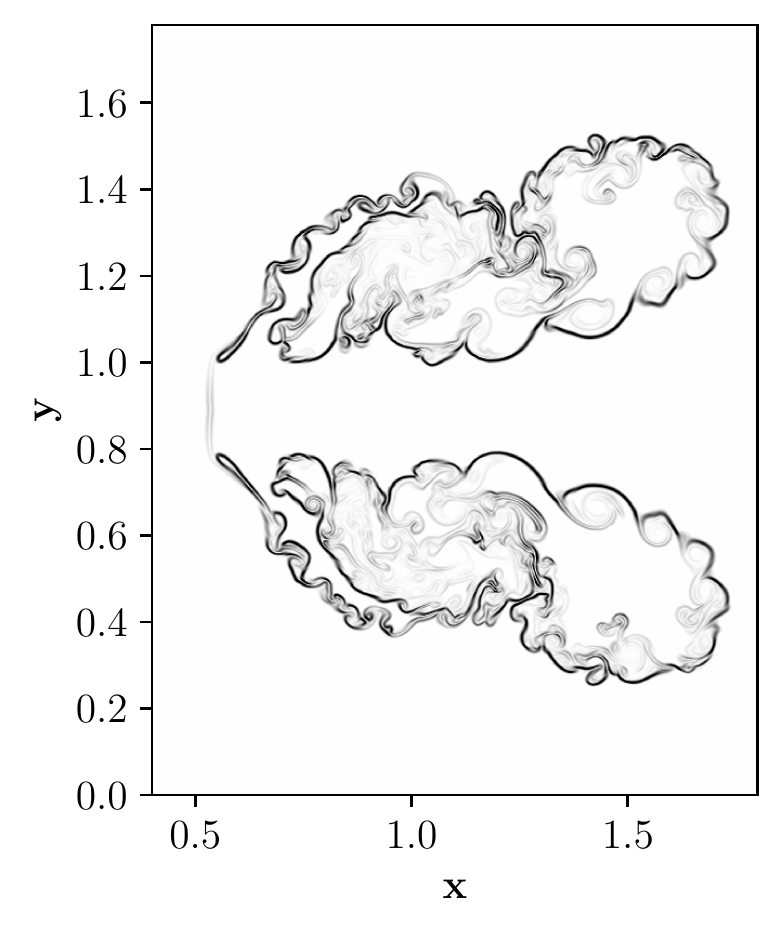}
\label{fig:MP5_695_finest}}
\subfigure[t = 3.25]{\includegraphics[width=0.28\textwidth]{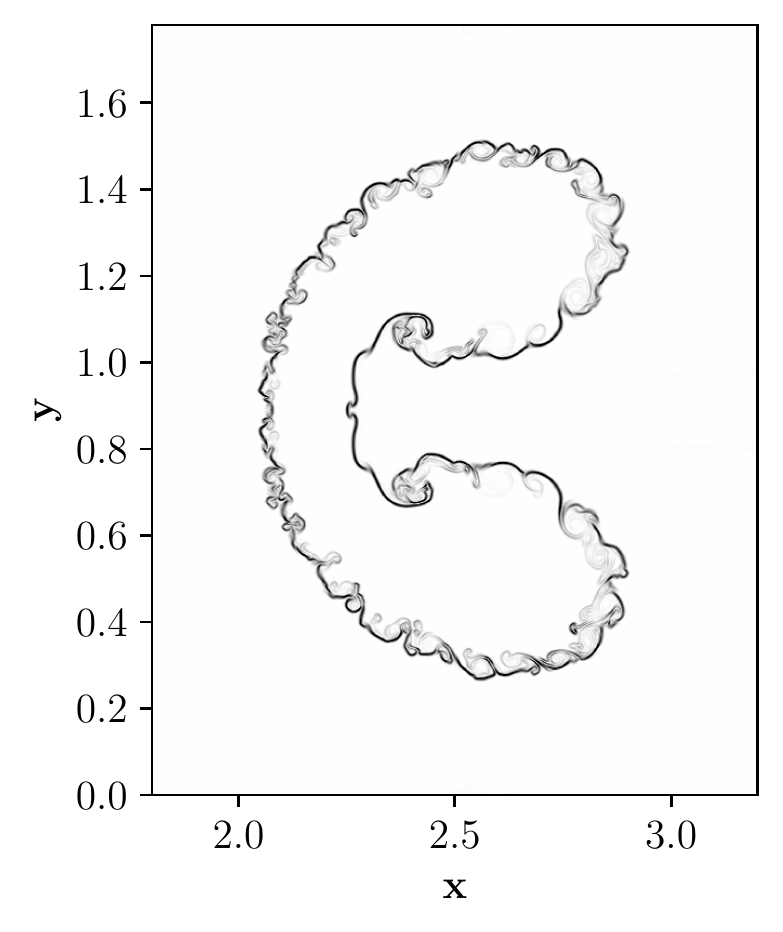}
\label{fig:IG6_325_finest}}
\subfigure[t = 4.90]{\includegraphics[width=0.28\textwidth]{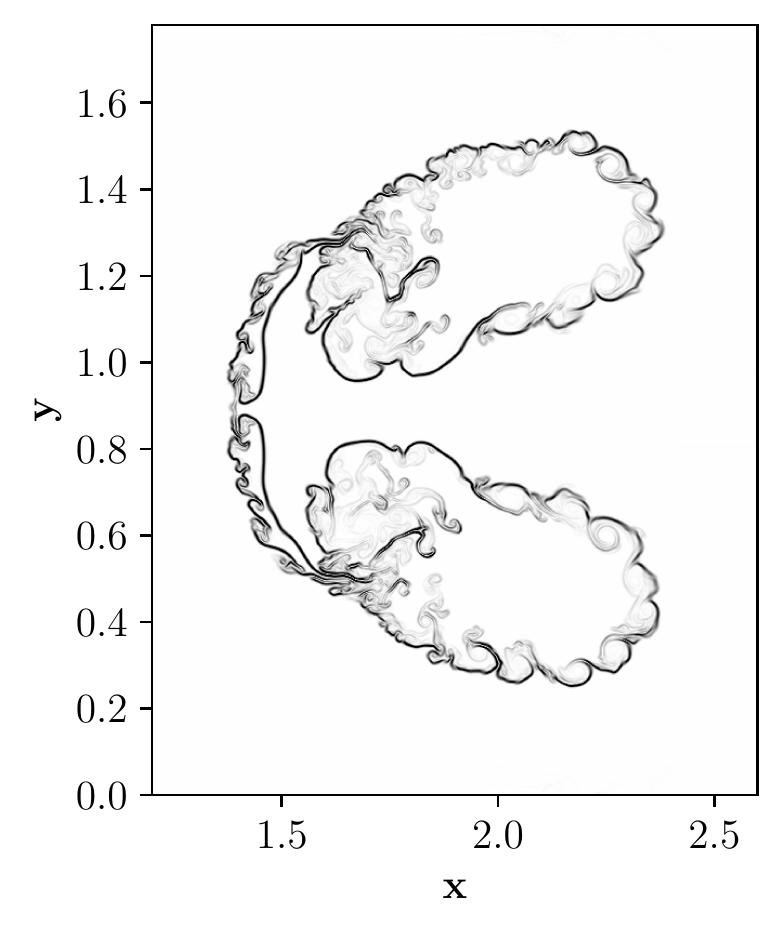}
\label{fig:IG6_490_finest}}
\subfigure[t = 6.90]{\includegraphics[width=0.28\textwidth]{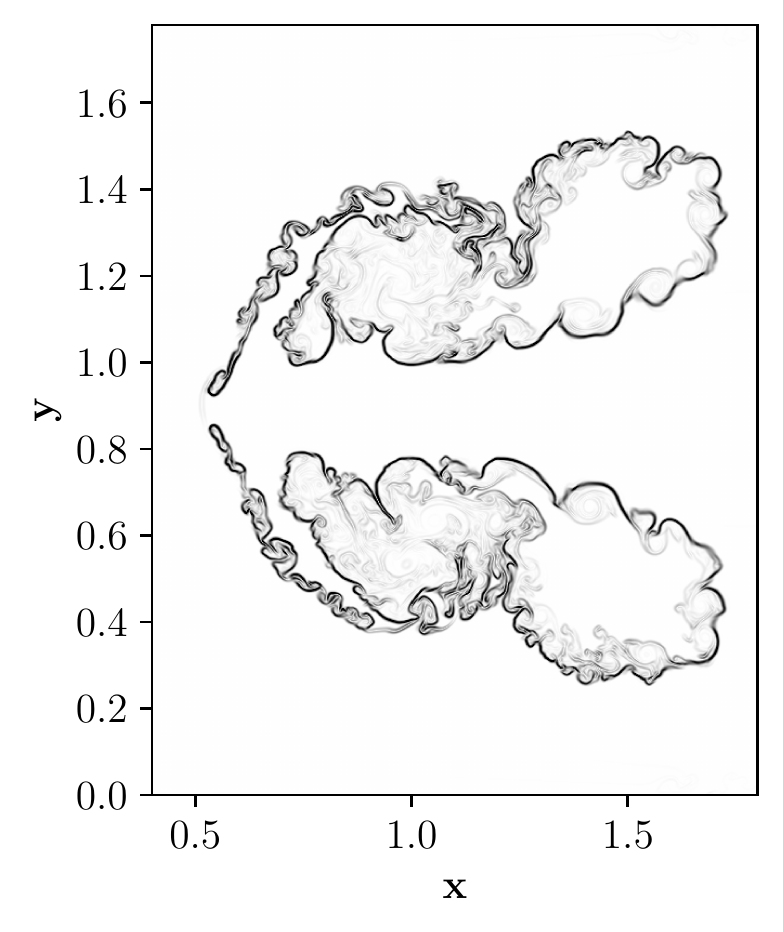}
\label{fig:IG6_695_finest}}
\subfigure[t = 3.25]{\includegraphics[width=0.28\textwidth]{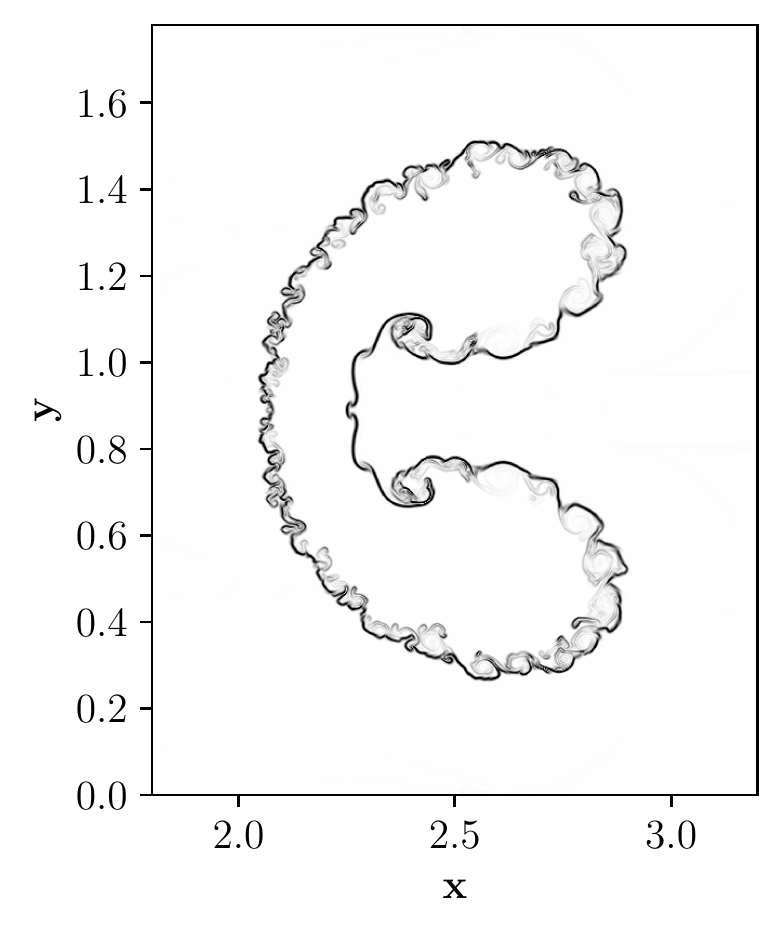}
\label{fig:IG4_325_finest}}
\subfigure[t = 4.90]{\includegraphics[width=0.28\textwidth]{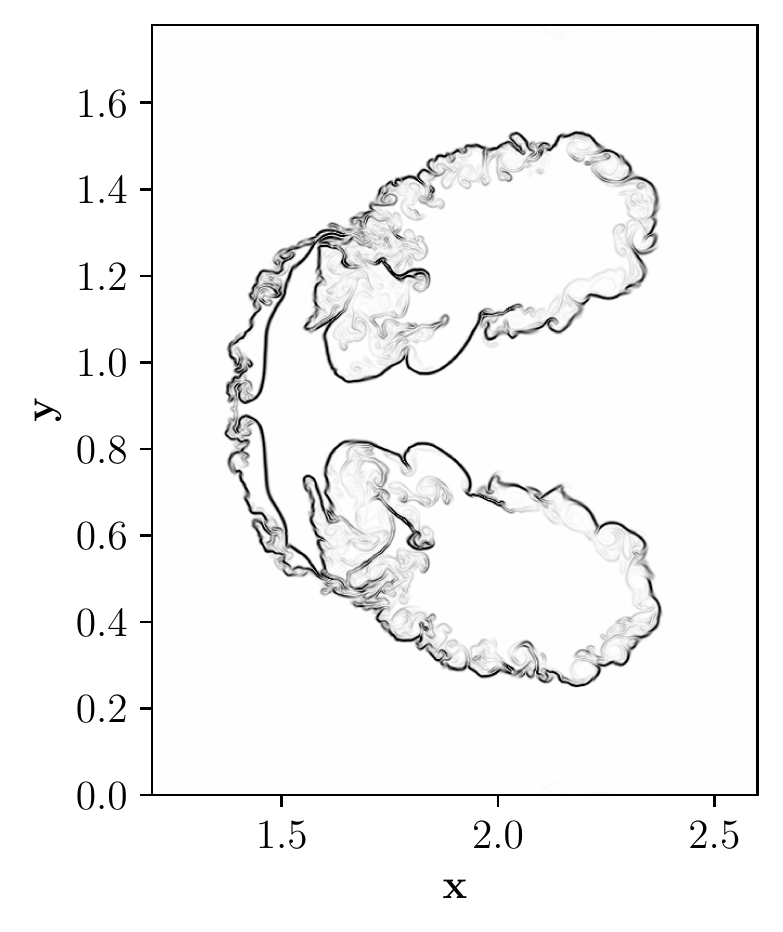}
\label{fig:IG4_490_finest}}
\subfigure[t = 6.90]{\includegraphics[width=0.28\textwidth]{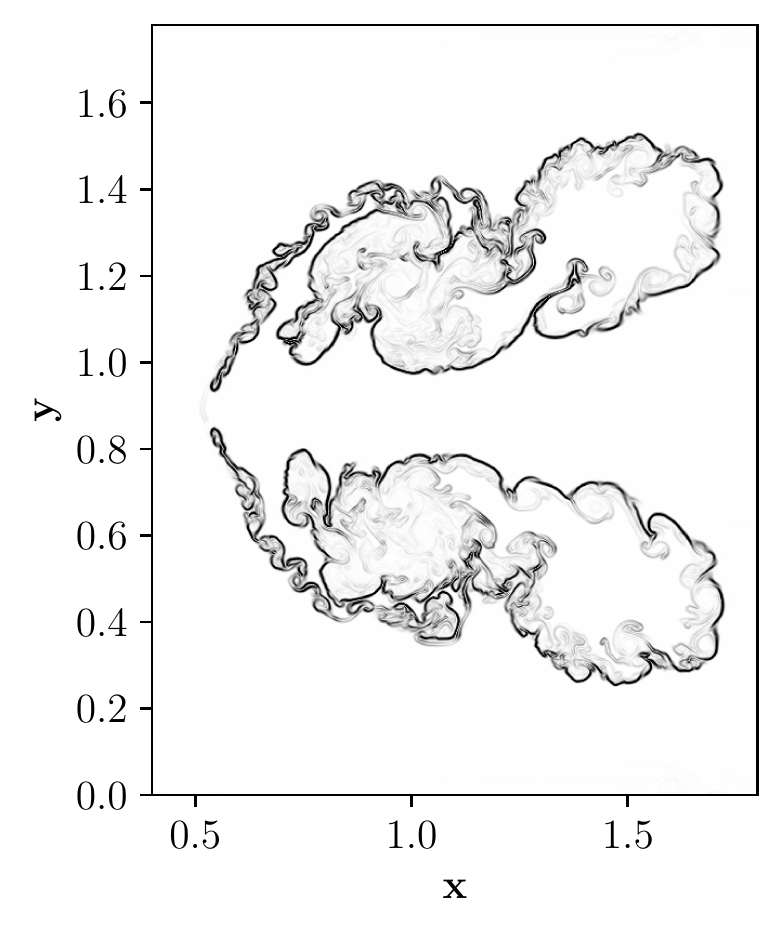}
\label{fig:IG4_695_finest}}
\caption{Comparison of normalized density gradient magnitude, $\phi$, contours for two-dimensional shock-bubble interaction problem in Example \ref{ex:He-bubble} on a grid resolution of 2600 $\times$ 712. Contours are from 1 to 1.7 at different times $t$ using different schemes. Top row: MP5; Middle row: IG6MP and bottom row: IG4MP.}
\label{fig_bubble_He_finest}
\end{onehalfspacing}
\end{figure}

\begin{example}\label{ex:RM-viscous}{Two dimensional viscous Richtmeyer-Meshkov instability}
\end{example}
In this last test case, the two-dimensional viscous Richtmeyer-Meshkov (RMI) instability is computed. Numerical simulations of inviscid test cases are carried out by Nonomura et al.\cite{Nonomura2012}, and Kawai and Terashima \cite{kawai2011high} whereas Yee and B. Sjögreen \cite{yee2007simulation} have included physical viscosity. RMI occurs when an incident shock accelerates an interface between two fluids of different densities. The schematic of the test case's initial condition is shown in Fig. \ref{rminstability}.  The computational domain for this test case extends from $0.0 \leq x \leq 16.0\lambda$ and $0.0 \leq y \leq 1.0\lambda$ where $\lambda$ is the initial perturbation wavelength and the initial shape of the interface is given by

\begin{equation}
\frac{x}{\lambda}=0.4-0.1 \sin \left(2 \pi\left(\frac{y}{\lambda}+0.25\right)\right),
\end{equation}
 and the initial conditions are as follows:

\begin{equation}
(\rho, u, v, p, \gamma)=\left\{\begin{array}{ll}
(1,1.24,0,1 / 1.4,1.4), & \text { for pre-shock air } \\
(1.4112,0.8787,0,1.6272 / 1.4,1.4), & \text { for post-shock air } \\
(5.04,1.24,0,1 / 1.4,1.093), & \text { for } S F_{6}
\end{array}\right.
\end{equation}

\begin{figure}[H]
\centering
 \includegraphics[width=1.0\textwidth]{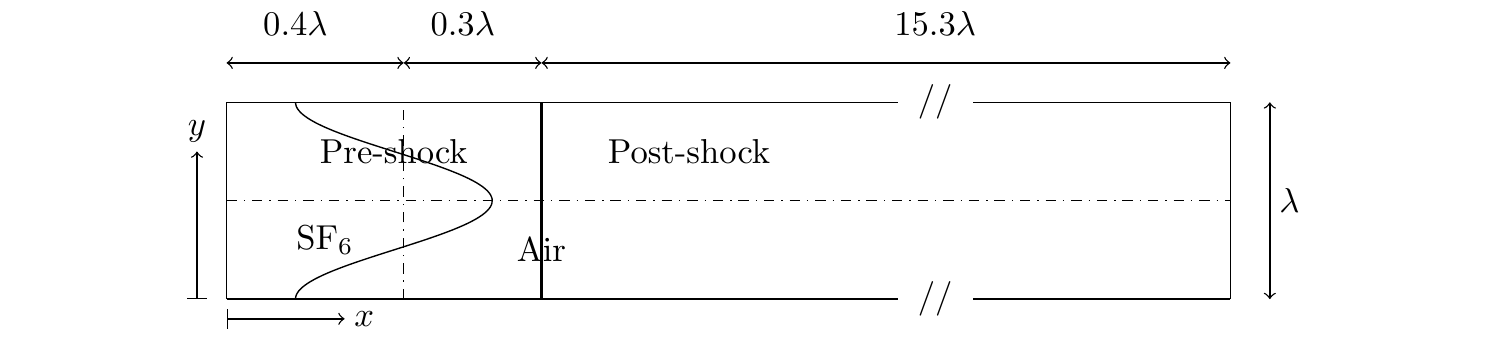}
\caption{Schematic of initial condition of Ricthmyer-Meshkov instability, Example \ref{ex:RM-viscous}.}
\label{rminstability}
\end{figure}

where, $\lambda$ =1. In the current simulations we considered  constant dynamic viscosity $\mu$ = 1.0 $\times$ $10^{-4}$ (resulting in $Re= 10^4$) \cite{yee2007simulation}. Periodic boundary conditions are imposed at the top, and bottom boundaries of the domain and the initial values are fixed at the left and right boundaries. Due to the cartesian grid considered in the present simulation, secondary instabilities will be generated at the material interface that are not observed in experimental conditions \cite{jacobs2005experiments}. To avoid these instabilities, the initial perturbation is smoothened by adding an artificial diffusion layer given by Ref. \cite{Wong2017}:

\begin{equation}
\begin{array}{c}
f_{sm}=\frac{1}{2}\left(1+\operatorname{erf}\left(\frac{\Delta D}{E_{i} \sqrt{\Delta x \Delta y}}\right)\right) \\
u=u_{L}\left(1-f_{sm}\right)+u_{R} f_{sm}
\end{array}
\end{equation}
where $u$ represents the primitive variables near the initial interface. The parameter $E_i$ introduces additional thickness to the initial material interface, $\Delta D$ is the distance from the initial perturbed material interface and subscripts $L$ and $R$ denote the left and right interface conditions. Parameter $E_i$ is chosen as 5 in this test case. Simulations are conducted on two different grid sizes, 2048 $\times$128 and 4096 $\times$ 256 cells, with a constant CFL of 0.1. Computational results of normalized density gradient magnitude $\phi = $exp$(|\nabla \rho|/|\nabla \rho|_{max} )$ obtained at $t$ = 11.0 by various schemes are shown in Fig. \ref{fig_RM_viscous}. It can be seen from these figures that there are no noticeable spurious oscillations in both the IGMP schemes. The IG4MP scheme is clearly less dissipative and therefore produces more small scale structures along the material interface than IG6MP and MP schemes.

\begin{figure}[H]
\begin{onehalfspacing}
\centering\offinterlineskip
\subfigure[MP5]{\includegraphics[width=0.3\textwidth]{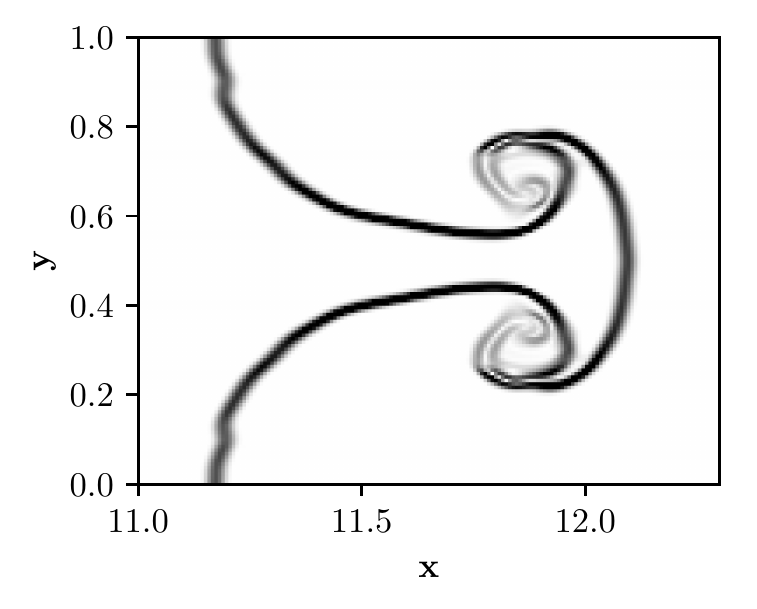}
\label{fig:MP5_rmv}}
\subfigure[IG6MP]{\includegraphics[width=0.3\textwidth]{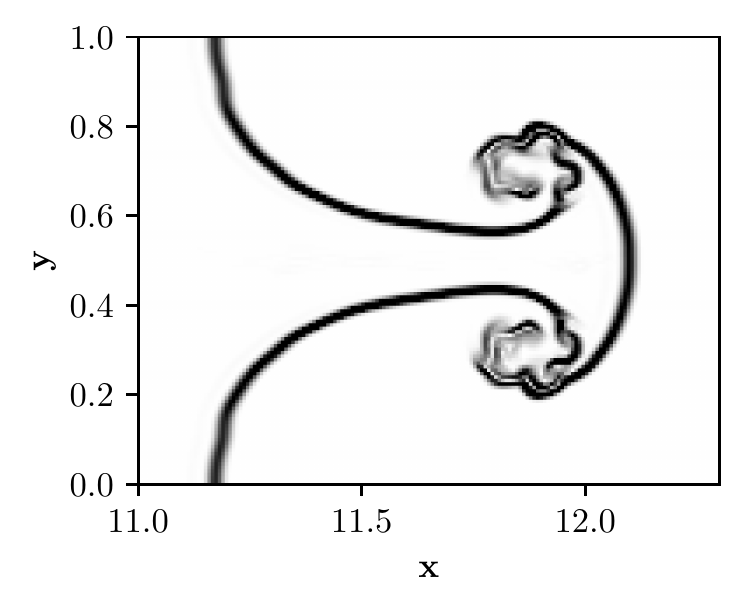}
\label{fig:IG6MP_rmv}}
\subfigure[IG4MP]{\includegraphics[width=0.3\textwidth]{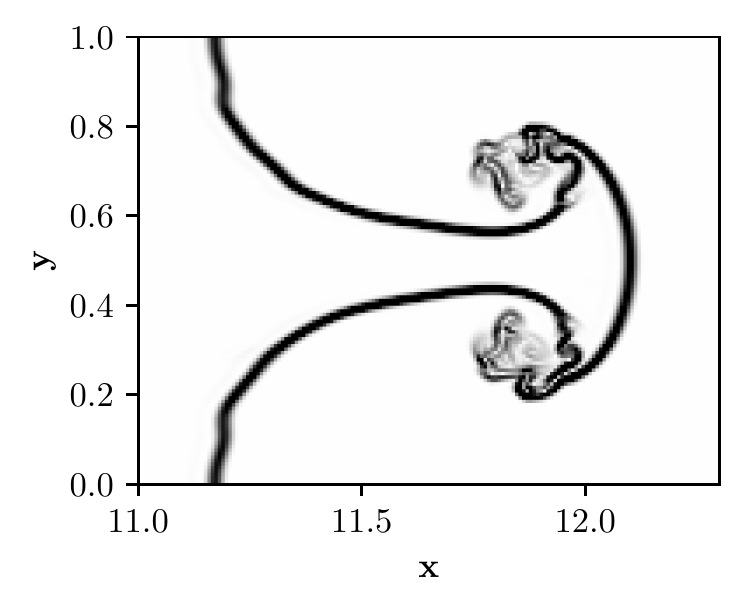}
\label{fig:IG4MP_rmv}}
\subfigure[MP5]{\includegraphics[width=0.3\textwidth]{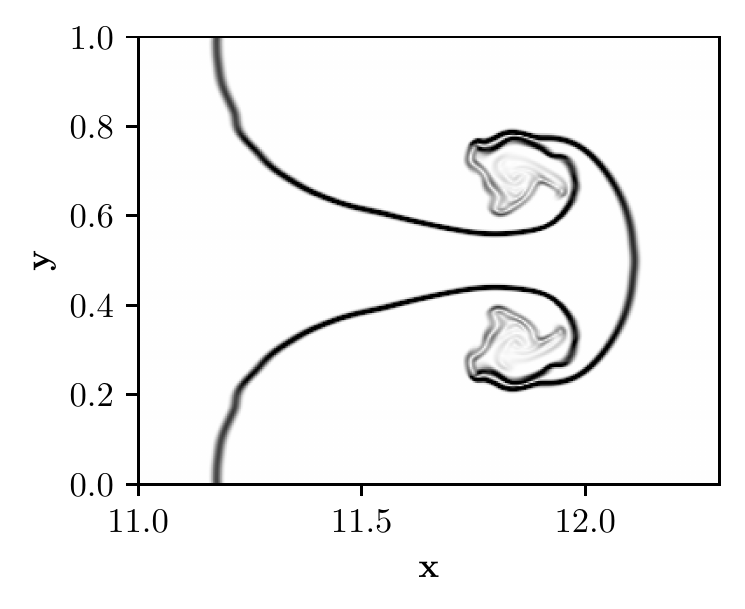}
\label{fig:MP5_rmvf}}
\subfigure[IG6MP]{\includegraphics[width=0.3\textwidth]{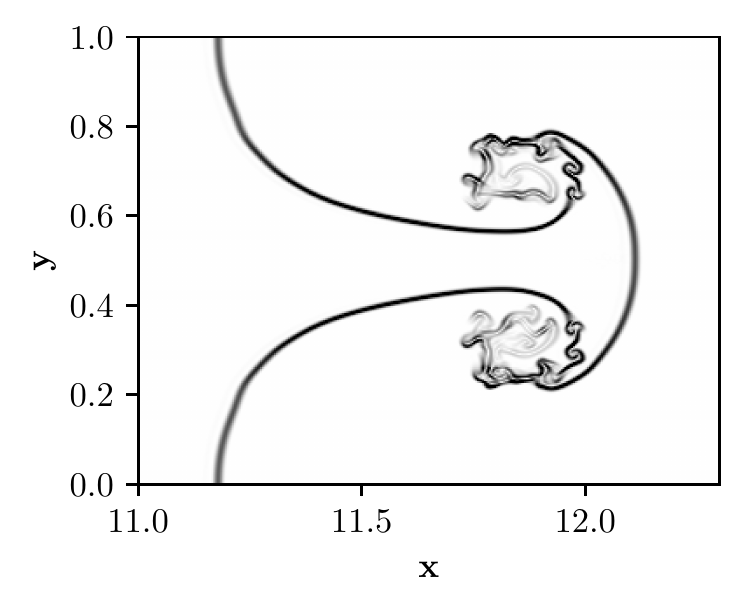}
\label{fig:IG6MP_rmvf}}
\subfigure[IG4MP]{\includegraphics[width=0.3\textwidth]{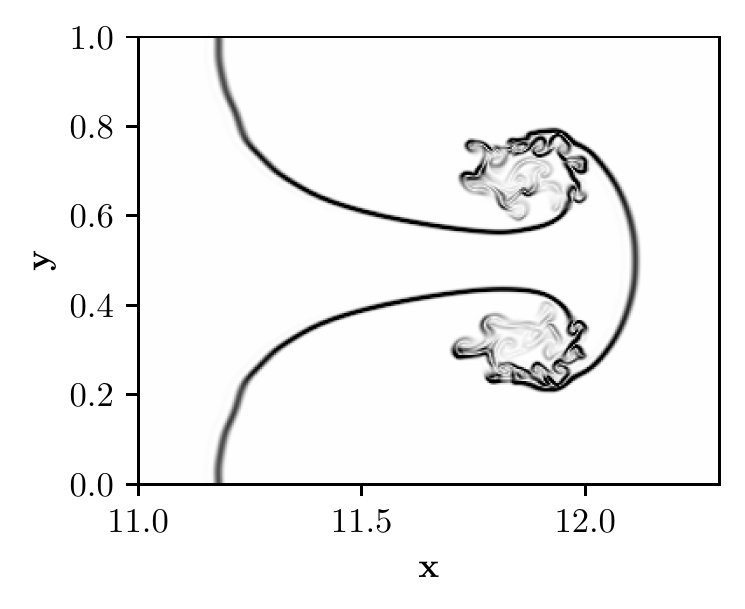}
\label{fig:IG4MP_rmvf}}
\caption{Comparison of normalized density gradient magnitude, $\phi$, contours for two-dimensional shock-bubble interaction problem in Example \ref{ex:RM-viscous} on a grid resolution top row: 2048 $\times$ 128 and bottom row: 4096 $\times$ 256. Contours are from 1 to 1.7 at time t$=$11.0 s using different schemes.}
\label{fig_RM_viscous}
\end{onehalfspacing}
\end{figure}
\section{Conclusions}\label{sec-5}

In this paper, high-order implicit gradient finite volume schemes have been developed and the important contributions and observations of the paper are summarized as follows

\begin{enumerate}
\item We proposed a novel approach of computing the cell interface values in finite volume framework where the gradients of reconstruction polynomials are computed by compact finite differences. The resulting upwind schemes have superior dissipation and dispersion than the compact reconstruction approach. 
\item These gradients are re-used in viscous flux computation and post-processing, which improved efficiency and indicated that significant improvements could be obtained in flow features for the viscous flows.
\item Problem independent shock-capturing technique via the BVD algorithm gave superior results for several benchmark test cases involving shocks, material interfaces and small scale features.
\item Results obtained by IG4MP are superior to IG6MP scheme, even though they they share the same stencil, due its superior dispersion and dissipation properties. 
\end{enumerate}

Implicit gradients are expensive to compute and may be replaced with explicit higher-order ones. It is possible to replace the second derivatives with explicit gradients, which may improve the scheme's overall computational efficiency. Despite the BVD algorithm's excellent performance in capturing discontinuities, a native shock-capturing method can be developed to reduce computational expenses further. Future work will focus on such capabilities.

\section*{Acknowledgements}
A.S. is supported by Technion Fellowship during this work. A.S. would like to give special thanks to Jyothi and Arya for the motivation to complete this paper.

\section*{Appendix}
\renewcommand{\thesubsection}{\Alph{subsection}}

\subsection{\textcolor{black}{Eigenvalues and eigenvectors}} \label{sec-appa}
\onehalfspacing
 The left and right eigenvectors of the two-dimensional single-component equations, denoted by $\bm{L_n}$ and $\bm{R_n}$, used for characteristic variable projection are as follows:
\begin{align}\label{eqn:leftright}
	\bm{R_n} = \begin{bmatrix} 
   		 1 & 1 & 1 & 0 \\
		 \\
   		 -\frac{n_xc}{\rho} & 0 & \frac{n_xc}{\rho} & \frac{l_x}{\rho} \\
		 \\
   		 -\frac{n_yc}{\rho} & 0 & \frac{n_yc}{\rho} & \frac{l_y}{\rho} \\
		 \\
     c^2 & 0 & c^2 & 0 \\
   			 \end{bmatrix}, \enskip
	\bm{L_n} = \begin{bmatrix}
   		 0 & -\frac{n_x\rho}{2c} &-\frac{n_y\rho}{2c} & \frac{1}{2c^2} \\
		 \\
   		 1 & 0 & 0 & -\frac{1}{c^2} \\
		 \\
   		 0 & -\frac{n_x\rho}{2c} & \frac{n_y\rho}{2c} & \frac{1}{2c^2} \\
		 \\
     0 & \rho l_x & \rho l_y & 0 \\
   			 \end{bmatrix}, &&\\\nonumber
\end{align}
where $\bm{n}$ = $[n_x \ n_y]^t$ and $[l_x \ l_y]^t$ is a tangent vector (perpendicular to $\bm{n}$) such as $[l_x \ l_y]^t$ = $[-n_y \ n_x]^t$. By taking $\bm{n}$ = $[1, 0]^t$ and $[0, 1]^t$ we obtain the corresponding eigenvectors in $x$ and $y$ directions. Similarly, left and right eigenvectors of the two-dimensional multi-component equations, denoted by $\bm{L_n}$ and $\bm{R_n}$, used for characteristic variable projection are as follows:

\begin{align}\label{eqn:leftright-multi}
\onehalfspacing
	\bm{R_n} = \begin{bmatrix} 
   		\frac{\alpha_{1} \rho_{1}}{c^{2} \rho} & 1 & 0 & 0 & 0 & \frac{\alpha_{1} \rho_{1}}{c^{2} \rho} \\ 
\\
\frac{\alpha_{2} \rho_{2}}{c^{2} \rho} & 0 & 1 & 0 & 0 & \frac{\alpha_{2} \rho_{2}}{c^{2} \rho}\\
\\
- \frac{n_x}{c \rho} & 0 & 0 & n_y & 0 & \frac{n_x}{c \rho}\\
\\
- \frac{n_y}{c \rho} & 0 & 0 & n_x & 0 & \frac{n_y}{c \rho}\\
\\
1 & 0 & 0 & 0 & 0 & 1\\
\\
0 & 0 & 0 & 0 & 1 & 0
   			 \end{bmatrix}, \enskip
	\bm{L_n} = \begin{bmatrix}
   		0 & 0 & - \frac{n_x c \rho}{2} & - \frac{n_y c \rho}{2} & \frac{1}{2} & 0\\
\\
1 & 0 & 0 & 0 & - \frac{\alpha_{1} \rho_{1}}{c^{2} \rho} & 0\\
\\
0 & 1 & 0 & 0 & - \frac{\alpha_{2} \rho_{2}}{c^{2} \rho} & 0\\
\\
0 & 0 & n_y & n_x & 0 & 0\\
\\
0 & 0 & 0 & 0 & 0 & 1\\
\\
0 & 0 & \frac{n_x c \rho}{2} &  \frac{n_y c \rho}{2} & \frac{1}{2} & 0
   			 \end{bmatrix}. &&\\\nonumber
\end{align}


\bibliographystyle{elsarticle-num}
\bibliography{bvd_ref}
\end{document}